\documentclass[amscd,amssymb,epic,eepic,reqno,12pt]{amsart}
\usepackage{latexsym,amsmath,amsfonts,amscd,amssymb}
\renewcommand\today{July 9, 2003}

\theoremstyle{plain}  % default
\newtheorem{theorem}{Theorem}[section]
\newtheorem*{theoremA}{Theorem A}
\newtheorem*{theoremB}{Theorem B}

\newtheorem*{theorem*}{Theorem}

\newtheorem{corollary}[theorem]{Corollary}
\newtheorem{lemma}[theorem]{Lemma}
\newtheorem{proposition}[theorem]{Proposition}

\newtheorem{tech-lemma}[theorem]{Technical Lemma}
\theoremstyle{definition}
\newtheorem{definition}[theorem]{Definition}
\theoremstyle{remark}

\newtheorem{remark}[theorem]{Remark}

\newtheorem*{claim*}{Claim}
\numberwithin{equation}{section}
\renewcommand{\leq}{\leqslant}
\renewcommand{\le}{\leqslant}
\renewcommand{\geq}{\geqslant}
\renewcommand{\ge}{\geqslant}
\renewcommand{\setminus}{\smallsetminus}
\newcommand{\lto}{\longrightarrow}
\newcommand{\R}{\mathbb{R}}
\newcommand{\Z}{\mathbb{Z}}
\newcommand{\C}{\mathbb{C}}
\newcommand{\HH}{\mathbb{H}}
\newcommand{\dbar}{\bar{\partial}}
\newcommand{\into}{\hookrightarrow}
\newcommand{\abs}[1]{\lvert#1\rvert}
\newcommand{\norm}[1]{\lVert#1\rVert}
\newcommand{\lie}{\mathfrak}

\newcommand{\suchthat}{\;\;|\;\;}

\newcommand{\PU}{\mathrm{PU}}
\newcommand{\PGL}{\mathrm{PGL}} 
 
\newcommand{\SU}{\mathrm{SU}}
\newcommand{\U}{\mathrm{U}}
\newcommand{\GL}{\mathrm{GL}}
\newcommand{\GCD}{\mathrm{GCD}}

\DeclareMathOperator{\ad}{ad}
\DeclareMathOperator{\Ad}{Ad}

\DeclareMathOperator{\tr}{tr}
\DeclareMathOperator{\rk}{rk}
\DeclareMathOperator{\im}{im}

\DeclareMathOperator{\Hom}{Hom}
\DeclareMathOperator{\End}{End}

\DeclareMathOperator{\Id}{Id}

\begin{document}
\pagestyle{empty}

\noindent 
{\Large\textbf{Surface group representations \\ and $\U(p,q)$-Higgs 
    bundles}} \bigskip

\noindent 
{\textbf{Steven B. Bradlow}}
\footnotemark[1]$^{,}$\footnotemark[2] \\
%\textit{\nocorr
Department of Mathematics, \\
University of Illinois, \\
Urbana, IL 61801,
USA \\
E-mail: \texttt{bradlow@math.uiuc.edu}
%}
\medskip

\noindent 
\textbf{Oscar Garc{\'\i}a--Prada}
\footnotemark[1]$^{,}$\footnotemark[3]$^{,}$\footnotemark[5]$^{,}$\footnotemark[6]  \\
%\textit{\nocorr
Instituto de Matem\'aticas y F\'{\i}sica Fundamental,\\
Consejo Superior de Investigaciones Cient\'{\i}ficas, \\
Serrano 113 bis,
28006 Madrid, Spain \\
E-mail: \texttt{oscar.garcia-prada@uam.es}
%}

\medskip
\noindent 
\textbf{Peter B. Gothen}
\footnotemark[1]$^{,}$\footnotemark[4]$^{,}$\footnotemark[5]   \\
%\textit{\nocorr
Departamento de Matem{\'a}tica Pura, \\
Faculdade de Ci{\^e}ncias,
Universidade do Porto, \\
Rua do Campo Alegre 687,
4169-007 Porto,
Portugal \\
E-mail: \texttt{pbgothen@fc.up.pt}
%}
\footnotetext[1]{Members of VBAC (Vector Bundles on Algebraic
  Curves), which is partially supported by EAGER (EC FP5 Contract no.\ 
  HPRN-CT-2000-00099) and by EDGE (EC FP5 Contract no.\ 
  HPRN-CT-2000-00101).}  \footnotetext[2]{Partially supported by the
  National Science Foundation under grant DMS-0072073 }
\footnotetext[3]{Partially supported by the Ministerio de Ciencia y
  Tecnolog\'{\i}a (Spain) under grant BFM2000-0024} \footnotetext[4]{
  Partially supported by the Funda{\c c}{\~a}o para a Ci{\^e}ncia e a
  Tecnologia (Portugal) through the Centro de Matem{\'a}tica da
  Universidade do Porto and through grant no.\ SFRH/BPD/1606/2000.}
\footnotetext[5]{Partially supported by the Portugal/Spain bilateral
  Programme Acciones Integradas, grant nos.\ HP2000-0015 and AI-01/24}
\footnotetext[6]{Partially supported by a British EPSRC grant
  (October-December 2001)}

\bigskip
\noindent
\textbf{\today} \vfill

\newpage

\pagestyle{headings}

\paragraph{\textbf{Abstract.}}
Using the $L^2$ norm of the Higgs field as a Morse function, we study 
the moduli spaces of $\U(p,q)$-Higgs bundles over a Riemann surface. 
We require that the genus of the surface be at least two, but place 
no constraints on $(p,q)$.  A key step is the identification of the 
function's local minima as moduli spaces of holomorphic triples.  In 
a companion paper \cite{bradlow-garcia-gothen:2002:triples} we prove 
that these moduli spaces of triples are non-empty and irreducible.

Because of the relation between flat bundles and fundamental group
representations, we can interpret our conclusions as results about the
number of connected components in the moduli space of semisimple
$\PU(p,q)$-representations.  The topological invariants of the flat
bundles are used to label subspaces. These invariants are bounded by 
a Milnor--Wood type inequality. For each allowed value of the 
invariants satisfying a certain coprimality condition, we prove that 
the corresponding subspace is non-empty and connected. If the 
coprimality condition does not hold, our results apply to the closure 
of the moduli space of irreducible representations. 

\newpage
%\tableofcontents
%\newpage

%%%%%%%%%%%%%%%%%%%%%
\section{Introduction} \label{sec:introduction}
%%%%%%%%%%%%%%%%%%%%%%

The relation between Higgs bundles and fundamental group
representations provides a vivid illustration of the interaction
between geometry and topology. On the topological side we have a
closed oriented surface $X$ and the moduli space (or character
variety) of representations of $\pi_1X$ in a Lie group $G$. We cross
over to complex geometry by fixing a complex structure on $X$, thereby
turning it into a Riemann surface. The space of representations, or
equivalently the space of flat $G$-bundles, then emerges as a complex
analytic moduli space of $G$-Higgs bundles. In this guise, the moduli
space carries a natural proper function whose restriction to the
smooth locus is a Morse-Bott function. We can therefore use this
function to determine topological properties of the moduli space of
representations. Our goal in this paper is to pursue these ideas in
the case where the group $G$ is the real Lie group $\PU(p,q)$, the
adjoint form of the non-compact group $\U(p,q)$.

The relevant Higgs bundles in our situation are $\U(p,q)$-Higgs 
bundles. These can be seen as a special case of the 
$G$-Higgs bundles defined by Hitchin in \cite{hitchin:1987}, 
where $G$ is a real form  of a complex reductive Lie group. Such 
objects provide a natural generalization of holomorphic vector 
bundles, which correspond to the case $G=\U(n)$ and zero Higgs field. 
In particular, they permit an extension to other groups of the 
Narasimhan and Seshadri theorem (\cite{narasimhan-seshadri:1965}) on 
the relation between unitary representations of $\pi_1 X$ and stable 
vector bundles.  By embedding $\U(p,q)$ in $\GL(p+q)$ we can give a 
concrete description 
  of a $\U(p,q)$-Higgs bundle as a pair
 \begin{equation}\label{eqtn:intro1}
 (V\oplus W,\Phi=\left(
   \begin{smallmatrix}
     0 & \beta \\
     \gamma & 0
   \end{smallmatrix}
   \right))\
 \end{equation}
 \noindent where $V$ and $W$ are holomorphic vector bundles
 of rank $p$ and $q$ respectively, $\beta$ is a section in
 $H^0(\Hom(W,V) \otimes K)$, and $\gamma \in H^0(\Hom(V,W) \otimes K)$,
 so that $\Phi\in H^0(\mathrm{End}(V\oplus W) \otimes K)$. 
 
 \par By the work
  of Hitchin \cite{hitchin:1987,hitchin:1992} Donaldson \cite{donaldson:1987},
 Simpson \cite{simpson:1988,simpson:1992,simpson:1994a,simpson:1994b} and 
 Corlette \cite{corlette:1988}, we can define moduli spaces of
 polystable Higgs bundles, and these can be identified with moduli spaces of
 solutions to natural gauge theoretic equations. Moreover, since the
 gauge theory equations amount to a projective flatness condition,
 these moduli spaces correspond to moduli spaces of flat
 structures. In the case of
 $\U(p,q)$-Higgs bundles, the flat structures correspond to semi-simple
 representations of $\pi_1 X$ into the group $\PU(p,q)$. The Higgs
 bundle moduli spaces can thus be used, in a way which we make precise
 in Sections
 \ref{sec:background} and \ref{sec:higgs-bundles},  to study
 the representation variety
 $$\mathcal{R}(\PU(p,q)) = \Hom^+(\pi_1 X, \PU(p,q)) / \PU(p,q)\ ,$$
 \noindent where $\Hom^+(\pi_1 X, \PU(p,q))$ denotes the set of semi-simple
 representations of $\pi_1X$ in $\PU(p,q)$, and the quotient is by the
 adjoint action.
 
 \par Our main tool for studying the topology of the Higgs moduli space is 
 the function which measures the $L^2$-norm of the Higgs field. When the
 moduli space is smooth, this turns out to provide a suitably
 non-degenerate Bott-Morse function which is, moreover, a proper map.
 In some cases (cf.\ \cite{hitchin:1987,gothen:1995,gothen:2002}) the critical
 submanifolds are well enough understood to allow the extraction of
 topological information as detailed as the Poincar\'e polynomial.  In
 our case our understanding is confined to the local minima of the function.  
 This is sufficient to allow us to count the number of
 components of the Higgs moduli spaces, and thus of the representation
 varieties.  A trivial but important observation is that the
 properness of the function allows us to draw conclusions about
 connected components also in the non-smooth case; we shall
 henceforth, somewhat imprecisely, refer to the function as the
 ``Morse Function'', whether or not the moduli space is smooth.
  
The criterion we use for finding the local minima can be applied more 
generally, for instance if $\U(p,q)$ is replaced by any real form of 
a complex reductive group.  This should provide an important tool for 
future research. In the present case, this criterion allows us to 
identify the subspaces of local minima as moduli spaces in their own 
right, namely as moduli spaces of the holomorphic triples introduced 
in \cite{bradlow-garcia-prada:1996}.  In a companion paper 
\cite{bradlow-garcia-gothen:2002:triples} we 
 develop the theory of such objects and their moduli spaces. Using the
 results of \cite{bradlow-garcia-gothen:2002:triples} we are able to
 deduce several results about the Higgs moduli spaces and also about
 the corresponding representation spaces.
 
 The relation between Higgs bundles and surface group representations
 has been successfully exploited by others, going back originally to
 the work of Hitchin and Simpson on complex reductive groups. The use
 of Higgs bundle methods to study $\mathcal{R}(G)$ for real $G$ was pioneered 
 by Hitchin in \cite{hitchin:1992}, and further developed in 
 \cite{gothen:1995,gothen:2001}. It has also been used  by Xia and Xia-Markman 
 (in \cite{xia:1997, xia:1999, xia:2000, markman-xia:2001}) to study various 
 special cases of $G=\PU(p,q)$. None of these, though, address the general 
case of $\PU(p,q)$, as we do in this paper. 
 
 We now give a brief summary of the contents and main results of this paper.
 
 In Sections \ref{sec:background} and \ref{sec:higgs-bundles} we give
 some background and describe the basic objects of our study.  In
 Section \ref{sec:background} we describe the natural invariants
 associated with representations of $\pi_1 X$ into $\PU(p,q)$. We also 
 discuss the invariants associated with representations of $\Gamma$, 
 the universal central extensions of
 $\pi_1$, into $\U(p,q)$. The space of such representations is denoted
 by $\mathcal{R}_{\Gamma}(\U(p,q))$. In both cases, these involve a pair of
 integers $(a,b)$ which can be interpreted respectively as degrees of
 rank $p$ and rank $q$ vector bundles over $X$. In the case of the
 $\PU(p,q)$ representations, the pair is well defined only as a class
 in a quotient
 $\Z\oplus\Z/(p,q)\Z$. This leads us to define subspaces
 $\mathcal{R}[a,b]\subset \mathcal{R}(\PU(p,q))$ and
 $\mathcal{R}_{\Gamma}(a,b)\subset \mathcal{R}_{\Gamma}(\U(p,q))$. For
 fixed $(a,b)$, the space $\mathcal{R}_{\Gamma}(a,b)$ fibers over
 $\mathcal{R}[a,b]$ with connected fibers.
 
 In section \ref{sec:higgs-bundles} we define $\U(p,q)$-Higgs bundles
 and their moduli spaces and establish their essential properties.
 Thinking of a $\U(p,q)$-Higgs bundle as a pair $(V\oplus W,\Phi)$,
 the parameters $(a,b)$ appear here as the degrees of the bundles $V$
 and $W$.  The moduli space of polystable $\U(p,q)$-Higgs bundles with
 $\deg(V)=a$ and $\deg W=b$, which we denote by $\mathcal{M}(a,b)$, is
 the space that can be identified with the component
 $\mathcal{R}_{\Gamma}(a,b)$ of $\mathcal{R}_{\Gamma}(\U(p,q))$. This,
 together with the fibration over $\mathcal{R}_{\Gamma}(\U(p,q))$ are
 the crucial links between the Higgs moduli and the surface group
 representation varieties.

Fixing $p,q,a$ and $b$, we begin the Morse theoretic analysis of 
$\mathcal{M}(a,b)$ in Section  \ref{sec:morse-theory}. The basic 
results we need (cf.\ Proposition \ref{prop:topology-exercise}) are 
that the $L^2$-norm of the Higgs field has a minimum on each 
connected component of ${\mathcal{M}}(a,b)$, and hence if the 
subspace of local minima is connected then so is $\mathcal{M}(a,b)$. 
We identify the local minima, the loci of which we denote by 
 $\mathcal{N}(a,b)$,  and prove (cf.\ Theorem
 \ref{thm:minima} and Proposition \ref{prop:vanishing}) that these
 correspond precisely to holomorphic triples in the sense of
 \cite{bradlow-garcia-prada:1996}. A full treatment of holomorphic triples and
 their moduli spaces is given in \cite{bradlow-garcia-gothen:2002:triples}.
 We summarize the  salient features of these moduli spaces in Section 
\ref{sec:stable-triples}.

In section \ref{sec:main-results} we knit together all the strands. 
Using the properties of the moduli spaces of triples, we establish 
the key (for our purposes) topological properties of the strata 
$\mathcal{N}(a,b)$. These lead directly to our main results for the 
moduli spaces 
$\mathcal{M}(a,b)$. Some of the results depend on $(a,b)$ only in
the combination 
 $$\tau=\tau(a,b)=2\frac{aq-bp}{p+q}\ ,$$
known as the Toledo invariant. Indeed, $(a,b)$ is constrained by the 
bounds  $0\le|\tau|\le\tau_M$,  where $\tau_M=2\min\{p,q\}(g-1)$. 
Originally proved by Domic and Toledo in \cite{domic-toledo:1987}, 
these bounds emerge naturally from our point of view (cf.\ Corollary 
\ref{cor:toledo} and Remark  \ref{remark:MWbound}). Bounds on 
invariants of this type, for representations of finitely generated 
groups in $\U(p,q)$, have also recently been studied using techniques 
from ergodic theory (see \cite{burger-iozzi}). Summarizing our main 
results, we prove 

\begin{theoremA}[Theorems \ref{thm:summary-M(a,b)} and \ref{thm:coprime-M(a,b)}]
  
 Fix positive integers $(p,q)$.  Take $(a,b)\in\Z\oplus\Z$ and let
  $\tau(a,b)$ be the Toledo invariant. Let 
  $\mathcal{M}^s(a,b)\subseteq \mathcal{M}(a,b)$ denote the 
  moduli space of {\it strictly stable} $\U(p,q)$-Higgs bundles.
\begin{itemize}
  
\item[$(1)$] $\mathcal{M}(a,b)$ is non-empty
if and only if  $0\le|\tau(a,b)|\le\tau_M$. If 
$\tau(a,b)=0$, or $|\tau(a,b)|=\tau_M$ and $p\ne q$ then $\mathcal{M}^s(a,b)$
is empty; otherwise it is non-empty whenever $\mathcal{M}(a,b)$ is 
non-empty.

\item[$(2)$] If $|\tau(a,b)|=0$ or $|\tau(a,b)|=\tau_M$ and $p\ne q$
then $\mathcal{M}(a,b)$ is connected. 

\item[$(3)$] Whenever non-empty, the moduli space 
$\mathcal{M}^s(a,b)$ is a smooth manifold of the expected 
dimension (i.e.\ $1+(p+q)^2(g-1)$), with connected closure 
$\bar{\mathcal{M}}^s(a,b)\subseteq\mathcal{M}(a,b)$. In these cases,
if $\mathcal{M}(a,b)$ has more than one connected component, then 
$\GCD(p+q,a+b)\ne 1$  and, if $p=q$, $0<|\tau|\le (p-1)(2g-2)$. 
\end{itemize}
\end{theoremA}

\begin{theoremB}[Theorem \ref{prop:rigidity}]
  
Suppose that $p\ne q$ and $(a,b)\in\Z\oplus\Z$ are such that 
$|\tau(a,b)|=\tau_M$. To  be specific, suppose that $p < q$ and 
$\tau(a,b)=p(2g-2)$. Then every element in 
$\mathcal{M}(a,b)$ decomposes as the 
direct sum of a polystable $\U(p,p)$-Higgs bundle with maximal Toledo 
invariant and a polystable vector bundle of rank $q-p$. Thus 
  \begin{equation}
    \mathcal{M}(p,q,a,b)\cong \mathcal{M}(p,p,a,a-p(2g-2))
\times M(q-p, b-a +p(2g-2)).
  \end{equation}
  In particular, the smooth locus in
  $\mathcal{M}(p,q,a,b)$ has dimension $2+ (q^2+5p^2-2pq)(g-1)$.  
  This is strictly smaller than the expected dimension if $g\ge 2$. 
  
  (A similar result holds if $p>q$, as well as if $\tau=-p(2g-2)$).
\end{theoremB}

Since we identify $\mathcal{M}(a,b)=\mathcal{R}_\Gamma(a,b)$, we can 
translate these results directly into statements about 
$\mathcal{R}_\Gamma(a,b)$ (given in Theorems 
\ref{thm:results-R-Gamma} and \ref{thm:RgammaTmax}).  The subspace
in $\mathcal{R}_\Gamma(a,b)$ which corresponds to 
$\mathcal{M}^s(a,b)\subseteq \mathcal{M}(a,b)$ is denoted by 
$\mathcal{R}^*_\Gamma(a,b)$. The representations it labels include all the 
simple representations. Defining 
$\mathcal{R}^*_{\Gamma}(\U(p,q))\subset\mathcal{R}_{\Gamma}(\U(p,q))$
to be the union over all $(a,b)$ of the components 
$\mathcal{R}^*_\Gamma(a,b)$ we thus obtain

\bigskip

\noindent{\bf Theorem C (Corollary \ref{cor:number-of-[a,b]}) } {\it 
The moduli space 
$\mathcal{R}^*_{\Gamma}(\U(p,q))$ has 
$$2(p+q)\min\{p,q\}(g-1)+\GCD(p,q)$$
\noindent connected components.} 

\bigskip

\noindent Since $\mathcal{R}_{\Gamma}(a,b)$ fibers over  
$\mathcal{R}[a,b]$ with connected fibers, we can apply our results to the
latter. The results are given in Theorems \ref{thm:results-R1} and 
\ref{thm:results-R2}.

The above results fall just short of saying that the full moduli 
spaces $\mathcal{M}(a,b)$ ($ = \mathcal{R}(a,b)$) and 
$\mathcal{R}[a,b]$ are connected for all allowed choices of $(a,b)$. 
They show however that if any one is not connected then it has one 
(non-empty) connected component which contains all the irreducible 
objects. Any other components must thus consist entirely of reducible 
(or strictly semisimple) elements. Theorem B and its analogs for 
$\mathcal{R}_{\Gamma}(a,b)$ and $\mathcal{R}[a,b]$
 generalize rigidity results of Toledo 
\cite{toledo:1989} (when 
$p=1)$ and Hern\'andez 
\cite{hernandez:1991} (when $p=2$).

This paper, together with its companion 
\cite{bradlow-garcia-gothen:2002:triples} form a substantially revised
version of the preprint \cite{bradlow-garcia-gothen:2002}. The main 
results proved in this paper were announced in the note 
\cite{bradlow-garcia-gothen:2001}.  In that note we claim (without proof)
that the connectedness results for the moduli spaces 
$\mathcal{R}(a,b)$ and $\mathcal{R}[a,b]$  hold without the above 
qualifications.  This is a reasonable conjecture, which we hope to 
come back to in a future publication. 

We note, finally, that our methods surely apply more widely than to
$\U(p,q)$-Higgs bundles and $\PU(p,q)$ representations (see, for
example, Remark~\ref{rem:HiggsG-minima}). Moreover, careful scrutiny
of the Lie algebra properties used in the proofs suggests certain
aspects can be generalized to representations in any real group $G$
for which $G/H$ is hermitian symmetric, where $H\subset G$ is a
maximal compact subgroup.  This will be addressed in a future
publication.

\textbf{Acknowledgements.} We thank the mathematics departments
of the University of Illinois at Urbana-Champaign, the University 
Aut{\'o}noma of  Madrid and the University of Aarhus, the Department 
of Pure Mathematics of the University of Porto, the Mathematical 
Sciences Research Institute of Berkeley and the Mathematical 
Institute of the University of Oxford, and the Erwin
Schr\"odinger International Institute for Mathematical
Physics in Vienna for their hospitality during 
various stages of this research.  We thank Fran Burstall, Bill 
Goldman, Nigel Hitchin, Eyal Markman, S. Ramanan, Domingo Toledo, and 
Eugene Xia, for many insights and patient explanations.

%%%%%%%%%%%%%%%%%%%%%%%%%%%%%%%%%%%%%%%%%%%%%%%%%%%%%%%%%%%%%%
\section{Representations of surface groups}\label{sec:background}
%%%%%%%%%%%%%%%%%%%%%%%%%%%%%%%%%%%%%%%%%%%%%%%%%%%%%%%%%%%%%%

In this section we record some general facts about representations of 
a surface group in $\U(p,q)$ or $\PU(p,q)$ and set up our notation. A 
very useful reference for the general theory is Goldman's paper 
\cite{goldman:1985}. 

%%%%%%%%%%%%%%%%%%%%%%%%%
\subsection{Moduli spaces of representations}
%%%%%%%%%%%%%%%%%%%%%%%%%

Let $X$ be a closed oriented surface of genus $g \geq 2$. By 
definition 
$\U(p,q)$ is the subgroup of $\GL(n,\C)$ (with 
$n=p+q$) which leaves invariant a hermitian form of signature $(p,q)$.
It is a non-compact real form of $\GL(n,\C)$ with center $\U(1)$ and
maximal compact subgroup $\U(p)\times\U(q)$.  The quotient
$\U(p,q)/(\U(p)\times\U(q))$ is a hermitian symmetric space. The
adjoint form $\PU(p,q)$ is given by the exact sequence of groups
$$
1\longrightarrow \U(1)\longrightarrow
\U(p,q)\longrightarrow\PU(p,q)\longrightarrow 1\ ,
$$
and we have a standard inclusion $\PU(p,q)\subset \PGL(n,\C)$.

\begin{definition}
\label{def:RG}
By a \emph{representation} of $\pi_1 X$ in $\PU(p,q)$ we mean a 
homomorphism 
$\rho \colon \pi_1 X \to \PU(p,q)$.  We say that a representation of $\pi_1
X$ in $\PU(p,q)$ is \emph{semi-simple} if the induced (adjoint) 
representation on the Lie algebra of $\PU(p,q)$ is semi-simple. The 
group $\PU(p,q)$ acts on the set of representations via conjugation. 
Restricting to the semi-simple representations, we get the 
\emph{moduli space} of representations,
\begin{equation}\label{eqn:RGdef}
  \mathcal{R}(\PU(p,q)) = \Hom^+(\pi_1 X, \PU(p,q)) / \PU(p,q)\ .
\end{equation}

\end{definition}
The moduli space of representations can be described more concretely
as follows.  {}From the standard presentation
\begin{displaymath}
  \pi_1 X = \langle A_{1},B_{1}, \ldots, A_{g},B_{g} \suchthat
  \prod_{i=1}^{g}[A_{i},B_{i}] = 1 \rangle
\end{displaymath}
we see that $\Hom^{+}(\pi_1 X, \PU(p,q))$ can be embedded in 
$\PU(p,q)^{2g}$ via 
\begin{align*}
  \Hom^{+}(\pi_1 X, \PU(p,q)) &\to \PU(p,q)^{2g} \\
  \rho &\mapsto (\rho(A_1), \ldots \rho(B_g)).
\end{align*}
We  give 
$\Hom^{+}(\pi_1 X, \PU(p,q))$ the subspace topology and 
$\mathcal{R}(\PU(p,q))$ the quotient topology.  This topology is Hausdorff 
because we have restricted attention to semi-simple representations. 

Clearly any representation of $\pi_1 X$ in $\U(p,q)$ gives rise to a 
representation in $\PU(p,q)$; however, not all representations in 
$\PU(p,q)$ lift to $\U(p,q)$.  We are thus motivated 
to consider representations of the central extension
\begin{equation}\label{eq:gamma}
0\longrightarrow\mathbb{Z}\longrightarrow\Gamma\longrightarrow\pi_1 
X\longrightarrow 1 \ .
\end{equation} 
\noindent Such extensions are defined (as in
\cite{atiyah-bott:1982}) by the generators $A_{1},B_{1}, \ldots,
A_{g},B_{g}$ and a central element $J$ subject to the relation 
$\prod_{i=1}^{g}[A_{i},B_{i}] = J$.  With $\Gamma$ thus defined, 
any representation of $\pi_1 X$ in 
$\PU(p,q)$ can be lifted to a representation of $\Gamma$ in
$\U(p,q)$. 

In analogy with Definition~\ref{def:RG} we make the following
definition.
\begin{definition}
\label{def:RgG}
We define the \emph{moduli space} of
semi-simple representations of $\Gamma$ in $\U(p,q)$ by
\begin{equation}\label{eq:RgGdef}
  \mathcal{R}_\Gamma(\U(p,q)) = \Hom^+(\Gamma, \U(p,q)) / \U(p,q)\ ,
\end{equation}
where semi-simplicity is defined with respect to the induced adjoint 
representation. This space is topologized in the same way as 
$\mathcal{R}(\PU(p,q))$. 
\end{definition}

%%%%%%%%%%%%%%%%%%%%%
\subsection{Invariants}\label{subs:invar}
%%%%%%%%%%%%%%%%%%%%%%%

Our basic objective is to study the number of connected components of
the spaces $\mathcal{R}(\PU(p,q))$ and $\mathcal{R}_\Gamma(\U(p,q))$.
The first step in the study of topological properties of these spaces
is to identify the appropriate topological invariant of a
representation $\rho \colon \pi_1 X \to G$.  For a general connected
Lie group $G$ the relevant invariant is an obstruction class in
$H^2(X,\pi_1 G) \cong \pi_1 G$ (see Goldman
\cite{goldman:1985,goldman:1988}).  In the following we give an
explicit description of this invariant in our case, using
characteristic classes of the flat bundles associated to
representations of the fundamental group.  In fact we shall not need
the more general description of the invariant.

We begin by considering the case $G = \U(p,q)$.  By the same argument
as in \cite{atiyah-bott:1982}\footnotemark \footnotetext{While
  \cite{atiyah-bott:1982} gives the argument for $\U(n)$ and $\PU(n)$,
  there are no essential changes to be made in order to adapt for the
  case of $\U(p,q)$ and $\PU(p,q)$.}, $\mathcal{R}_{\Gamma}(\U(p,q))$
can be identified with the moduli space of connections with central
curvature on a fixed $\U(p,q)$-bundle on $X$.  Taking a reduction to the
maximal compact $\U(p)\times\U(q)$, we thus associate to each class
$\tilde{\rho}\in \mathcal{R}_{\Gamma}(\U(p,q))$ a vector bundle of the
form $V\oplus W$, where $V$ and $W$ are rank $p$ and $q$ respectively,
and thus a pair of integers $(a,b) =(\deg(V),\deg(W))$.  There is thus
a map
\begin{displaymath}
  \tilde{c} \colon  \mathcal{R}_{\Gamma}(\U(p,q)) \to \Z\oplus\Z
\end{displaymath}
given by $\tilde{c}(\tilde{\rho}) = (a,b)$.  The corresponding map on $\Hom^+(\Gamma, \U(p,q))$ is 
clearly continuous and thus locally constant. Since 
$\U(p,q)$ is connected, the map $\tilde{c}$ is likewise 
continuous and thus constant on connected components.
We make the following
definition.

\begin{definition}
  The subspace of $\mathcal{R}_{\Gamma}(\U(p,q))$ corresponding to
  representations with invariants $(a,b)$ is denoted by
\begin{align*}
  \mathcal{R}_\Gamma(a,b) &= \tilde{c}^{-1}(a,b)\\
  &= \{ \tilde{\rho} \in \mathcal{R}_\Gamma(\U(p,q)) \suchthat
  \tilde{c}(\tilde{\rho})=(a,b) \in \Z \oplus \Z\}\ .
\end{align*}
\end{definition}
Note that $\mathcal{R}_\Gamma(a,b)$ is a union of connected
components, because $\tilde{c}$ is constant on each connected
component.

Next we consider the case $G = \PU(p,q)$.  Any flat $\PU(p,q)$-bundle
lifts to a $\U(p,q)$-bundle with a connection with constant central
curvature.  This lift is, however, not uniquely determined: in fact
two such $\U(p,q)$-bundles give rise to the same flat
$\PU(p,q)$-bundle if and only if one can be obtained from the other by
twisting with a line bundle $L$ with a unitary connection of constant
curvature.  If the invariant of the $\U(p,q)$-bundle is $(a,b)$ and
the degree of $L$ is $l$, then the invariant associated to the 
twisted bundle is $(a+pl,b+ql)$.  There is thus a well defined map 
\begin{equation}\label{eq:invariants}
c\colon \mathcal{R}(\PU(p,q))\longrightarrow (\Z \oplus \Z) / (p,q)\Z\ ,
\end{equation}
where $(\Z \oplus \Z) / (p,q)\Z$ denotes the quotient of $\Z \oplus
\Z$ by the $\Z$-action $l\cdot (a,b)=(a+pl,b+ql)$. 
Notice that $(\Z \oplus \Z)/(p,q)\Z$ can be identified 
with $\pi_1(\PU(p,q))$. The invariant defined by $c$ is 
the same as the obstruction class defined by Goldman 
\cite{goldman:1985,goldman:1988}.

\begin{definition} Denote the 
image of $(a,b)$ in $(\Z \oplus \Z) / (p,q)\Z$ by  $[a,b]$.  The 
subspace of $\mathcal{R}(\PU(p,q))$ corresponding to 
  representations with invariant $[a,b]$ is denoted by
\begin{align*}
  \mathcal{R}[a,b] &= c^{-1}[a,b]\\
  &= \{ \rho \in \mathcal{R}(\PU(p,q)) \suchthat c(\rho)=[a,b] \in (\Z
  \oplus\Z)/(p,q)\Z \}\ .
\end{align*}
\end{definition}
The space $\mathcal{R}[a,b]$ is a union of connected components in 
the same way as $\mathcal{R}_\Gamma(a,b)$. In order to compare the 
spaces $\mathcal{R}_\Gamma(a,b)$ and 
$\mathcal{R}[a,b]$ notice that we have surjective maps
\begin{equation}
  \label{eq:principal-jac}
  \mathcal{R}_\Gamma(a,b) \to \mathcal{R}[a,b].
\end{equation}
\noindent Moreover, the preimage 
\begin{equation}
\pi^{-1}(\mathcal{R}[a,b])=\bigcup_{(a,b)}\mathcal{R}_\Gamma(a,b) 
\end{equation}
\noindent where the union is over all $(a,b)$ in the class
$[a,b]\in (\Z \oplus\Z)/(p,q)\Z$. As mentioned above, tensoring by
line bundles of degree $l$ with constant curvature connections gives 
an isomorphism 
\begin{displaymath}
  \mathcal{R}_\Gamma(a,b)
  \xrightarrow{\cong}
  \mathcal{R}_\Gamma(a+pl,b+ql)\ .
\end{displaymath}
Notice that if $c(\rho)=[a,-a]$ for a representation
$\rho\in\mathcal{R}(\PU(p,q))$, then the associated $\U(p,q)$-bundle
can be taken to have degree zero and the projectively flat connection 
is actually flat.  Then $\rho$ defines a representation of $\pi_1 X$ 
in $\U(p,q)$. Under the correspondence between 
$\mathcal{R}(\PU(p,q))$ and $\mathcal{R}_{\Gamma}(\U(p,q))$, 
$\rho$ corresponds to a $\Gamma$ representation in which the central
element $J$ acts trivially. Furthermore, the subspaces 
$\mathcal{R}_\Gamma(a,-a)\subset 
\mathcal{R}_{\Gamma}(\U(p,q))$ can be identified with components
of $\mathcal{R}(\U(p,q))$ (the moduli space for representations 
of 
$\pi_1 X$ in $\U(p,q)$). Indeed, defining 
\begin{equation}\label{eqn:RUpq}
\mathcal{R}(a) = \mathcal{R}_\Gamma(a,-a)\ , 
\end{equation}
we see that $\mathcal{R}(\U(p,q))$ is a union over $a\in\Z$ of the 
subspaces  $\mathcal{R}(a)$. 

Finally, we observe that the moduli space of flat degree zero line 
bundles acts by tensor product of bundles on 
$\mathcal{R}_\Gamma(a,b)$.  Since this moduli space 
is isomorphic to the torus $\U(1)^{2g}$, we get the
following relation between connected components.

\begin{proposition}\label{prop:principal-jac} The map 
  $\mathcal{R}_\Gamma(a,b) \to \mathcal{R}[a,b]$ given in
  \eqref{eq:principal-jac} defines a $\U(1)^{2g}$-fibration. Thus the
  subspace $\mathcal{R}[a,b] \subseteq \mathcal{R}(\PU(p,q))$ is
  connected if $\mathcal{R}_\Gamma(a,b)$ is connected.  \hfill\qed
\end{proposition}

%%%%%%%%%%%%%%%%%%%%%%%%%%%%%%%%%%%%%%%%%%%%%%%%
\section{Higgs bundles and flat connections}
\label{sec:higgs-bundles}
%%%%%%%%%%%%%%%%%%%%%%%%%%%%%%%%%%%%%%%%%%%%%%%%
%%%%%%%%%

We study the moduli spaces of representations by choosing a complex
structure on $X$.  This allows us to identify these spaces with 
certain moduli spaces of Higgs bundles.  In this section we explain 
this correspondence and recall some general facts about Higgs 
bundles. Following this, we describe the special class of Higgs 
bundles relevant for the study of representations in $\PU(p,q)$ and 
$\U(p,q)$ and derive some basic results about these moduli spaces. 

%%%%%%%%%%%%%%%%%%%%%%%%%%%%
\subsection{$\GL(n,\C)$-Higgs bundles}
%%%%%%%%%%%%%%%%%%%%%%%%%%%%%
Give $X$ the structure of a Riemann surface.  We recall (from
\cite{corlette:1988,donaldson:1987,hitchin:1987,simpson:1988,
simpson:1994a,simpson:1994b})
the following definition and basic facts about $\GL(n,\C)$-Higgs 
bundles. 
\begin{definition}
\label{defn:GLnHiggs}

(1) A \emph{$\GL(n,\C)$-Higgs bundle} on $X$ is a pair
  $(E,\Phi)$, where $E$ is a rank $n$ holomorphic vector bundle over
  $X$ and $\Phi \in H^0(\End(E) \otimes K)$ is a holomorphic
  endomorphism of $E$ twisted by the canonical bundle $K$ of $X$.

\bigskip

\noindent (2) The $\GL(n,\C)$-Higgs bundle $(E,\Phi)$ is \emph{stable} if
  the slope stability condition
\begin{equation}\label{eq:stability}
\mu(E') < \mu(E)
\end{equation} 
holds for all proper $\Phi$-invariant subbundles $E'$ of $E$.  Here
the \emph{slope} is defined by $\mu(E)=\deg(E)/\rk(E)$ and
\emph{$\Phi$-invariance} means that $\Phi(E')\subset E'\otimes K$.
\emph{Semistability} is defined by replacing the above strict
inequality with a weak inequality. A Higgs bundle is called
\emph{polystable} if it is the direct sum of stable Higgs bundles with
the same slope.

\bigskip

\noindent (3) Given a hermitian metric on $E$, let $A$ denote the
  unique unitary connection compatible with the holomorphic structure, and let
  $F_A$ be its curvature.
  \emph{Hitchin's equations} on $(E,\Phi)$ are
\begin{equation}
  \label{eq:hitchin1}
  \begin{aligned}
    F_A + [\Phi,\Phi^*] &= -\sqrt{-1}\mu \text{Id}_E \omega, \\
    \dbar_{A} \Phi &=0,
  \end{aligned}
\end{equation}
where $\mu$ is a constant, $\text{Id}_E$ is the identity on $E$,
$\dbar_A$ is the anti-holomorphic part of the covariant derivative
$d_A$ and $\omega$ is the K\"ahler form on $X$.  If we normalize
$\omega$ so that $\int_X \omega = 2\pi$ then, taking the trace and
integrating over $X$ in the first equation, one sees that $\mu =
\mu(E)$. A solution to Hitchin's equations is \emph{irreducible} if
there is no proper subbundle of $E$ preserved by $A$ and $\Phi$.
\end{definition}

\begin{theorem}\label{thm:HK}
%\begin{itemize}
%\item[$(1)$] 
(1) Let $(E,\Phi)$ be a $\GL(n,\C)$-Higgs bundle. Then $(E,\Phi)$ is 
polystable if and only if it admits a hermitian metric such that 
Hitchin's equations \eqref{eq:hitchin1} are satisfied.  Moreover,
$(E,\Phi)$ is stable if and only if the corresponding solution is
irreducible. 

%\item[$(2)$]
(2) Fix a hermitian metric in a smooth rank $n$ complex
  vector bundle on $X$, then there is a gauge theoretic moduli space 
  of pairs $(A,\Phi)$, consisting of a unitary
  connection $A$ and an endomorphism valued $(1,0)$-form $\Phi$, which
  are solutions to Hitchin's equations \eqref{eq:hitchin1}, modulo
  $\U(n)$-gauge equivalence.

%\item[$(3)$]
(3) The moduli space of rank $n$ degree $d$ polystable
  Higgs bundles is a quasi-projective variety of complex
  dimension $2(1+n^2(g-1))$.  There is a map from the gauge theoretic
  moduli space to this moduli space given by taking a solution
  $(A,\Phi)$ to Hitchin's equations to the Higgs bundle $(E,\Phi)$,
  where the holomorphic structure on $E$ is given by $\dbar_{A}$.
  This map is a homeomorphism, and a diffeomorphism on the smooth 
  locus.

%\item[$(4)$]
(4) If we define a Higgs connection (as in
  \cite{simpson:1994a}) by
\begin{equation}\label{eqn:Dhiggs} 
D=d_A+\theta
\end{equation}
\noindent where $\theta=\Phi+\Phi^* $, then Hitchin's equations are 
equivalent to the conditions
\begin{equation}
  \label{eq:harmonic1}
  \begin{aligned}
  F_D =& -\sqrt{-1}\mu \text{Id}_E \omega, \\
  &d_A \theta =0, \\
  &d_A^* \theta =0.
  \end{aligned}
\end{equation}

In particular, $D$ is a projectively flat connection. If $\deg(E)=0$ 
then $D$ is actually 
  flat. It follows that in this case the pair $(E,D)$ defines a
  representation of $\pi_1 X$ in $\GL(n,\C)$. If $\deg(E)\ne 0$, then the
  pair $(E,D)$ defines a representation of $\pi_1 X$ in $\PGL(n,\C)$, or
  equivalently, a representation of $\Gamma$ in $\GL(n,\C)$.  By the
  theorem of Corlette (\cite{corlette:1988}), every semisimple
  representation of $\Gamma$ (and therefore every semisimple
  representation of $\pi_1 X$) arises in this way.

%\item[$(5)$]
(5) This correspondence gives rise to a homeomorphism between
  the moduli space of polystable Higgs bundles of rank $n$ and the
  moduli space of semisimple representations of $\Gamma$ in
  $\GL(n,\C)$.  If the degree of the Higgs bundle is zero, then the
  moduli space is homeomorphic to the moduli space of
  representations of $\pi_1 X$ in $\GL(n,\C)$.
%\end{itemize}
\end{theorem}

%%%%%%%%%
\subsection{$\U(p,q)$-Higgs bundles}\label{subs:upq}
%%%%%%%%%
If we fix integers $p$ and $q$ such that $n=p+q$, then we can isolate 
a special class of $\GL(n,\C)$-Higgs bundles by the requirements that 
\begin{equation} \label{upq-higgs-bundle}
  \begin{aligned}
  E &= V \oplus W \\
  \Phi &=
  \left(
  \begin{smallmatrix}
    0 & \beta \\
    \gamma & 0
  \end{smallmatrix}
  \right)
  \end{aligned}
\end{equation}
where $V$ and $W$ are holomorphic vector bundles of rank $p$ and $q$ 
respectively and the non-zero components in the Higgs field are 
$\beta \in H^0(\Hom(W,V) 
\otimes K)$, and 
$\gamma 
\in H^0(\Hom(V,W) 
\otimes K)$. 

The form of the Higgs field is determined by the Lie theory of the 
symmetric space $\U(p,q)/(\U(p)\times\U(q))$. Recall that for any real 
form $G$ of a complex reductive group $G^{\C}$, with maximal compact 
subgroup $H$, there is an $\Ad$-invariant decomposition 

$$ \lie{g}=\lie{h} + \lie{m} $$

\noindent where $\lie g$=Lie($G$), $\lie h$=Lie($H$) is the $+1$ eigenspace of 
the Cartan involution and $\lie{m}$ is the $-1$ eigenspace. This 
induces a decomposition  

\begin{equation}\label{eqtn:cartandecomp}
 \lie{g}^{\C}=\lie{h}^{\C} + \lie{m}^{\C}
\end{equation}

\noindent of $\lie g^{\C}$=Lie($G^{\C}$). In the case of
$G=\U(p,q)$, where $H=\U(p)\times\U(q)$ and thus $\lie{h}^{\C}=
\lie{gl}(p,\C)\oplus\lie{gl}(q,\C)$, the decomposition 
(\ref{eqtn:cartandecomp}) becomes

\begin{equation}\label{eqtn:cartandecomp-upq}
 \lie{gl}(n,\C)=(\lie{gl}(p,\C)\oplus\lie{gl}(q,\C)) 
 + \lie{m}^{\C}\ .
\end{equation}
If we identify $\lie{gl}(p,\C)\oplus\lie{gl}(q,\C)$ with the block 
diagonal elements in $\lie{gl}(n,\C)$, then $\lie{m}^{\C}$ 
corresponds to the off diagonal matrices. 

We can now describe the above Higgs bundles more intrinsically as 
follows. Let 
$P_{\GL(p,\C)}$ and $P_{\GL(q,\C)}$ be the principal frame bundles for $V$
and $W$ respectively.  Let $P=P_{\GL(p,\C)}\times P_{\GL(q,\C)}$ be 
the fiber product, and let $\Ad P=P\times_{\Ad}\mathfrak{gl}(n,\C)$ 
be the adjoint bundle, where 
$\GL(p,\C)\times\GL(q,\C)\subset\GL(n,\C)$ acts by the adjoint action 
on the Lie algebra of $\GL(n,\C)$.  This defines a subbundle 
\begin{equation}\label{eq:Ppq}
P_{\lie{m}^{\C}} = P\times_{\Ad}\lie{m}^{\C} \subset \Ad P\ .
\end{equation}
\noindent We can then make the following definition.
\begin{definition} 
\label{defn:upq}
A \emph{$\U(p,q)$-Higgs bundle\footnote{The reason for the name is 
explained by Remark~\ref{rem:HiggsG} and Lemma \ref{lemma:upq}}} on 
$X$ is a pair $(P,\Phi)$ where 
$P$ is a holomorphic principal $\GL(p,\C)\times\GL(q,\C)$ bundle, and 
$\Phi$ is a holomorphic section of the vector bundle 
$P_{\lie{m}^{\C}}\otimes K$ (where $P_{\lie{m}^{\C}}$ is the bundle
defined in (\ref{eq:Ppq})).
\end{definition}
\begin{remark} We can always write $P=P_{\GL(p,\C)}\times P_{\GL(q,\C)}$. 
If we let $V$ and $W$ be the standard vector bundles associated to 
$P_{\GL(p,\C)}$ and $P_{\GL(q,\C)}$ respectively, then any 
$\Phi\in H^0(P_{\lie{m}^{\C}}\otimes K)$ can be written as
in (\ref{upq-higgs-bundle}). We will usually adopt the vector bundle 
description of $\U(p,q)$-Higgs bundles. 
\end{remark} 
\begin{remark}\label{rem:HiggsG}  Definition \ref{defn:upq} is compatible
  with the definitions in \cite{hitchin:1992} and \cite{gothen:1995}, 
  where $G$-Higgs bundles
  are defined for any real form $G$ of a complex reductive Lie group
  $G^{\C}$. There, using the above notation,
  a $G$-Higgs bundle is a pair $(P,\Phi)$,
  where $P$ is a principal $H^{\C}$-bundle and $\Phi$ is a
  holomorphic section of $(P \times_{\Ad} \lie{m}^{\C}) \otimes K$.
  From a different perspective, Definition~\ref{defn:upq} defines an 
  example of a
  principal pair in the sense of \cite{banfield} and \cite{mundet}.
  Strictly speaking, since the canonical bundle $K$ plays the role of
  a fixed `twisting bundle', what we get is a principal pair in the
  sense of \cite{BGM}. The defining data for the pair are then the
  principal $\GL(p,\C)\times\GL(q,\C)\times\GL(1)$-bundle $P_{\GL(p,\C)}\times
  P_{\GL(q,\C)}\times P_K$ (where $P_K$ is the frame bundle for $K$), and
  the associated vector bundle $P_{\lie{m}^{\C}}\otimes K$.
\end{remark}
\begin{lemma}\label{lemma:upq} Let $(E = V\oplus W,\Phi)$ be a
  $\U(p,q)$-Higgs bundle with a hermitian metric such that $V\oplus W$
  is a unitary orthogonal decomposition.  Let $A$ be a unitary
  connection and let $D=d_A+\theta$ be the
  corresponding Higgs connection, where $\theta=\Phi+\Phi^*$. Then $D$
  is a $\U(p,q)$-connection, i.e.\ in any unitary local frame the
  connection 1-form takes its values in the Lie algebra of $\U(p,q)$.
\end{lemma}
\begin{proof} Fix a local unitary frame.  Then $D=d+A+\theta$, where
  $A$ takes its values in $\mathfrak{u}(p)\oplus
  \mathfrak{u}(q)\subset \mathfrak{u}(p,q)$,
  while $\theta$ takes its values in $\lie{m}$, where
  
  $$\lie{u}(p,q)=\mathfrak{u}(p)\oplus
  \mathfrak{u}(q)+\lie{m}$$
  
  \noindent is the eigenspace decomposition of the Cartan involution.
 \end{proof}
 
\begin{definition}\label{defn:U(pq)Higgs}
Let $(E,\Phi)$ be a $\U(p,q)$-Higgs bundle
with $E=V\oplus W$ and $\Phi=\left(
  \begin{smallmatrix}
    0 & \beta \\
    \gamma & 0
  \end{smallmatrix}\right)$. 
We say $(E,\Phi)$ is a \emph{stable} $\U(p,q)$-Higgs bundle if the  
slope stability condition 
$\mu(E') < \mu(E)$, is satisfied for all 
$\Phi$-invariant subbundles of the form $E'=V'\oplus W'$, i.e.\ for all
subbundles $V'\subset V$ and 
$W'\subset W$ such that 
\begin{align}
\beta &:W'\longrightarrow V'\otimes K\\
\gamma &:V'\longrightarrow W'\otimes K\ .
\end{align}
\emph{Semistability} for $\U(p,q)$-Higgs bundles is defined by
replacing the above strict inequality with a weak inequality, and
$(E,\Phi)$ is \emph{polystable} if it is a direct sum of stable
$\U(p,q)$-Higgs bundles all of the same slope.  We shall say that a
polystable $\U(p,q)$-Higgs bundle which is not stable is
\emph{reducible}.  A \emph{morphism} between two $\U(p,q)$-Higgs
bundles $(V\oplus W,\Phi)$ and $(V'\oplus W',\Phi')$ is given by maps
$g_V:V\rightarrow V'$ and $g_W:W\rightarrow W'$ which intertwine
$\Phi$ and $\Phi'$, i.e.\ such that $(g_V\oplus g_W)\otimes
I_K\circ\Phi=\Phi'\circ(g_V\oplus g_W)$ where $I_K$ is the identity on
$K$.  In particular we have a natural notion of isomorphism of
$\U(p,q)$-Higgs bundles.
\end{definition}

\begin{remark}\label{rem:stability-equiv} The  stability condition for
  a $\U(p,q)$-Higgs bundle is \emph{a priori} weaker than the stability
  condition given in Definition \ref{defn:GLnHiggs} for $\GL(n,\C)$-Higgs
  bundles.  However, it is shown in \cite[Section 2.3]{gothen:2001} that the
  weaker condition is in fact equivalent to the ordinary stability of
  $(E,\Phi)$.
\end{remark}
\begin{proposition}
 \label{prop:U(pq)HK}
Let $(E,\Phi)$ be a $\U(p,q)$-Higgs bundle with $E=V\oplus W$ and 
$\Phi=\left( 
  \begin{smallmatrix}
    0 & \beta \\
    \gamma & 0
  \end{smallmatrix}\right)$. Then $(E,\Phi)$ is polystable 
if and only if it admits a hermitian metric such that $E = V \oplus W$
is an orthogonal decomposition and such that Hitchin's equations
(\ref{eq:hitchin1}) are satisfied.
\end{proposition}
\begin{proof} This is a special case of
  the correspondence invoked in \cite{hitchin:1992} for $G$-Higgs
  bundles where $G$ is a real form of a reductive Lie group. By Remark
  \ref{rem:HiggsG} it can also be seen as a special case of the
  Hitchin--Kobayashi correspondence for principal pairs (cf.\ 
  \cite{banfield} and \cite{mundet} and \cite{BGM}).  We note finally
  that in one direction the result follows immediately from Theorem
  \ref{thm:HK} (1): if $(V\oplus W, \Phi)$ supports a
  compatible metric such that (\ref{eq:hitchin1}) is satisfied, then
  it is polystable as a $\GL(n,\C)$-Higgs bundle, and hence it is
  $\U(p,q)$-polystable.
\end{proof}

\begin{definition}\label{defn:M(a,b)}
We define $\mathcal{M}(a,b)$ to be the moduli space of polystable 
$\U(p,q)$-Higgs bundles with $\deg(V)=a$ and $\deg W=b$. We denote by
$\mathcal{M}^s(a,b)$ the subspace parameterizing the
 strictly stable $\U(p,q)$-Higgs bundles.
\end{definition}

  The construction of $\mathcal{M}(a,b)$ is essentially the same as in
  section~\S 9 of \cite{simpson:1994b}. There the moduli space of
  $G$-Higgs bundles is constructed for any reductive group $G$. We
  take $G=\GL(p,\C)\times \GL(q,\C)$. The difference between a
  $\U(p,q)$-Higgs bundle and a $\GL(p,\C)\times \GL(q,\C)$-Higgs bundle is
  entirely in the nature of the Higgs fields. Taking the standard
  embedding of $\GL(p,\C)\times \GL(q,\C)$ in $\GL(p+q,\C)$ we see that in a
  $\GL(p,\C)\times \GL(q,\C)$-Higgs bundle the Higgs field $\Phi$ takes its
  values in the subspace
  $(\mathfrak{gl}(p)\oplus\mathfrak{gl}(q))\subset\mathfrak{gl}(p+q)$,
  while in a $\U(p,q)$-Higgs bundle the Higgs field $\Phi$ takes its
  values in the complementary subspace
  $\lie{m}^{\C}$ (as in (\ref{eqtn:cartandecomp-upq})). Since both
  subspaces are invariant under the adjoint action of
  $\GL(p,\C)\times\GL(q,\C)$, the same method of construction works for the
  moduli spaces of both types of Higgs bundle. 
 
We can describe the gauge theory version of the moduli space 
$\mathcal{M}(a,b)$ using standard methods; see Hitchin \cite{hitchin:1987} 
for a construction in the case of  ordinary rank $2$ Higgs bundles. 
To adapt to our case we proceed as follows. Let $E = V \oplus W$ be a 
smooth complex vector bundle with a hermitian metric such that the 
direct sum decomposition is orthogonal.  We let $\mathcal{A}$ denote 
the space of connections on $E$ which are direct sums of unitary 
connections on 
$V$ and $W$ and we let $\boldsymbol{\Omega}$ denote the space of 
Higgs fields 
$\Phi \in \Omega^{1,0}(\End(E))$ of the form 
$\Phi = \left(
    \begin{smallmatrix}
      0 & \beta \\
      \gamma & 0
    \end{smallmatrix} \right)$.
  The correspondence between unitary connections and holomorphic
  structures via $\dbar$-operators turns $\mathcal{A} \times
  \boldsymbol{\Omega}$ into a complex affine space which acquires a
  hermitian metric using the metric on $E$ and integration over $X$.
  The group $\mathcal{G}$ of $\U(p)\times \U(q)$-gauge transformations
  acts on the configuration space $\mathcal{C} \subseteq \mathcal{A}
  \times \boldsymbol{\Omega}$ of solutions $(A,\Phi)$ to Hitchin's
  equations \eqref{eq:hitchin1}.  The quotient
  $\mathcal{C}/\mathcal{G}$ is, by definition, the gauge theory moduli space.
  As in \cite{hitchin:1987}, the open subset of $\mathcal{C}/\mathcal{G}$ 
  corresponding to irreducible solutions has a K\"ahler manifold structure 
  
  To see that the gauge theory moduli space is homeomorphic to
  $\mathcal{M}(a,b)$ we can consider this latter space from the
  complex analytic point of view (cf.\ Remark~\ref{rem:dolbeault}
  below): consider triples $(\dbar_V,\dbar_W,\Phi)$, where $\dbar_V$
  and $\dbar_W$ are $\dbar$-operators on $V$ and $W$, respectively,
  and $\Phi\in\boldsymbol{\Omega}$.  Let $\mathcal{C}_{\C}$ be the set
  of such triples for which $\Phi$ is holomorphic and the associated
  $\U(p,q)$-Higgs bundle is polystable.  We can then view
  $\mathcal{M}(a,b)$ as the quotient of $\mathcal{C}_{\C}$ by the
  complex gauge group.  We clearly have an inclusion $\mathcal{C}
  \into \mathcal{C}_{\C}$ which descends to give a continuous map from
  the gauge theory moduli space to $\mathcal{M}(a,b)$.  The
  Hitchin-Kobayashi correspondence of Proposition~\ref{prop:U(pq)HK}
  now shows that this map is in fact a homeomorphism.
  
  For a third perspective, we observe that provided that $V$ and $W$
  are not isomorphic bundles, i.e.\ provided $p\ne q$ or $a\ne b$, we
  can view $\mathcal{M}^s(a,b)$ as a subvariety of a moduli space of
  stable $\GL(p+q)$-Higgs bundle. If $V\simeq W$, then
  $\mathcal{M}^s(a,b)$ is a finite cover of a subvariety in the larger
  moduli space:

\begin{proposition}With $n=p+q$ and $d=a+b$, let $\mathcal{M}^s(d)$ denote the 
  moduli space of stable $\GL(n,\C)$-Higgs bundles of degree $d$. If
  $p\ne q$ or $a\ne b$ then $\mathcal{M}^s(a,b)$ embeds as a closed
  subvariety in $\mathcal{M}^s(d)$. If $p=q$ and $a=b$, then there is
  an involution on $\mathcal{M}^s(a,a)$ such that the quotient injects
  into $\mathcal{M}^s(d)$.
\end{proposition}

\begin{proof}Let $[V\oplus W,\Phi]_{p,q}$ denote the point in 
  $\mathcal{M}^s(a,b)$ represented by the $\U(p,q)$-Higgs bundle
  $(V\oplus W,\Phi)$. Then $(E=V\oplus W,\Phi)$ is a stable
  $\GL(n,\C)$-Higgs bundle and the map
  $\mathcal{M}^s(a,b)\rightarrow\mathcal{M}(d)$ is defined by
$$[V\oplus W,\Phi]_{p,q}\mapsto [E,\Phi]_n\ ,$$
\noindent where $[,]_n$ denotes the isomorphism class in 
$\mathcal{M}(d)$. The only question is whether this map is injective.
Suppose that $(E=V\oplus W,\Phi)$ and $(E'=V'\oplus W',\Phi')$ are 
isomorphic as $\GL(n,\C)$-Higgs bundles. Let the isomorphism be given 
by a complex gauge transformation $g:E\rightarrow E'$. If $g$ is not 
of the form $(\begin{smallmatrix} 
  g_V & 0 \\
  0 & g_W
  \end{smallmatrix})$
then the off diagonal components determine morphisms 
$\xi:V\rightarrow W'$ and $\sigma:W\rightarrow V'$. 
Let $N=\ker(\xi)\oplus\ker (\sigma)$ be the subbundle of $V\oplus W$ 
determined by the kernels of $\xi$ and $\sigma$. If $p\ne q$ then 
$N$ is a non-trivial proper subbundle. Moreover, using the fact that 
$g\Phi=\Phi' g$, we see that it is $\Phi$-invariant. Since $(V\oplus 
W,\Phi)$ is stable, it follows that 

\begin{equation}\label{eqtn:mu(K)}
\mu(N)<\mu(E)\ .
\end{equation}

Similarly, the images of $\xi$ and $\sigma$ determine a proper 
$\Phi'$-invariant subbundle of $E'$, say $I$, for which

\begin{equation}\label{eqtn:mu(I)}
\mu(I)<\mu(E')\ .
\end{equation}

But if $\mu(E)=\mu(E')$ then (\ref{eqtn:mu(K)}) and  
(\ref{eqtn:mu(I)}) cannot both be satisfied. Thus $\xi$ and $\sigma$ 
must both vanish and hence $[V\oplus W,\Phi]_{p,q}=[V'\oplus 
W',\Phi']_{p,q}$. 

If $p=q$, then this argument can fail, but only if $\xi$ and 
$\sigma$ are both isomorphisms. In that case, $N=0$ and $I=E$.  This 
also requires $a=b$. Under these conditions, if $V$ and $W$ are 
non-isomorphic, then 
$[V\oplus W 
  ,(\begin{smallmatrix}
    0 & \beta \\
    \gamma & 0
  \end{smallmatrix})]_n=[W\oplus V 
  ,(\begin{smallmatrix}
    0 & \gamma \\
    \beta & 0
  \end{smallmatrix})]_n$ but the Higgs bundles are not isomorphic as 
  $\U(p,q)$-Higgs bundles.  Hence the last statement of the
  Proposition follows taking the involution
  $[V\oplus W 
  ,(\begin{smallmatrix}
    0 & \beta \\
    \gamma & 0
  \end{smallmatrix})] \mapsto [W\oplus V 
  ,(\begin{smallmatrix}
    0 & \gamma \\
    \beta & 0
  \end{smallmatrix})]$
  on $\mathcal{M}^s(a,a)$.
\end{proof}
 
 \begin{proposition}\label{prop:Ms=M}
If $\GCD(p+q,a+b)=1$ then $\mathcal{M}^s(a,b)=\mathcal{M}(a,b)$.
\end{proposition}
\begin{proof} If $\GCD(p+q,a+b)=1$ then for 
  purely numerical reasons there are no strictly semistable
  $\U(p,q)$-Higgs bundles in $\mathcal{M}(a,b)$. 
\end{proof}

The link to moduli spaces of representations is provided by the next 
result.

 \begin{proposition}\label{prop:R=M}
  There is a homeomorphism
    $\mathcal{M}(a,b)\cong \mathcal{R}_\Gamma(a,b)$. 
 \end{proposition}

\begin{proof}
  Suppose that $(E=V\oplus W,\Phi)$ represents a point in
  $\mathcal{M}(a,b)$, i.e.\ suppose that it is a $\U(p,q)$-polystable
  Higgs bundle, and suppose that $E$ has a hermitian metric such that
  the direct sum decomposition is orthogonal and Hitchin's equations
  (\ref{eq:hitchin1}) are satisfied. Rewriting the equations in terms
  of the Higgs connection $D=d_A+\theta$, where $A$ is the metric
  connection and $\theta=\Phi+\Phi^*$, we see that
  $D$ is projectively flat. By Lemma \ref{lemma:upq} it is a
  projectively flat $\U(p,q)$-connection, and thus defines a point in
  $\mathcal{R}_\Gamma(a,b)$. Conversely by Corlette's theorem
  \cite{corlette:1988}, every representation in $\Hom^+(\pi_1
  X,\PU(p,q))$, or equivalently every representation in
  $\Hom^+(\Gamma,\U(p,q))$, arises in this way.  The fact that this
  correspondence gives a homeomorphism follows by the same argument as
  the one given in \cite{simpson:1994b} for ordinary Higgs bundles.
\end{proof}

\begin{definition}\label{defn:R*}
 Define the subspace 
$\mathcal{R}^*_\Gamma(a,b)$ to be the subspace corresponding 
to $\mathcal{M}^s(a,b)$ via the homeomorphism in Proposition
\ref{prop:R=M}.  Using the fibration of $\mathcal{R}_\Gamma(a,b)$ 
over $\mathcal{R}[a,b]$, define 
$\mathcal{R}^*[a,b]\subset\mathcal{R}[a,b]$ to 
be the image of $\mathcal{R}^*_\Gamma(a,b)$.
\end{definition}

%
% New version of Remark 6.3
%
\begin{remark}\label{rem:R*}Thus 
$\mathcal{R}^*_\Gamma(a,b)$ parameterizes the representations which 
give rise to stable 
$\U(p,q)$-Higgs bundles.  Recall from Remark~\ref{rem:stability-equiv} 
that a $\U(p,q)$-Higgs bundle is stable (in the sense of 
Definition~\ref{defn:U(pq)Higgs}) if  and only if its is stable as an 
ordinary $\GL(n,\C)$-Higgs bundle. Now, a $\GL(n,\C)$-Higgs bundle is 
stable if and only if the  corresponding representation of $\Gamma$ 
on $\C^n$ is irreducible (cf.\ Corlette \cite{corlette:1988}).  Hence 
we see that the 
  subspace $\mathcal{R}^*_\Gamma(a,b)$ corresponds to the
  representations of $\Gamma$ in $\U(p,q)$ which are irreducible as
  $\GL(n,\C)$ representations.  Similarly, the subspace
  $\mathcal{R}^*[a,b]$ corresponds to the representations of $\pi_1 X$
  which are irreducible as $\PGL(n,\C)$ representations.
  
We point out, moreover, that the subspace $\mathcal{R}^*_\Gamma(a,b)$ 
includes as a dense open set the representations whose induced 
adjoint representations on the Lie algebra of $\PU(p,q)$ are 
irreducible. It may also contain some representations whose induced 
adjoint representation is reducible for the following reason. If 
$(E=V\oplus W,\Phi)$ is the 
$\U(p,q)$-Higgs bundle corresponding to a representation in 
$\mathcal{R}^*_\Gamma(a,b)$, then $(\End(E),\Phi)$ is a polystable
Higgs bundle but it is not necessarily stable. The representations 
with reducible induced adjoint representation are the ones for which 
$(\End(E),\Phi)$ is strictly polystable. 
\end{remark}
%
% Old version of Remark 6.3
%
% \begin{remark}\label{rem:R*}
%   
% \end{remark}

%%%%%%%%%%%%%%%%%%%%%%%%%%%%%%%%%%%%%%%%%%%%%%%%
\subsection{Deformation theory}
\label{sec:higgs-deformation-theory}
%%%%%%%%%%%%%%%%%%%%%%%%%%%%%%%%%%%%%%%%%%%%%%%%

The results of Biswas and Ramanan 
\cite{biswas-ramanan:1994} and Hitchin \cite{hitchin:1992} readily adapt 
to describe the deformation theory of $\U(p,q)$-Higgs bundles.

\begin{definition}
\label{def:u-notation}
Let $(E = V \oplus W, \Phi)$ be a
$\U(p,q)$-Higgs bundle.  We introduce the following notation:
  \begin{align*}
  U & = \End(E) \\ 
U^+ &= \End(V) \oplus \End(W)\ , \\ 
U^- &= \Hom(W,V)  \oplus \Hom(V,W)\ . 
\end{align*}
\end{definition}

With this notation, $U=U^+\oplus U^-$, $\Phi \in H^0(U^-\otimes 
K)$, and $\ad(\Phi)$ interchanges $U^+$ and 
$U^-$. We consider the complex of sheaves
\begin{equation}
  \label{eq:tangentspace}
  C^{\bullet} : 
U^+  \xrightarrow{\ad(\Phi)}  U^- \otimes K\ .  
\end{equation}

\begin{lemma}
 \label{lem:stable-upq-vanishing}
  Let $(E,\Phi)$ be a stable $\U(p,q)$-Higgs bundle.  Then
\begin{align}
  \label{eq:ker-phi-u}
  \ker\bigl(\ad(\Phi) \colon H^0(U^+) \to
  H^0(U^-\otimes K)\bigr) &= \C\ , \\
  \label{eq:ker-phi-u-perp}
  \ker\bigl(\ad(\Phi) \colon H^0(U^-) \to
  H^0(U^+\otimes K)\bigr) &= 0\ . 
\end{align}
\end{lemma}
\begin{proof}
  By Remark~\ref{rem:stability-equiv} $(E,\Phi)$ is stable as a
  $\GL(n,\C)$-Higgs bundle.  Hence it is simple, that is, its only
  endomorphisms are the non-zero scalars.  Thus,
$$
\ker\bigl(\ad(\Phi) \colon H^0(U) \to
H^0(U\otimes K) \bigr) = \C\ . 
$$ 
Since 
$U = U^+ \oplus U^-$ and
$\ad(\Phi)$ interchanges these two summands, the statements of the Lemma
follow. 
\end{proof}

\begin{proposition}[Biswas-Ramanan \cite{biswas-ramanan:1994}]
  \label{prop:upq-deformation} $\mathrm{}$
\begin{itemize}
\item[$(1)$] The space of endomorphisms of $(E,\Phi)$ is isomorphic to
  the zeroth hypercohomology group $\HH^0(C^\bullet)$.
\item[$(2)$] The space of infinitesimal deformations of $(E,\Phi)$ is
  isomorphic to the first hypercohomology group $\HH^1(C^\bullet)$.
\item[$(3)$] There is a long exact sequence
\begin{multline}  \label{eq:long-exact-tangent}
  0 \lto  \mathbb{H}^0(C^{\bullet}) \lto H^0(U^+) \lto
  H^0(U^-\otimes K) \lto  \mathbb{H}^1(C^{\bullet})  \\ 
  \lto  H^1(U^+) \lto  H^1(U^-\otimes  K) \lto   
\mathbb{H}^2(C^{\bullet})\lto 0\ ,
\end{multline}
where the maps $H^i(U^+) \lto  H^i(U^-\otimes  K)$ are induced by
$\ad(\Phi)$. 
\end{itemize}
\hfill\qed
\end{proposition}

\begin{proposition}
\label{prop:stable-upq-vanishing}
Let $(E,\Phi)$ be a stable $\U(p,q)$-Higgs bundle, then
\begin{enumerate}
\item[$(1)$] $\mathbb{H}^0(C^{\bullet}) = \C$ (in other words $(E,\Phi)$ is 
simple) and
\item[$(2)$] $\mathbb{H}^2(C^{\bullet}) = 0$. 
\end{enumerate}
\end{proposition}
\begin{proof}
  (1) Follows immediately from
  Lemma~\ref{lem:stable-upq-vanishing} and
  Proposition~\ref{prop:upq-deformation} (3).
  
  (2) We have natural
  $\ad$-invariant isomorphisms $U^+ \cong (U^+)^*$ and $U^- \cong
  (U^-)^*$.  Thus
  $$
  \ad(\Phi) \colon H^1(U^+) \to H^1(U^-\otimes K)
  $$
  is Serre dual to $\ad(\Phi) \colon H^0(U^-) \to H^0(U^+\otimes
  K)$.  Hence Lemma~\ref{lem:stable-upq-vanishing} and $(3)$ of
  Proposition~\ref{prop:upq-deformation} show that
  $\mathbb{H}^2(C^{\bullet}) = 0$.
\end{proof}

\begin{proposition}\label{prop:smoothness-higgs}
The moduli space  of  stable  $\U(p,q)$-Higgs bundles
is a  smooth complex variety of dimension $1+(p+q)^2 (g-1)$.
\end{proposition}
\begin{proof}
  By Proposition \ref{prop:stable-upq-vanishing} (2)
 $\mathbb{H}^2(C^{\bullet}) = 0$ at all points
  in the moduli space  of  stable  $\U(p,q)$-Higgs bundles.
  Smoothness is thus a consequence of the results of
  \cite{biswas-ramanan:1994}, as follows.  Let $e \in
  \mathcal{M}(a,b)$ be the point corresponding to a stable
  $\U(p,q)$-Higgs bundle $(E,\Phi)$ and let $\mathcal{F}$ be the
  infinitesimal deformation functor of $(E,\Phi)$ as in
  \cite{biswas-ramanan:1994}.  Then the completion of the local ring
  $\mathcal{O}_e$ pro-represents $\mathcal{F}$ (cf.\ Schlessinger
  \cite{schlessinger:1968}).  Now
  Proposition~\ref{prop:stable-upq-vanishing} and Theorem~3.1 of
  \cite{biswas-ramanan:1994} show that the completion of
  $\mathcal{O}_e$ is regular and hence $\mathcal{O}_e$ is itself
  regular.  Thus $\mathcal{M}(a,b)$ is smooth at $e$.

Using $(2)$ and $(3)$ of Proposition~\ref{prop:upq-deformation},
Proposition~\ref{prop:stable-upq-vanishing} and the Riemann-Roch
Theorem, the dimension of the moduli space is given by
\begin{align*}
\dim \mathbb{H}^1(C^{\bullet}) &= 1-\chi(U^+)+\chi(U^-\otimes K)\\
                               &= 1+(p^2+q^2)(g-1)+ 2pq (g-1)\\
                               &= 1+(p+q)^2 (g-1)\ .
\end{align*}
\end{proof} 
\begin{remark} The dimension of the moduli space of stable $\U(p,q)$-Higgs
bundles is half that of the moduli space of stable 
$\GL(p+q,\C)$-Higgs bundles.
\end{remark} 
\begin{remark}By Proposition~\ref{prop:Ms=M} $\mathcal{M}(a,b)$ is smooth if
$\GCD(p+q,a+b)=1$.
\end{remark}

\begin{remark}\label{rem:dolbeault}
  As an alternative to the algebraic arguments of
  \cite{biswas-ramanan:1994}, the fact that the deformation theory of
  a $\U(p,q)$-Higgs bundle is controlled by the complex of sheaves
  \eqref{eq:tangentspace} can be seen from the complex analytic point
  of view as follows.  As in the gauge theory construction of
  $\mathcal{M}(a,b)$ (cf.\ Section~\ref{subs:upq} let $V \oplus W$ be a
  smooth complex vector bundle, and consider a $\U(p,q)$-Higgs bundle
  as being given by a triple $(\dbar_V,\dbar_W,\Phi)$.  Now write down
  a Dolbeault resolution of the complex $C^{\bullet}$:
  \begin{displaymath}
    \begin{CD}
      \Omega^{0}(U^{+}) @>\ad(\Phi)>> \Omega^{1,0}(U^{-}) \\
      @VV\dbar V  @VV-\dbar V \\
      \Omega^{0,1}(U^{+}) @>\ad(\Phi)>> \Omega^{1,1}(U^{-}) \\
      @VVV  @VVV \\
      0 @>>> 0
    \end{CD}\ .
  \end{displaymath}
  Consider the associated total complex $\mathbf{C}^{0}
  \overset{D^0}{\to} \mathbf{C}^{1} \overset{D^1}{\to}
  \mathbf{C}^{2}$.  Then  $\mathbf{C}^{0}$ is
  the Lie algebra of the $\GL(p,\C)\times\GL(q,\C)$-gauge group and
  $\mathbf{C}^{1}$ is the tangent space to the affine space of
  triples $(\dbar_V,\dbar_W,\Phi)$.  Furthermore, $D^0$ is the
  infinitesimal action of the complex gauge group, while $D^1$ is the
  derivative of the holomorphicity condition: this gives the desired
  interpretation of the deformation complex $C^{\bullet}$ in complex
  analytic terms.
  
  To conclude this line of thought we give an alternative argument
  for the smoothness of the moduli space of stable $\U(p,q)$-Higgs
  bundles: suppose that
  $(\dbar_V,\dbar_W,\Phi)$ corresponds to a stable $\U(p,q)$-Higgs
  bundle $(E,\Phi)$.  Proposition~\ref{prop:stable-upq-vanishing}
  shows that $\mathbb{H}^0(C^{\bullet}) = \C$ and
  $\mathbb{H}^{2}(C^{\bullet})=0$.  The differential of the
  holomorphicity condition is thus surjective and $(E,\Phi)$ has no non-trivial
  automorphisms.  It follows by standard arguments that the
  moduli space can be constructed as a smooth complex manifold near
  $(E,\Phi)$.
\end{remark}

%%%%%%%%%%%%%%%%%%%%%%%%%%%%%%%%%%%%%%%%%%%%%%%%%
\subsection{Bounds on the topological invariants}
\label{sec:topological-bounds}
%%%%%%%%%%%%%%%%%%%%%%%%%%%%%%%%%%%%%%%%%%%%%%%%%
In this section we show how the Higgs bundle point of view provides 
an easy proof of a result of Domic and Toledo 
\cite{domic-toledo:1987} which allows us to bound the topological
invariants $\deg(V)$ and $\deg(W)$ for which $\U(p,q)$-Higgs bundles
may exist.  The lemma is a slight variation
on the results of \cite[Section 3]{gothen:2001} (cf.\ also Lemma 3.6
of Markman and Xia \cite{markman-xia:2001}).
\begin{lemma}
  \label{lem:slope-bound}
  Let $(E,\Phi)$ be a semistable $\U(p,q)$-Higgs bundle.  Then
  \begin{align}
    \label{eq:muE1-mu}
    p(\mu(V) - \mu(E)) &\leq \rk(\gamma)(g-1), \\
    \label{eq:muE2-mu}
    q(\mu(W) - \mu(E)) &\leq \rk(\beta)(g-1).
  \end{align}
  If equality occurs in \eqref{eq:muE1-mu} then either $(E,\Phi)$ is
  strictly semistable or $p = q$ and $\gamma$ is an isomorphism. 
  If equality occurs in \eqref{eq:muE2-mu} then either $(E,\Phi)$ is
  strictly semistable or $p = q$ and $\beta$ is an isomorphism. 
\end{lemma}
\begin{proof}
  If $\gamma = 0$ then $V$ is $\Phi$-invariant. By stability,
  $\mu(V) \leq \mu(E)$ and equality can only occur if 
  $(E,\Phi)$ is strictly semistable.  This proves \eqref{eq:muE1-mu} 
  in the case $\gamma = 0$. We may therefore assume that $\gamma 
  \neq 0$.   Let $N = \ker(\gamma) \subseteq V$ and let $I = \im(\gamma)
  \otimes K^{-1} \subseteq W$.  Then
  \begin{equation}
    \rk(N) + \rk(I) = p
    \label{eq:rkG+rkH}
  \end{equation}
  and, since $\gamma$ induces a non-zero section of $\det((V/N)^* \otimes
  I \otimes K)$,
  \begin{equation}
    \deg(N) + \deg(I) + \rk(I)(2g-2) \geq \deg(V).
    \label{eq:degG+degH-1}
  \end{equation}
  The bundles $N$ and $V \oplus I$ are $\Phi$-invariant subbundles 
  of $E$ and hence we obtain by semistability that $\mu(N) \leq \mu(E)$ 
  and $\mu(V \oplus I) \leq \mu(E)$ or, equivalently, that
  \begin{align}
    \label{eq:degG}
    \deg(N) &\leq \mu(E) \rk(N), \\
    \label{eq:degH}
    \deg(I) &\leq \mu(E)(p + \rk(I)) - \deg(V).
  \end{align}
  Adding \eqref{eq:degG} and \eqref{eq:degH} and using 
  \eqref{eq:rkG+rkH} we obtain
  \begin{equation}
    \deg(N) + \deg(I) \leq 2\mu(E)p - \deg(V).
    \label{eq:degG+degH-2}
  \end{equation}
  Finally, combining \eqref{eq:degG+degH-1} and 
  \eqref{eq:degG+degH-2} we get 
  \begin{displaymath}
    \deg(V) - \rk(I)(2g-2) \leq 2 \mu(E)p - \deg(V),
  \end{displaymath}
  which is equivalent to \eqref{eq:muE1-mu} since $\rk(\gamma) = \rk(I)$.  
  Note that equality can only occur if we have equality in 
  \eqref{eq:degG} and \eqref{eq:degH} and thus either $(E,\Phi)$ is 
  strictly semistable or neither of the subbundles $N$ and $V \oplus 
  I$ is proper and non-zero.  In the latter case, clearly $N = 0$ and 
  $I = W$ and therefore $p = q$; furthermore we must also 
  have equality in \eqref{eq:degG+degH-1} implying that $\gamma$ is an 
  isomorphism.
  An analogous argument applied to $\beta$ proves \eqref{eq:muE2-mu}.
\end{proof}
\begin{remark}
  The proof also shows that if we have equality in, say,
  \eqref{eq:muE1-mu} then $\gamma \colon V/N \to I\otimes K$ is an
  isomorphism.  In particular, if $p <q$ and $\mu(V) -
  \mu(E) = g-1$ then $\gamma \colon V \xrightarrow{\cong} I \otimes K$.
\end{remark}
We can re-formulate Lemma~\ref{lem:slope-bound} to obtain the 
following corollary. 
\begin{corollary}
  \label{cor:mu-muE_i}
  Let $(E,\Phi)$ be a semistable $\U(p,q)$-Higgs bundle.  Then
  \begin{align}
    \label{eq:mu-mu_2}
    q(\mu(E) - \mu(W)) &\leq \rk(\gamma)(g-1), \\
    \label{eq:mu-mu_1}
    p(\mu(E) - \mu(V)) &\leq \rk(\beta)(g-1).
  \end{align}
\end{corollary}
\begin{proof}
 Use $\mu(W) - \mu(E) =
  \frac{p}{q}\bigl( \mu(E) - \mu(V)\bigr)$ to see that 
  \eqref{eq:mu-mu_2} is equivalent to \eqref{eq:muE1-mu}.  Similarly
  \eqref{eq:mu-mu_1} is equivalent to \eqref{eq:muE2-mu}.
\end{proof} 
An important corollary of the lemma above is the 
following Milnor--Wood type inequality for $\U(p,q)$-Higgs bundles (due
to Domic and Toledo \cite{domic-toledo:1987}, improving on a bound
obtained by Dupont \cite{dupont:1978} in the case $G = \SU(p,q)$). 
This result gives bounds on the possible values of the topological
invariants $\deg(V)$ and $\deg(W)$.
\begin{corollary}
  \label{cor:toledo}
  Let $(E,\Phi)$ be a semistable $\U(p,q)$-Higgs bundle.  Then
  \begin{equation}
    \label{eq:toledo}
    \frac{pq}{p+q} \abs{\mu(V) - \mu(W)} 
    \leq \min\{p,q\}(g-1).
  \end{equation}
\end{corollary}
\begin{proof}
  Since $\mu(E) = \frac{p}{p+q}\mu(V) + \frac{q}{p+q}\mu(W)$ 
  we have $\mu(V) - \mu(E) = \frac{q}{p+q}(\mu(V) - 
  \mu(W)$ and therefore \eqref{eq:muE1-mu} gives
  \begin{displaymath}
    \frac{pq}{p+q} (\mu(V) - \mu(W)) \leq 
    \rk(\gamma)(g-1).
  \end{displaymath}
  A similar argument using \eqref{eq:muE2-mu} shows that 
  \begin{displaymath}
    \frac{pq}{p+q} (\mu(W) - \mu(V)) \leq 
    \rk(\beta)(g-1).
  \end{displaymath}
  But, obviously, $\rk(\beta)$ and $\rk(\gamma)$ are both less than or equal to 
  $\min\{p,q\}$.
\end{proof}

\begin{definition}\label{defn:toledo} The \emph{Toledo invariant} of 
the representation corresponding to $(E=V\oplus W,\Phi)$ is 
  \begin{equation}
\tau=\tau(a,b)    = 2\frac{qa - pb}{p + q}
\label{toledo-invariant}
  \end{equation}
\noindent where $a=\deg(V)$  
and 
$b=\deg(W)$.
 \end{definition}

\begin{remark}
\label{rem:tau-M}
Since
$$
    \tau = 2\frac{pq}{p+q} (\mu(V) - \mu(W))
    = -2p(\mu(E) - \mu(V))=2q(\mu(E) - \mu(W))\ ,
$$
the inequalities in Lemma~\ref{lem:slope-bound} and
Corollary~\ref{cor:mu-muE_i} can be written as 
  \begin{align}
    \label{eq:rkgamma}
   \frac{\tau}{2} &\leq \rk(\gamma)(g-1), \\
    \label{eq:rkbeta}
  - \frac{\tau}{2}&\leq \rk(\beta)(g-1).
  \end{align}
\noindent Similarly the inequality
\eqref{eq:toledo} can be written 
$\abs{\tau} \leq \tau_M$, where
\begin{equation}\label{eqn:tau-M}
\tau_M=\min\{p,q\}(2g-2)\ .
\end{equation}
\end{remark}

%%%%%%%%%%%%%%%%%%%%%%%%%%%%%%%%%%%%%%%%%%%%%%%%%%%%%%%%%%%%%%%%%
\subsection{Rigidity and extreme values of the Toledo 
invariant}\label{subs:rigidity}
%%%%%%%%%%%%%%%%%%%%%%%%%%%%%%%%%%%%%%%%%%%%%%%%%%%%%%%%%%%%%%%%%

If $|\tau|=\tau_M$ then the moduli space $\mathcal{M}(a,b)$ has 
special features. These depend on whether $p=q$ or $p\ne q$. 

Consider first the case $p = q$.  Notice that if $p=q$ then 
$\tau(a,b)=a-b$ and $\tau_M=p(2g-2)$. We thus examine the moduli space 
$\mathcal{M}(a,b)$ when $|a-b|=p(2g-2)$. Before giving a description
we review briefly the notion of 
$L$-twisted Higgs pairs. Let $L$ be a line bundle. An $L$-twisted 
Higgs pair $(V,\theta)$ consists of a holomorphic vector bundle 
$V$ and an $L$-twisted homomorphism $\theta:V\lto V\otimes L$.
The notions of stability, semistability and polystability are defined 
as for Higgs bundles. The moduli space of semistable $L$-twisted 
Higgs pairs has been constructed by Nitsure using Geometric Invariant 
Theory \cite{nitsure:1991}. Let $\mathcal{M}_L(n,d)$ be the moduli
space of polystable $L$-twisted Higgs pairs of rank $n$ and degree 
$d$. 

\begin{proposition} \label{prop:p=q-toledo-max}
Let $p=q$ and $|a-b|=p(2g-2)$. Then 
$$
\mathcal{M}(a,b)\cong \mathcal{M}_{K^2}(p,a)\cong \mathcal{M}_{K^2}(p,b).
$$
\end{proposition}

\begin{proof}
Let $(E=V\oplus W, \Phi)\in \mathcal{M}(a,b)$. Suppose for 
definiteness  that $b-a=p(2g-2)$. {}From (\ref{eq:muE2-mu}) it 
follows that $\gamma: V\lto W\otimes K$ is an isomorphism. We can 
then compose  $\beta: W\lto V\otimes K$ with 
$\gamma\otimes \Id_K: V\otimes K\lto W\otimes K^2$
to obtain a $K^2$-twisted Higgs pair 
$\theta_W: W\lto W\otimes K^2$. 
Similarly, twisting $\beta: W\lto V\otimes K$ with $K$ and composing 
with $\gamma$, we obtain a $K^2$-twisted Higgs pair 
$\theta_V: V\lto V\otimes K^2$. Conversely, given an isomorphism
$\gamma: V\lto W\otimes K$, we can recover $\beta$ from $\theta_V$
as well as from $\theta_W$. It is clear that the (poly)stability of 
$(E,\Phi)$ is equivalent to the (poly)stability of $(V,\theta_V)$ and 
to the (poly)stability of 
$(W,\theta_W)$, proving the claim.
\end{proof}

\begin{remark}
  The moduli space $\mathcal{M}_{K^2}(p,a)$ contains an open
  (irreducible) subset consisting of a vector bundle over the moduli
  space of stable bundles of rank $p$ and degree $a$. This is because
  the stability of $V$ implies the stability of any $K^2$-twisted
  Higgs pair $(V,\theta_V)$, and $H^1(\End V\otimes K^2)=0$. The rank
  of the bundle is determined by the Riemann--Roch Theorem.
\end{remark}

Now consider the case $p\ne q$. For definiteness, we assume $p<q$. We 
use the more precise notation 
$\mathcal{M}(p,q,a,b)$ for the moduli space of $\U(p,q)$-Higgs bundles
such that $\deg(V)=a$, and $\deg(W)=b$, and write the Toledo 
invariant as 
  \begin{equation}
\tau= \tau(p,q,a,b) 
    = 2\frac{qa - pb}{p + q}.
\label{toledo-invariant-a-b}
\end{equation}

\begin{theorem}\label{prop:rigidity}
Suppose $(p,q,a,b)$ are such that  $p<q$ and 
$|\tau(p,q,a,b)|=p(2g-2)$. Then every element in  
$\mathcal{M}(p,q,a,b)$ is strictly semistable and 
decomposes as the direct sum of a polystable $\U(p,p)$-Higgs bundle 
with maximal Toledo invariant and a polystable vector bundle of rank 
$(q-p)$. If $\tau=p(2g-2)$, then 
  \begin{equation}
    \label{rigidity}
    \mathcal{M}(p,q,a,b)\cong \mathcal{M}(p,p,a,a-p(2g-2))
\times M(q-p, b-a +p(2g-2)),
  \end{equation}
\noindent where $M(q-p, b-a +p(2g-2))$ denotes the moduli space of polystable
bundles of degree $q-p$ and rank $b-a +p(2g-2)$.  In particular, the 
dimension at a smooth point in 
$\mathcal{M}(p,q,a,b)$ is $2+ (p^2+5q^2-2pq)(g-1)$, and it is
hence strictly smaller than the expected dimension. 

(A similar result holds if $\tau=-p(2g-2)$ and also if $p>q$.)
\end{theorem}
\begin{proof}
Let  $(E=V\oplus W,\Phi)\in\mathcal{M}(p,q,a,b)$ and suppose that 
$\tau(p,q,a,b)=p(2g-2)$. Then 
$\mu(V)-\mu(E)=g-1$ and $\mu(E)-\mu(W)=\frac{p}{q}(g-1)$. 
Since $\rk(\beta)$ and $\rk(\gamma)$ are at most $p$, it follows from
(\ref{eq:muE1-mu}) and (\ref{eq:mu-mu_1}) that 
$\rk(\beta)=\rk(\gamma)=p$. Let 
$W_\gamma=\im(\gamma)\otimes K^{-1}$ and let $W_\beta=\ker(\beta)$. 
Then $V\oplus W_\gamma$ is a $\Phi$-invariant subbundle of 
$V\oplus W$, and $\mu(V\oplus W_\gamma)=\mu(0\oplus W_{\beta})=\mu(E)$. 
We see that $(E,\Phi)$
is strictly semistable (as we already knew from Lemma 
\ref{lem:slope-bound}). Since it is polystable it must split as 
$$
(V\oplus W_\gamma,\Phi)\oplus (0\oplus W/W_{\gamma},0).
$$
It is clear that $(V\oplus W_\gamma,\Phi)\in 
\mathcal{M}(p,p,a,a-p(2g-2))$ and that 
$(V\oplus W_\gamma,\Phi)$ has maximal Toledo invariant, that is,
$\tau(p,p,a,a-p(2g-2))=2p(g-1)$. Also, using 
\begin{displaymath}
  0 \lto \ker(\Phi) \lto V\oplus W \lto (V\oplus W_{\gamma})\otimes K \lto 0.
\end{displaymath}
we see that  $W/W_{\gamma} \in M(q-p, b-a +p(2g-2))$. To complete the 
proof we observe that 
\begin{align*}
\dim \mathcal{M}^s(p,p,a,a-p(2g-2))+\dim  M^s(q-p, b-a +p(2g-2))\\
=1+(2p)^2(g-1)+1+(q-p)^2(g-1)=2+(p^2+5q^2-2pq)(g-1).
\end{align*}
Since $q>1$,  this is smaller than $1+(p+q)^2(g-1)$, the dimension of 
$\mathcal{M}(p,q,a,b)$ when the Toledo invariant is not maximal.
\end{proof}

\begin{remark}\label{remark:rigidity} The fact 
the moduli space has smaller dimension than expected may be viewed as 
a certain kind of rigidity. This phenomenon (for large Toledo 
invariant) has been studied from the point of view of representations 
of the fundamental group by D. Toledo 
\cite{toledo:1989} when $p=1$ and L. Hern\'andez 
\cite{hernandez:1991} when $p=2$. We deal here with the general case 
which, as far as we know, has not appeared previously in the 
literature. 
\end{remark}

\begin{corollary}
Fix $(p,q,a,b)$ such that $p<q$ and $\tau(p,q,a,b)=p(2g-2)$. Then 
$$
    \mathcal{M}(p,q,a,b)\cong \mathcal{M}_{K^2}(p,a-p(2g-2))
\times M(q-p, b-a +p(2g-2)).
$$
\end{corollary}
\begin{proof}
  It follows from Theorem \ref{prop:rigidity} and Proposition
  \ref{prop:p=q-toledo-max}.
\end{proof}

%%%%%%%%%%%%%%%%%%%%%%%%%%%%
\section{Morse theory}
\label{sec:morse-theory}
%%%%%%%%%%%%%%%%%%%%%%%%%%%%

Morse theoretic techniques for studying the topology of moduli spaces 
of Higgs bundles were introduced by Hitchin 
\cite{hitchin:1987,hitchin:1992}.  Though standard Morse theory 
cannot be applied to 
$\mathcal{M}(a,b)$ when it is not smooth, as we shall see 
in the following, we can still use Morse theory ideas to count 
connected components. Throughout this section we assume that $p$ and 
$q$ are any positive integers and that $(a,b)\in\Z\oplus\Z$ is such 
that $|\tau| \le\tau_M$, where $\tau$ is as in Definition 
\ref{defn:toledo} and 
$\tau_M$ is given by (\ref{eqn:tau-M}).

%%%%%%%%%%%%%%%%%%%%%%%%%%%%%%%
\subsection{The Morse function}\label{subs:Morse}
%%%%%%%%%%%%%%%%%%%%%%%%%%%%%%%

Consider the moduli space $\mathcal{M}(a,b)$ from the gauge theory 
point of view (cf.\ Section \ref{subs:upq}).  We can then define a 
real positive function 
\begin{equation}
\label{eq:def-morse-function}
\begin{aligned}
  f \colon \mathcal{M}(a,b) &\to \R \\
   [A,\Phi] &\mapsto \tfrac{1}{\pi}\norm{\Phi}^2\ ,
\end{aligned}  
\end{equation}
where the $L^2$-norm of $\Phi$ is $\norm{\Phi}^2=
\frac{\sqrt{-1}}{2}\int_X \tr(\Phi\Phi^*)$.

\goodbreak

We have the following result due to Hitchin \cite{hitchin:1987}.

\begin{proposition}
  \begin{itemize}
  \item[$(1)$] The function $f$ is proper.
  \item[$(2)$] The restriction of $f$ to $\mathcal{M}^s(a,b)$ is a
    moment map (up to a constant) for the Hamiltonian circle action
    $[A,\Phi] \mapsto [A,e^{i \theta}\Phi]$.
\item[$(3)$] If $\mathcal{M}(a,b)$ is smooth, then $f$ is a perfect
  Bott-Morse function.
  \end{itemize}\hfill\qed
\end{proposition}
Thus, if the moduli space is smooth, then its number of connected
components is bounded by the number of connected components of the 
subspace of local minima of $f$.  However, even if 
${\mathcal{M}}(a,b)$ is not smooth, $f$ can be used to obtain
information about the connected components of ${\mathcal{M}}(a,b)$
using the following elementary result.

\begin{proposition}
\label{prop:real-topology-exercise}
Let $Z$ be a Hausdorff space and let $f \colon Z \to \R$ be proper 
and bounded below.  Then $f$ attains a minimum on each connected 
component of $Z$ and, furthermore, if the subspace of local minima of 
$f$ is connected then so is $Z$. \hfill\qed
\end{proposition}

\noindent In particular this applies to our situation, giving:

\begin{proposition}
\label{prop:topology-exercise}
The function $f \colon \mathcal{M}(a,b) \to \R$ defined in
\eqref{eq:def-morse-function} has a minimum on each
connected component of ${\mathcal{M}}(a,b)$.  Moreover, if the
subspace of local minima of $f$ is connected then so is
$\mathcal{M}(a,b)$.  \hfill\qed
\end{proposition}
\begin{definition}\label{defn:N(a,b)}
Let 
\begin{equation}
 \mathcal{N}(a,b) =\{ (E,\Phi) \in \mathcal{M}(a,b)\suchthat
 \beta= 0 \;\;\mbox{or}\;\;\gamma=0\}. \label{minima}
\end{equation}
\end{definition}

\begin{proposition}\label{prop:absmin} For all $(E,\Phi)\in\mathcal{M}(a,b)$
\begin{equation}\label{eqn:absmin}
f(E,\Phi)\ge \frac{|\tau(a,b)|}{2}
\end{equation}
with equality if and only if $(E,\Phi)\in 
\mathcal{N}(a,b)$.
\end{proposition}

\begin{proof}
Writing out the first of Hitchin's equations 
\eqref{eq:hitchin1} for a $\U(p,q)$-Higgs bundle $(E,\Phi)$ in its
componenents on $V$ and $W$ we get the pair of equations 
\begin{align*}
    F(A_V) + \beta\beta^* + \gamma^*\gamma &=
    -\sqrt{-1}\mu\Id_V\omega\ ,  \\
    F(A_W) + \gamma\gamma^* + \beta^*\beta &=
    -\sqrt{-1}\mu\Id_W\omega\ , 
\end{align*}
where $A_V$ and $A_W$ are the components on $V$ and $W$, respectively,
of the unitary connection $A$ on $E=V \oplus W$.
Taking the trace and integrating over $X$ in the first of these
equations we get from  Chern-Weil theory 
\begin{displaymath}
  \deg(V) = \mu p -\tfrac{1}{\pi}\norm{\beta}^2 +
  \tfrac{1}{\pi}\norm{\gamma}^2\ ,
\end{displaymath}
where we have used $\int_X \omega = 2\pi$.  Since $\mu=\mu(E)$, this
is equivalent to
\begin{displaymath}
\tfrac{1}{\pi}\norm{\beta}^{2} -\tfrac{1}{\pi}\norm{\gamma}^{2} 
= p(\mu(E)- \mu(V)) 
=-\frac{\tau}{2}\ .
\end{displaymath}
But $f(E,\Phi)=\frac{1}{\pi}\norm{\beta}^{2}
  + \frac{1}{\pi}\norm{\gamma}^{2}$ and thus 
\begin{equation}\label{eqn:fabsmin}
\begin{aligned}
  f(E,\Phi)&=\tfrac{2}{\pi}\norm{\gamma}^{2}-\frac{\tau}{2}\\
&=\tfrac{2}{\pi}\norm{\beta}^{2}+\frac{\tau}{2}\ , 
\end{aligned}
\end{equation}
from which the result is immediate.
\end{proof}

\goodbreak

\noindent The above Proposition identifies $\mathcal{N}(a,b)$ as the set of 
global minima of $f$. The following Theorem, which is of fundamental 
importance to our approach, shows that there are no other local 
minima. 

\begin{theorem}
  \label{thm:minima}
Let $(E,\Phi)$ be  a polystable $\U(p,q)$-Higgs bundle 
in ${\mathcal{M}}(a,b)$. Then $(E,\Phi)$ is a local 
minimum of $f \colon  {\mathcal{M}}(a,b) \to \R$ if and only 
if $(E,\Phi)$ belongs to  $\mathcal{N}(a,b)$.
\end{theorem}
\begin{proof}
  This follows directly from Proposition~\ref{prop:absmin} above and
  Propositions~\ref{prop:thm-minima-2} and \ref{lem:thm-minima-3},
  which are given in Sections \ref{sec:stable-minima} and
  \ref{sec:reducible-minima}, respectively.
\end{proof}
\begin{remark}
  \label{rem:minima-for-lower-ranks}
  This Theorem was already known to hold when $p,q \leq 2$ (by the
  results of \cite{gothen:2001}, Hitchin \cite{hitchin:1987}, and Xia
  \cite{xia:2000}), and also when $p=q$ and $(p-1)(2g-2)<|\tau|\leq
  p(2g-2)$ by Markman-Xia \cite{markman-xia:2001}.
\end{remark}

\noindent Which section actually vanishes for a 
minimum is given by the following.
 
\begin{proposition} 
\label{prop:vanishing}
Let $(E,\Phi)\in \mathcal{N}(a,b)$.  Then
  \begin{enumerate}
  \item[$(1)$] $\gamma = 0$ if and only if $a/p \leq b/q$ 
  (i.e.\ $\tau \leq 0$). In this case, 
  $$f(\mathcal{N}(a,b)) = b - q(a+b)/(p+q)=-\frac{\tau}{2}\ .$$
  \item[$(2)$] $\beta = 0$ if and only if $a/p \geq b/q$ 
  (i.e.\ $\tau \geq 0$). In this case, 
  $$f(\mathcal{N}(a,b)) = a -p(a+b)/(p+q)=\frac{\tau}{2}\ .$$
\end{enumerate}
In particular, $\beta = \gamma = 0$ if and only if $a/p = b/q$ (i.e.
$\tau =0$) and, in this case, $f(E,\Phi) = 0$.
  
\end{proposition}
\begin{proof}
The relation between the conditions on $\tau$ and those on
$a/p-b/q$ follows directly from the definition of 
$\tau$ (cf.\ \eqref{toledo-invariant}). The rest follows immediately 
from (\ref{eqn:fabsmin}) and the fact that $f$ is, by definition, 
non-negative. Alternatively one can argue algebraically, using 
Lemma~\ref{lem:slope-bound} and polystability. 
\end{proof}

\begin{corollary}\label{minima-toledo-0}
If $a/p = b/q$ then $\mathcal{N}(a,b)\cong M(p,a)\times M(q,b)$.
\end{corollary}

\begin{proof} If  $a/p = b/q$, then any $(E,\Phi)\in \mathcal{N}(a,b)$
has $E=V\oplus W$ and $\Phi=0$. Polystability of $(E,\Phi)$ is thus 
equivalent to the polystability of $V$ and $W$. 
\end{proof}

%%%%%%%%%%%%%%%%%%%%%%%%%%%%%%%%%%%%%%%%%%%%%%%%%%%
\subsection{Critical points of the Morse function}
\label{sec:critical-points}
%%%%%%%%%%%%%%%%%%%%%%%%%%%%%%%%%%%%%%%%%%%%%%%%%%%
In this section we recall Hitchin's method
\cite{hitchin:1987,hitchin:1992} for determining the local minima of
$f$ and spell out how this works in the case of $\U(p,q)$-Higgs
bundles.  

Since $f$ is a moment map, a smooth point of the moduli space is a
critical point if and only if it is a fixed point of the circle
action.  To determine the fixed points, note that, if $(A,\Phi)$
represents a fixed point then there must be a 1-parameter family of
gauge transformations $g(\theta)$ taking $(A,\Phi)$ to
$(A,e^{i\theta}\Phi)$.  This gives an infinitesimal $\U(p)\times
\U(q)$-gauge transformation $\psi = \dot{g}$ which is covariantly
constant (i.e.\ $d_{A}\psi = 0$) and such that $[\psi,\Phi] = i\Phi$.
(Note that we can take $\psi$ to be trace-free.)  It follows that we
can decompose $E$ in holomorphic subbundles $F_{\lambda}$ on which
$\psi$ acts as $i\lambda$ and furthermore that $\Phi$ maps
$F_{\lambda}$ to $F_{\lambda + 1} \otimes K$.  We thus have the
following result.

\begin{proposition}
  \label{prop:variation-of-hodge}
  A $\U(p,q)$-Higgs bundle $(E,\Phi)$ in $\mathcal{M}(a,b)$
  represents a fixed point of the circle action if and only if it is a
  system of Hodge bundles, that is,
\begin{equation}
  E = F_1 \oplus \cdots \oplus F_m
  \label{eq:variation-of-hodge}
\end{equation}
for holomorphic vector bundles $F_i$ such that the restriction 
$$
\Phi_{i}:=\Phi_{|F_i} \in H^0(\Hom(F_i,F_{i+1})\otimes K)\ ,
$$
and the $F_i$ are direct sums of bundles contained in $V$ and $W$.
Furthermore, each $F_i$ is an eigenbundle for an infinitesimal
trace-free gauge transformation $\psi$.  If $\Phi_{i} \neq 0$, then
the weight of $\psi$ on $F_{i+1}$ is one plus the weight of $\psi$ on
$F_{i}$.  Moreover, if $(E,\Phi)$ is stable, then each restriction
$\Phi_{i}$ is non-zero and the $F_i$ are alternately contained in
$V$ and $W$.
\end{proposition}

\begin{proof}
  Only the last statement requires a proof.  But if some component of
  $\Phi$ vanished, or if some $F_{i}$ had a non-zero component in both
  $V$ and $W$, then $(E,\Phi)$ would be reducible and hence not
  stable.
\end{proof}

When $(E,\Phi)$ is stable
the decomposition $E = F_1 \oplus \cdots \oplus F_m$ gives a
corresponding decomposition of the bundle $U = \End(E)$ into
eigenbundles for the adjoint action of $\psi$:
\begin{displaymath} 
U = \bigoplus_{k=-m+1}^{m-1} U_k\ ,
\end{displaymath}
where
\begin{math}
  U_k = \bigoplus_{i-j = k} \Hom(F_j,F_i)
\end{math}
is the eigenbundle corresponding to the eigenvalue $ik$.

By Hitchin's calculations in \cite[\S 8]{hitchin:1992} (see also 
\cite[Section 2.3.2]{gothen:1995}) the eigenvalues of the 
Hessian of $f$ at a smooth critical point can be determined in the 
following way. 
\begin{proposition}
\label{prop:weight-k-space}
  Let $(E,\Phi)$ be a stable $\U(p,q)$-Higgs bundle which represents a
  critical point of $f$.  Then the eigenspace of the Hessian of $f$
  corresponding to the eigenvalue $-k$ is $\mathbb{H}^1$ of the
  following complex:
\begin{equation}
  \label{eq:weight-k-space}
  C^{\bullet}_{k} : 
  U_{k}^{+}
  \xrightarrow{\ad(\Phi)}
  U_{k+1}^{-} \otimes K,
\end{equation}
where we use the notation
\begin{align*}
  U_k^{+} &= U_{k} \cap U^+\ , \\
  U_k^{-} &= U_{k} \cap U^-\ ,
\end{align*}
with $U^+$ and $U^-$ as defined in Definition~\ref{def:u-notation}.
In particular $(E,\Phi)$ corresponds to a local minimum of $f$ if and
only if
\begin{displaymath}
\mathbb{H}^1( C^{\bullet}_{k}) = 0
\end{displaymath}
for all ${k \geq 1}$.
\hfill\qed
\end{proposition}

\begin{remark}
\label{rem:u-even-odd}
When $(E,\Phi)$ is a stable $\U(p,q)$-Higgs bundle, we know from
Proposition~\ref{prop:variation-of-hodge} that the $F_i$ are
alternately contained in $V$ and $W$. Thus we have
\begin{equation}
      \label{eq:u-even-odd}
U^+= \bigoplus_{\text{$k$ even}} U_k\ ; \quad
U^- = \bigoplus_{\text{$k$ odd}} U_k\ .
\end{equation}
In particular all the eigenvalues of the Hessian of $f$ are even.
\end{remark}

\begin{remark}
  The description in Proposition~\ref{prop:weight-k-space} of the
  eigenspace of the Hessian of $f$ gives rise to the long exact
  sequence
\begin{multline*}
  0 \lto \mathbb{H}^0(C^{\bullet}_k) \lto H^0(U^+_k) \lto
  H^0(U^-_{k+1}\otimes K) \lto  \mathbb{H}^1(C^{\bullet}_k) \\ 
    \lto H^1(U^+_k) \lto  H^1(U^-_{k+1}\otimes  K) \lto   
\mathbb{H}^2(C^{\bullet}_k)\lto 0\ .
\end{multline*} 
Suppose that $(E,\Phi)$ is a stable $\U(p,q)$-Higgs bundle.  The
vanishing result of Proposition~\ref{prop:stable-upq-vanishing} shows
that $\mathbb{H}^0(C^{\bullet}_{k}) = \mathbb{H}^2(C^{\bullet}_{k}) =
0$ for $k \neq 0$ (while $\mathbb{H}^0(C^{\bullet}_{0}) = \C$ and
$\mathbb{H}^2(C^{\bullet}_{0}) = 0)$. Hence one can use this exact
sequence, Remark~\ref{rem:u-even-odd}, and the Riemann--Roch formula 
to calculate the dimension of $\mathbb{H}^1(C^{\bullet}_{k})$ for any 
$k$ in terms of the ranks and the degrees of the $F_i$. This provides a
method for calculating the Morse index of $f$ at a critical point.  
However, we shall omit the formula since we have no need for it. 
\end{remark}

%%%%%%%%%%%%%%%%%%%%%%%%%%%%%%%%%%%%%%%%%%%%%%%%%%%
\subsection{Local minima and the adjoint bundle}
\label{sec:local-minima-adjoint}
%%%%%%%%%%%%%%%%%%%%%%%%%%%%%%%%%%%%%%%%%%%%%%%%%%%
In this section we give a criterion for $(E,\Phi)$ to be a local 
minimum in terms of the adjoint bundle.  This is the key step in the 
proof of Theorem~\ref{thm:minima}. We use the notation introduced in 
Section 
\ref{sec:critical-points}. 

Consider the complex $C^{\bullet}_{k}$ defined in 
\eqref{eq:weight-k-space} and let
\begin{displaymath}
  \chi(C^{\bullet}_{k}) = \dim \mathbb{H}^{0}(C^{\bullet}_{k})
  - \dim \mathbb{H}^{1}(C^{\bullet}_{k})
  + \dim \mathbb{H}^{2}(C^{\bullet}_{k}).
\end{displaymath}
\begin{proposition}
  \label{prop:adjoint-minima}
  Let $(E,\Phi)$ be a polystable $\U(p,q)$-Higgs bundle which is a fixed
  point of the $S^1$-action on $\mathcal{M}(a,b)$.
  Then $\chi(C^{\bullet}_k) \leq 0$ and equality holds if and only if 
  \begin{displaymath}
    \ad(\Phi) \colon U_{k}^{+} \to U_{k+1}^{-} \otimes K
  \end{displaymath}
  is an isomorphism.
\end{proposition}
\begin{proof}
  For simplicity we shall adopt the notation
  \begin{displaymath}
    \Phi_{k}^{\pm}
     = \ad(\Phi)_{| U_{k}^{\pm}} \colon U_{k}^{\pm} \lto U_{k+1}^{\mp}
     \otimes K.
  \end{displaymath}
  The key fact we need is that there is a natural $\ad$-invariant
  isomorphism $U \cong
  U^{*}$ under which we have $U^+ \cong
  (U^+)^*$, $U^- \cong (U^-)^*$ and $U_{k}^{\pm} \cong (U_{-k}^{\pm})^{*}$. 
  Since $\ad(\Phi)^{t} = \ad(\Phi) \otimes 1_{K^{-1}}$ under
  this isomorphism we have
  \begin{equation}
    \label{eq:Phi-symmetric}
    (\Phi_{k}^{\pm})^{t} = \Phi_{-k-1}^{\mp} \otimes 1_{K^{-1}}.
  \end{equation}
  We have the short exact sequence
  \begin{displaymath}
    0 \lto \ker(\Phi_{k}^{+}) \lto (U_{k+1}^{-} \otimes K)^* \lto
    \im(\Phi_{k}^{+}) \lto 0.
  \end{displaymath}
 {}From \eqref{eq:Phi-symmetric} we have $\ker(\Phi_{k}^{+,t}) \cong
  \ker(\Phi_{-k-1}^{-}) \otimes K^{-1}$. Thus, tensoring the above
  sequence by $K$, we obtain the short exact sequence
  \begin{displaymath}
    0 \lto \ker(\Phi_{-k-1}^{-}) \lto (U_{k+1}^{-})^* \lto
    \im(\Phi_{k}^{+}) \otimes K \lto 0.
  \end{displaymath}
  It follows that 
  \begin{displaymath}
    \deg(\im(\Phi_k^{+})) \leq \deg(U_{k+1}^{-}) +
    (2g-2)\rk(\Phi_{k}^{+}) + \deg(\ker(\Phi_{-k-1}^{-})).
  \end{displaymath}
  Combining this inequality with the fact that
  \begin{equation}
    \label {eq:degU_k^+}
      \deg(U_{k}^{+}) \leq
      \deg(\ker(\Phi_{k}^{+})) + \deg(\im(\Phi_{k}^{+})), 
  \end{equation}
  we obtain
  \begin{equation}
    \label{eq:deg-im-phi_k}
    \deg(U_k^{+}) \leq \deg(U_{k+1}^{-}) +
    (2g-2)\rk(\Phi_{k}^{+}) + \deg(\ker(\Phi_{-k-1}^{-}))
    + \deg(\ker(\Phi_{k}^{+})). 
  \end{equation}
  Since $(E,\Phi)$ is semistable, so is the Higgs bundle
  $(\End(E),\ad(\Phi))$.  Clearly the kernel $\ker(\Phi_{k}^{\pm}) \subseteq
  \End(E)$ is $\Phi$-invariant and hence, from semistability,
  \begin{displaymath}
    \deg(\ker(\Phi_{k}^{\pm})) \leq 0,
  \end{displaymath}
  for all $k$.  Substituting this inequality in
  \eqref{eq:deg-im-phi_k}, we obtain
  \begin{equation}
    \label{eq:deg-U_k^+}
    \deg(U_k^{+}) \leq \deg(U_{k+1}^{-}) +
    (2g-2)\rk(\Phi_{k}^{+}).
  \end{equation}
{}From the long exact sequence \eqref{eq:weight-k-space} and the
  Riemann--Roch formula we obtain
  \begin{align*}
    \chi(C^{\bullet}_k)
    &= \chi(U_{k}^{+}) - \chi(U_{k+1}^{-} \otimes K) \\
    &= (1-g)\bigl(\rk(U_{k}^{+}) + \rk(U_{k+1}^{-}) \bigr)
       +\deg(U_{k}^{+}) - \deg(U_{k+1}^{-}).
  \end{align*}
  Using this identity and the inequality \eqref{eq:deg-U_k^+} we see
  that 
  \begin{displaymath}
    \chi(C^{\bullet}_k) \leq
    (g-1)\bigl(2 \rk(\Phi_{k}^{+}) - \rk(U_{k}^{+}) - \rk(U_{k+1}^{-})
    \bigr). 
  \end{displaymath}
  Hence $\chi(C^{\bullet}_k) \leq 0$.  Furthermore, if equality holds
  we have 
  \begin{displaymath}
    \rk(\Phi_{k}^{+}) = \rk(U_{k}^{+}) = \rk(U_{k+1}^{-})
  \end{displaymath}
  and also equality must hold in \eqref{eq:deg-U_k^+} and so
  $\deg(\im(\Phi_{k}^{+})) = \deg(U_{k+1}^{-} \otimes K)$, showing
  that $\Phi_{k}^{+}$ is an isomorphism as claimed.
\end{proof}

\begin{corollary}
\label{cor:adjoint-minima}
  Let $(E,\Phi)$ be a stable $\U(p,q)$-Higgs bundle which represents a
  critical point of $f$.  This critical point is a local minimum 
  if and only if 
  \begin{displaymath}
    \ad(\Phi) \colon U_{k}^{+} \to U_{k+1}^{-} \otimes K
  \end{displaymath}
  is an isomorphism for all $k \geq 1$.
\end{corollary}

\begin{proof}
  By Proposition~\ref{prop:stable-upq-vanishing} we have
  $\mathbb{H}^0(C^{\bullet}_{k}) = \mathbb{H}^2(C^{\bullet}_{k}) = 0$
  for $k \geq 1$. Hence we have $-\chi(C^{\bullet}_{k}) =
  \mathbb{H}^1(C^{\bullet}_{k})$ and the result follows from
   Propositions \ref{prop:weight-k-space} and
  \ref{prop:adjoint-minima}.
\end{proof}

\begin{remark}
  \label{rem:HiggsG-minima}
  Let $(P,\Phi)$ be a $G$-Higgs bundle as defined in 
  Remark~\ref{rem:HiggsG} and define 
  \begin{align*}
    U &= P \times_{\Ad} \lie{g}^{\C}\ , \\
    U^+ &= P \times_{\Ad} \lie{h}^{\C}\ , \\
    U^- &= P \times_{\Ad} \lie{m}^{\C}\ . \\
  \end{align*}
  Then $U = U^+ \oplus U^-$ and if $(P,\Phi)$ is fixed under the 
  circle action we can write $U = \bigoplus U_{k}$ as a direct sum of 
  eigenbundles for an infinitesimal gauge transformation as before.  
  Thus we can define a complex
  $C^{\bullet}_{k}$ as in \eqref{eq:weight-k-space}.  If $(P,\Phi)$ is
  a stable $G$-Higgs bundle, then the Higgs vector bundle 
  $(U,\ad(\Phi))$ is semistable and so the proof of 
  Proposition~\ref{prop:adjoint-minima} goes through unchanged.  Thus 
  this key result is valid in the more general setting.
\end{remark}

%%%%%%%%%%%%%%%%%%%%%%%%%%%%%%%%%%%%%%%%%%%%%%%%%%%%%%%%%%%%%%%%%%%%%%
\subsection{Stable Higgs bundles}
\label{sec:stable-minima}
%%%%%%%%%%%%%%%%%%%%%%%%%%%%%%%%%%%%%%%%%%%%%%%%%%%%%%%%%%%%%%%%%%%%%%
In this section we prove Theorem~\ref{thm:minima} for stable Higgs 
bundles. The reducible (polystable) ones are dealt with in the next 
section.  We continue to use the notation of Section 
\ref{sec:critical-points}.
\begin{proposition}
  \label{prop:thm-minima-2}
  Let $(E,\Phi) = (F_{1} \oplus \cdots \oplus F_{m}, \Phi)$ be a stable
  $\U(p,q)$-Higgs bundle representing a critical point of $f$ such that 
  $m \geq 3$.  Then $(E,\Phi)$ is not a local minimum of $f$.
\end{proposition}
\begin{proof}
  Note that $U_{k} = 0$ for $\abs{k} \geq m$; in particular $U_{m} =
  0$.  We shall consider the cases when
  $m$ is odd and even separately.
  
  \emph{The case $m$ odd.}  In this case $m-1$ is even and so,
  using Remark~\ref{rem:u-even-odd} we see that $U_{m-1}^{+} = U_{m-1} \neq
  0$ while $U_{m}^{-} \subseteq U_m = 0$.  Hence $\ad(\Phi) \colon
  U_{m-1}^{+} \to U_{m}^{-} \otimes K$ cannot be an isomorphism and we
  are done by Corollary~\ref{cor:adjoint-minima}.

  \emph{The case $m$ even.} From Remark~\ref{rem:u-even-odd} we see that
  \begin{align*}
    U_{m-1}^{-} &= U_{m-1} = \Hom(F_{1},F_{m}) \\
    U_{m-2}^+ &= U_{m-2} = \Hom(F_{1},F_{m-1}) \oplus \Hom(F_{2},F_{m}).
  \end{align*}
  Thus,  by Corollary~\ref{cor:adjoint-minima} it
  suffices to prove that 
  \begin{displaymath}
    \ad(\Phi) \colon U_{m-2} \to U_{m-1} \otimes K
  \end{displaymath}
  is not an isomorphism.  In fact the restriction of $\ad(\Phi)$ to a
  fiber cannot even be injective.  Indeed, if it were, then its
  restriction to $\Hom(F_{1},F_{m-1})$ would be injective and hence
  $\Phi_{m-1}$ would also be injective.  Take a non-zero
  element $\eta \in \Hom(F_2, F_m)$ whose image is contained in the
  image of $\Phi_{m-1}$. Define $\zeta = \Phi^{-1}
  \eta \Phi \in \Hom(F_1,F_{m-1})$.  Then $\ad(\Phi)(\eta + \zeta) = 0$
  which is a contradiction.
\end{proof}

\begin{remark}
  \label{rem:morse-minima}
  Let $(E,\Phi)$ be a stable $\U(p,q)$-Higgs bundle with $\beta = 0$
  or $\gamma=0$.  Then, as pointed out above,
  Proposition~\ref{prop:absmin} shows that $(E,\Phi)$ is a local
  minimum of $f$.  This can also be seen from the Morse theory point
  of view, as follows.  Such a Higgs bundle either has $\beta
  =\gamma=0$ or it is a Hodge bundle of length $2$.  In the former
  case, clearly we have $\End(E) = U_{0}$.  In the
  latter case, $E = F_{1} \oplus F_{2}$ with $F_{1}= V$ and $F_{2} =
  W$ (if $\beta =0$) or vice-versa (if $\gamma=0$).  Hence $\End(E) =
  U_{-1} \oplus U_{0} \oplus U_{1}$.  Hence, in both cases $U_{k} = 0$ for
  $\abs{k} > 1$.  It follows that the complex $C^{\bullet}_{k}$ is
  zero for any $k>0$ and hence all eigenvalues of the Hessian of $f$
  are positive.
\end{remark}

%%%%%%%%%%%%%%%%%%%%%%%%%%%%%%%%%%%%%%%%%%%%%%%%%%%%%%%%%%%%%%%%%%%%%%%%%%%
\subsection{Reducible Higgs bundles}
\label{sec:reducible-minima}
%%%%%%%%%%%%%%%%%%%%%%%%%%%%%%%%%%%%%%%%%%%%%%%%%%%%%%%%%%%%%%%%%%%%%%%%%%%
In this section we shall finally conclude the proof of 
Theorem~\ref{thm:minima} by showing that it also holds for reducible 
Higgs bundles.  First we shall show that a reducible Higgs bundle 
which is not of the form given in Theorem~\ref{thm:minima} cannot be 
a local minimum of $f$; for this we use an argument similar to the one
given by Hitchin \cite[\S 8]{hitchin:1992} for the case of
$G=\mathrm{PSL}(n,\R)$. 

Let $(E,\Phi)$ be a strictly polystable $\U(p,q)$-Higgs bundle which
is a local minimum of $f$.  Since $f(E,\Phi)$ is the sum of the values
of $f$ on each of the stable direct summands (on the corresponding
lower rank moduli space), it follows that each stable direct summand
must be a local minimum in its moduli space and, therefore, a fixed
point of the circle action.  Hence $(E,\Phi)$ is itself fixed and thus
(cf.\ Proposition~\ref{prop:variation-of-hodge})
\begin{displaymath}
  E = \bigoplus F_{\lambda}\ ,
\end{displaymath}
where each $F_{\lambda}$ is an $i\lambda$-eigenbundle for an
infinitesimal trace-free $\U(p)\times \U(q)$-gauge transformation
$\psi$.  Moreover, if $\Phi_{|F_{\lambda}} \neq 0$, then its image is
contained in $F_{\lambda+1} \otimes K$. In analogy with the case of
stable $\U(p,q)$-Higgs bundles we write
\begin{displaymath}
  \End{E} = \bigoplus U_{\mu}\ ,
\end{displaymath}
where $U_{\mu}$ is the $i\mu$-eigenbundle for the adjoint action of
$\psi$. Let
\begin{align*}
  U_{\mu}^{+} &= U_{\mu} \cap U^+\ , \\
  U_{\mu}^{-} &= U_{\mu} \cap U^-\ ,
\end{align*}
then we can define a complex of sheaves
\begin{equation}
  \label{eq:pos-weight-space}
  C^{\bullet}_{>0} : 
  \bigoplus_{\mu>0}U_{\mu}^{+}
  \xrightarrow{\ad(\Phi)}
  \bigoplus_{\mu>1}U_{\mu}^{-} \otimes K\ .
\end{equation}
In this language Hitchin's criterion \cite[\S 8]{hitchin:1992} for
showing that a given fixed point is not a local minimum can be
expressed as follows.

\begin{lemma}
  \label{lem:minimum-criterion}
  Let $(E_t,\Phi_t)$ be a 1-parameter family of polystable
  $\U(p,q)$-Higgs bundles such that $(E_0,\Phi_0)$ is a fixed point of
  the circle action.  If the tangent vector $(\dot{E},\dot{\Phi})$ at
  $0$ is non-trivial and lies in the subspace
  \begin{displaymath}
    \mathbb{H}^1(C^{\bullet}_{>0})
  \end{displaymath}
  of the infinitesimal deformation space $\mathbb{H}^1(C^{\bullet})$
  of $(E_0,\Phi_0)$, then $(E_0,\Phi_0)$ is not a local minimum of
  $f$.
  \hfill\qed
\end{lemma}

\begin{proposition}
  \label{lem:thm-minima-3}
  Let $(E,\Phi)$ be a reducible $\U(p,q)$-Higgs bundle.  
  If $\beta \neq 0$ and $\gamma \neq 0$ then $(E,\Phi)$ is not a local minimum 
  of $f$.
\end{proposition}
\begin{proof}
  As we noted above, each stable direct summand of $(E,\Phi)$ is a
  local minimum on its moduli space and therefore (by Proposition
  \ref{prop:thm-minima-2}) it has $\beta=0$ or $\gamma=0$.  Hence we
  can choose two stable direct summands $(E' = V' \oplus W',\Phi')$
  and $(E''= V'' \oplus W'',\Phi'')$ such that $\gamma' \neq 0$ and
  $\beta'' \neq 0$ and $\beta'= \gamma''=0$.  It is clearly sufficient
  to show that $(E' \oplus E'',\Phi' \oplus \Phi'')$ is not a local
  minimum of $f$ on the corresponding moduli space and we can
  therefore assume that $(E,\Phi) = (E' \oplus E'',\Phi' \oplus
  \Phi'')$ without loss of generality.  We shall construct a family of
  deformations $(E_t,\Phi_t)$ of $(E,\Phi)$ satisfying the conditions
  of Lemma~\ref{lem:minimum-criterion}.
  
  By Lemma~\ref{lem:non-vanishing} both
  $H^1(\Hom(W'',W'))$ and $H^1(\Hom(V',V''))$ are non-vanishing, so
  let $\eta\in H^1(\Hom(V',V''))$ and $\sigma \in H^1(\Hom(W'',W'))$
  be non-zero.  We
  can then define a deformation of $(E,\Phi)$ by using that $\eta$
  defines an extension
  \begin{displaymath}
    0 \lto V'' \lto V^{\eta} \lto V' \lto 0\ ,
  \end{displaymath}
  while $\sigma$ defines an extension
  \begin{displaymath}
    0 \lto W' \lto W^{\sigma} \lto W'' \lto 0\ .
  \end{displaymath}
  Let $E^{(\eta,\sigma)} = V^{\eta} \oplus W^{\sigma}$ and define
  $\Phi^{(\eta,\sigma)}$ by the compositions
  \begin{align*}
    b^{(\eta,\sigma)} &\colon
      W^{\sigma} \lto W'' \overset{\beta''}{\lto} V'' \to V^{\eta}\ , \\
    c^{(\eta,\sigma)} &\colon
      V^{\eta} \lto V' \overset{\gamma'}{\lto} W' \lto W^{\sigma}.  \\
  \end{align*}
  Note that $(E^0,\Phi^0) = (E,\Phi)$ (the Higgs fields agree since
  $\beta' = \gamma'' = 0$).  It is then easy
  to see that $(E^{\eta,\sigma},\Phi^{\eta,\sigma})$ is stable: the
  essential point is that the destabilizing subbundles $V'$ and $W''$
  of $(E,\Phi)$ are not subbundles of the deformed Higgs bundle; we
  leave the details to the reader.
  
  Now define the family $(E_{t},\Phi_{t}) = (E^{(\eta t,\sigma t)},
  \Phi^{(\eta t,\sigma t)})$.  It is clear that the induced
  infinitesimal deformation of $E$ is
  $$
  \dot{E} = (\eta,\sigma)
  \in H^1(\Hom(V',V'')) \oplus H^1(\Hom(W'',W')) \subseteq
  H^1(\End(E))\ .
  $$
  Considering the holomorphic structure as given by a
  $\dbar$-operator on the underlying smooth bundle, our definition of
  $(E^{(\eta,\sigma)},\Phi^{(\eta,\sigma)})$ did not change the Higgs
  field but only the holomorphic structure on $E$.  Thus, taking a
  Dolbeault representative (cf.\ Remark~\ref{rem:dolbeault}) for
  $(\dot{E},\dot{\Phi}) \in \mathbb{H}^1(C^{\bullet})$ we see that the
  weights of $\psi$ on $(\dot{E},\dot{\Phi})$ are given by its weights
  on $\dot{E}$.  {}From Proposition~\ref{prop:variation-of-hodge} we
  have decompositions $E' = \bigoplus F_{k}'$ and $E'' = \bigoplus
  F_{k}''$ into eigenspaces of infinitesimal trace-free gauge
  transformations $\psi'$ and $\psi''$.  Note that the infinitesimal
  gauge transformation producing the decomposition of $E$ is $\psi =
  \psi' + \psi''$.  Clearly we have
  \begin{align*}
    F_{1}' &= V'\ ,  & F_{2}' &= W'\ , \\
    F_{1}'' &= W''\ , & F_{2}'' &= V'' \ .
  \end{align*}
  Let $\lambda_{V}'$ and $\lambda_{W}'$ be the weights of the action
  of $\psi'$ on $V'$ and $W'$ respectively, and analogously for $E''$.
  We then have that 
  \begin{align*}
    \lambda_{W}'  &= \lambda_{V}' + 1\ , &
    \lambda_{V}'' &= \lambda_{W}'' + 1\ .  
  \end{align*}
  and, since $\tr\psi' = \tr\psi'' = 0$,
  \begin{align*}
    \lambda_{V}' p' + \lambda_{W}' q' &= 0\ , \\
    \lambda_{V}'' p'' + \lambda_{W}'' q'' &= 0\ ,    
  \end{align*}
  where $p' = \rk(V')$, $q' = \rk(W')$, $p'' = \rk(V'')$ and $q'' =
  \rk(W'')$. 
{}From these equations we conclude that 
  \begin{align*}
    \lambda_{W}' - \lambda_{W}'' &=
      \frac{p'}{p' + q'} + \frac{p''}{p'' + q''} > 0\ , \\
    \lambda_{V}'' - \lambda_{V}' &=
      \frac{q''}{p'' + q''} + \frac{q'}{p' + q'} > 0\ .
  \end{align*}
  It follows that the weights of $\psi$ on $H^1(\Hom(W'',W'))$ and
  $H^1(\Hom(V',V''))$ are both positive and hence that
  $(\dot{E},\dot{\Phi})$ lies in a direct sum of positive weight
  spaces of $\psi$.  This concludes the proof of the Proposition.
\end{proof}
  \begin{lemma}
    \label{lem:non-vanishing}
    Let $(E' = V' \oplus W',\Phi')$ and $(E''= V'' \oplus W'',\Phi'')$
    be stable $\U(p,q)$-Higgs bundles of the same slope.  Then the
    cohomology groups $H^1(\Hom(V',V''))$ and $H^1(\Hom(W'',W'))$ are
    both non-vanishing.
  \end{lemma}
  \begin{proof}
    Since $\gamma'' = 0$, $V''$ is a 
    $\Phi$-invariant subbundle of $E''$. Thus $\mu(V'') <
    \mu(E'')$.  Using the Riemann--Roch formula and the equality
    $\mu(E'')  = \mu(E')$ we obtain
  \begin{align*}
    h^0(\Hom(V',V'') - h^1(\Hom(V',V'')
    &= p'p''(1-g + \mu(V'') - \mu(V')) \\
    &< p'p''(1-g + \mu(E') - \mu(V')). \\
  \end{align*}
  Since $\rk(\beta') \leq p'$ the inequality \eqref{eq:mu-mu_1} of 
  Corollary~\ref{cor:mu-muE_i} shows that $\mu(E') - \mu(V') \leq 
  g-1$ and we therefore deduce that 
  \begin{displaymath}
    h^0(\Hom(V',V'') - h^1(\Hom(V',V'') < 0,
  \end{displaymath}
  from which it follows that $H^1(\Hom(V',V'') \neq 0$.
    
  Similarly one sees that $H^1(\Hom(W'',W')) \neq 0$.
\end{proof}

%%%%%%%%%%%%%%%%%%%%%%%%%%%%%%%%%%%%%%%%%%%
\subsection{Local minima and connectedness}
\label{sec:minima-connectedness}
%%%%%%%%%%%%%%%%%%%%%%%%%%%%%%%%%%%%%%%%%%%

 In this section we obtain
connectedness results on $\mathcal{M}^{s}(a,b)$ and its closure
$\bar{\mathcal{M}}^{s}(a,b)$. We denote by 
$\mathcal{N}^{s}(a,b) \subseteq \mathcal{N}(a,b)$ the subspace
consisting of stable $\U(p,q)$-Higgs bundles, and denote its closure 
by $\bar{\mathcal{N}}^{s}(a,b)$.

The invariants $(a,b)$ will be fixed in the following and we shall
occasionally drop them from the notation and write $\mathcal{M} =
\mathcal{M}(a,b)$, etc.

\begin{proposition}
  \label{prop:closure-N^s}
  The closure of $\mathcal{N}^{s}$ in $\mathcal{M}$ coincides with
  $\bar{\mathcal{N}}^{s}$ and
  \begin{displaymath}
    \bar{\mathcal{N}}^{s} = \bar{\mathcal{M}}^{s} \cap \mathcal{N}\ .
  \end{displaymath}
\end{proposition}

\begin{proof}
  Clear.
\end{proof}

Now consider the restriction of the Morse function to
$\bar{\mathcal{M}}^{s}$,
\begin{displaymath}
  f \colon \bar{\mathcal{M}}^{s} \to \R\ .
\end{displaymath}

\begin{proposition}
  \label{prop:minima-closure}
  The restriction of $f$ to $\bar{\mathcal{M}}^{s}$ is proper and the
  subspace of local minima of this function coincides with
  $\bar{\mathcal{N}}^{s}$.
\end{proposition}

\begin{proof}
  Properness of the restriction follows from properness of $f$ and the
  fact that $\bar{\mathcal{M}}^{s}$ is closed in $\mathcal{M}$. By
  Proposition~\ref{prop:absmin} $f$ is constant on
  $\mathcal{N}$ and its value there is its global
  minimum on $\mathcal{M}$.  Thus $\bar{\mathcal{N}}^{s}$ is
  contained in the subspace of local minima of $f$.
  
  It remains to see that there are $f$ has no other local minima on
  $\bar{\mathcal{M}}^{s}$.  We already know that
  the subspace of local minima on $\mathcal{M}^{s}$ is
  $\mathcal{N}^{s}$.  Thus, since $\mathcal{M}^{s}$ is open in
  $\bar{\mathcal{M}^{s}}$, there cannot be any additional local
  minima on $\mathcal{M}^{s}$. We need to prove therefore 
  that there are no local minima in 
  $(\bar{\mathcal{M}^{s}} \setminus \mathcal{M}^{s}) \setminus
  \bar{\mathcal{N}}^{s}$.  So let $(E,\Phi)$ be a strictly poly-stable
  $\U(p,q)$-Higgs bundle representing a point in this space.  From 
  Proposition~\ref{prop:closure-N^s} we see that
  $\beta \neq 0$ and $\gamma \neq 0$.  In the proof of
  Proposition~\ref{lem:thm-minima-3} we constructed a family
  $(E_t,\Phi_t)$ of $\U(p,q)$-Higgs bundles such that $(E,\Phi) =
  (E_0,\Phi_0)$ and $(E_t,\Phi_t)$ is stable for $t \neq 0$.
  Furthermore we showed that the restriction of $f$ to this family
  does not have a local minimum at $(E_0,\Phi_0)$.  It follows that
  $(E,\Phi)$ is not a local minimum of $f$ on $\bar{\mathcal{M}^{s}}$.
\end{proof}

\begin{proposition}
  \label{prop:N-M-list}
  \begin{itemize}
  \item[$(1)$] If $\mathcal{N}(a,b)$ is connected, then so is
    $\mathcal{M}(a,b)$.
  \item[$(2)$] If $\mathcal{N}^s(a,b)$ is connected, then so is
    $\bar{\mathcal{M}}^s(a,b)$.
  \end{itemize}
\end{proposition}

\begin{proof}
  (1) In view of Proposition~\ref{prop:topology-exercise}, this follows
      from Theorem~\ref{thm:minima}.
      
  (2) If $\mathcal{N}^s(a,b)$ is connected, then so is its closure
      $\bar{\mathcal{N}}^s(a,b)$. But from
      Proposition~\ref{prop:minima-closure},
      $\bar{\mathcal{N}}^s(a,b)$ is the subspace of local minima
      of the proper positive map $f \colon
      \bar{\mathcal{M}}^s(a,b) \to \R$.  Hence the result follows
      from
      Proposition~\ref{prop:real-topology-exercise}.
\end{proof}

%%%%%%%%%%%%%%%%%%%%%%%%%%%%%%
\section{Local Minima as holomorphic triples}
\label{sec:stable-triples}
%%%%%%%%%%%%%%%%%%%%%%%%%%%%%%

The next step is to identify the spaces $\mathcal{N}(a,b)$ and  
$\mathcal{N}^s(a,b)$ as moduli spaces in their own right. By definition 
(cf.\ Definition \ref{defn:N(a,b)}), the Higgs bundles 
in $\mathcal{N}(a,b)$ all have $\beta=0$ or $\gamma=0$ in  their 
Higgs fields. Suppose first that $(E,\Phi)$ is a 
$\U(p,q)$-Higgs bundle with 
$\gamma = 0$.  Then $(E,\Phi)$ determines the triple 
$T = (E_{1},E_{2},\phi)$ where 
\begin{align*}
  E_{1} &= V\otimes K \\
  E_{2} &= W, \\
  \phi &= \beta\ .
\end{align*}
Conversely, given two holomorphic bundles $E_{1},E_{2}$ of rank $p$ 
and $q$ respectively, together with a bundle endomorphism $\Phi\in 
H^0(\Hom(E_2,E_1))$, we can use the above relations to define a
$\U(p,q)$-Higgs bundle with $\gamma = 0$.  Similarly, there is a 
bijective correspondence between $\U(p,q)$-Higgs bundles with 
$\beta = 0$ and holomorphic triples in which
\begin{align*}
  E_{1} &= W\otimes K, \\
  E_{2} &= V \\
  \phi &= \gamma.
\end{align*}

\noindent The triples $(E_1,E_2,\Phi)$ are examples of the holomorphic triples 
studied in \cite{bradlow-garcia-prada:1996} and
\cite{garcia-prada:1994}. 

%%%%%%%%%%%%%%%%%%
\subsection{Holomorphic triples}\label{subsec:triples}
%%%%%%%%%%%%%%%%%

We briefly recall the relevant definitions, referring to
\cite{bradlow-garcia-prada:1996} and \cite{garcia-prada:1994} for
details.  A \emph{holomorphic triple} on $X$, $T =
(E_{1},E_{2},\phi)$, consists of two holomorphic vector bundles
$E_{1}$ and $E_{2}$ on $X$ and a holomorphic map $\phi \colon E_{2}
\to E_{1}$.  Denoting the ranks $E_1$ and $E_2$ by $n_1$ and $n_2$,
and their degrees by $d_1$ and $d_2$, we refer to
$(\mathbf{n},\mathbf{d})=(n_1,n_2,d_1,d_2)$ as the {\em type} of the
triple.

A homomorphism from $T' = (E_1',E_2',\phi')$ to $T = (E_1,E_2,\phi)$   
is a commutative diagram 
\begin{displaymath}
  \begin{CD}
    E_2' @>\phi'>> E_1' \\
    @VVV @VVV  \\
    E_2 @>\phi>> E_1.
  \end{CD}
\end{displaymath}
$T'=(E_1',E_2',\phi')$ is a subtriple of $T = (E_1,E_2,\phi)$   
if the homomorphisms of sheaves $E_1'\to E_1$ and $E_2'\to E_2$ are 
injective. 

For any $\alpha \in \R$ the \emph{$\alpha$-degree} and 
\emph{$\alpha$-slope} of $T$ are 
defined to be 
\begin{align*}
  \deg_{\alpha}(T)
  &= \deg(E_{1}) + \deg(E_{2}) + \alpha
  \rk(E_{2}), \\ 
  \mu_{\alpha}(T)
  &=
  \frac{\deg_{\alpha}(T)}
  {\rk(E_{1})+\rk(E_{2})} \\ 
  &= \mu(E_{1} \oplus E_{2}) +
  \alpha\frac{\rk(E_{2})}{\rk(E_{1})+
    \rk(E_{2})}.
\end{align*}
The triple $T = (E_{1},E_{2},\phi)$ is 
\emph{$\alpha$-stable} if 
\begin{equation}\label{eqn:alpha-stability}
  \mu_{\alpha}(T')
  < \mu_{\alpha}(T)\ 
\end{equation}
for any proper sub-triple $T' = (E_{1}',E_{2}',\phi')$. Define  
\emph{$\alpha$-semistability} by replacing (\ref{eqn:alpha-stability})
with a weak inequality. A triple is called 
\emph{$\alpha$-polystable} if it is the direct sum of $\alpha$-stable
triples of the same $\alpha$-slope. It is {\it strictly 
$\alpha$-semistable (polystable)} if it is $\alpha$-semistable (polystable)
but not $\alpha$-stable. 

We denote the moduli space of isomorphism classes of $\alpha$-{\it 
polystable} triples of type $(n_1,n_2,d_1,d_2)$ by 
\begin{equation}\label{eqn:N(n,d)}
  \mathcal{N}_\alpha
  = \mathcal{N}_\alpha(\mathbf{n},\mathbf{d})
  = \mathcal{N}_\alpha(n_1,n_2,d_1,d_2)\ .
\end{equation}
Using Seshadri $S$-equivalence to define equivalence classes, this is 
the moduli space of equivalence classes of 
$\alpha$-semistable triples. The isomorphism 
classes of $\alpha$-{\it stable} triples form a subspace which we 
denoted by  
$\mathcal{N}_\alpha^s$.

\begin{proposition}[\cite{bradlow-garcia-prada:1996,garcia-prada:1994}]
  \label{prop:alpha-range} 
The moduli space   $\mathcal{N}_\alpha(n_1,n_2,d_1,d_2)$ is a complex 
analytic variety, which is projective when $\alpha$ is rational. A 
necessary condition for $\mathcal{N}_\alpha(n_1,n_2,d_1,d_2)$ to be 
non-empty is 
\begin{equation}
\begin{cases}
0\leq \alpha_m \leq \alpha \leq \alpha_M &\text{if $n_1\neq n_2$}\\ 
0\leq \alpha_m \leq \alpha &\text{if $n_1= n_2$} 
\end{cases}
\end{equation}
where 
\begin{align} 
 \alpha_m &= \mu_1-\mu_2, \label{alpha-bounds-m} \\
      \alpha_M &= 
      (1+ \frac{n_1+n_2}{|n_1 - n_2|})(\mu_1 - \mu_2) 
  \label{alpha-bounds-M} 
\end{align}
and $\mu_1=\frac{d_1}{n_1}$, $\mu_2=\frac{d_2}{n_2}$.
\end{proposition} 

Within the allowed range for $\alpha$ there is a discrete set of {\it 
critical values}. These are the values of 
$\alpha$ for which it is numerically possible to have a subtriple  
$T'=(E_1',E_2',\phi')$ such that $\mu(E'_1 \oplus E'_2)\ne\mu(E_1 
\oplus E_2)$ but 
$\mu_{\alpha}(T')=\mu_{\alpha}(T')$. All other values of $\alpha$ are 
called {\it generic}. The critical values of $\alpha$ are precisely 
the values for $\alpha$ at which the stability properties of a triple 
can change, i.e.\ there can be triples which are strictly 
$\alpha$-semistable, but either 
$\alpha'$-stable or $\alpha'$-unstable for $\alpha'\ne\alpha$. 

Strict $\alpha$-semistability can, in general, also occur at generic 
values for $\alpha$, but only if there can be subtriples with 
$\mu(E'_1 \oplus E'_2) = \mu(E_1 \oplus E_2)$ and
$\frac{n'_2}{n_1'+n_2'} = \frac{n_2}{n_1+n_2}$. 
In this case the triple is strictly
$\alpha$-semistable for all values of $\alpha$.  We refer to this
phenomenon as \emph{$\alpha$-independent semistability}. This cannot 
happen if $\GCD(n_2,n_1+n_2,d_1+d_2)=1$. 

%%%%%%%%%%%%%%%%%%%%%%%%%%
\subsection{Identification of $\mathcal{N}(a,b)$}
%%%%%%%%%%%%%%%%%%%%%%%%

The following result relates the stability conditions for holomorphic 
triples and that for $\U(p,q)$-Higgs bundles.

\begin{proposition}
  \label{prop:triple-higgs-stability}
  A $\U(p,q)$-Higgs bundle $(E,\Phi)$ with $\beta =0$ or $\gamma=0$ is
  \mbox{(semi)}stable if and only if the corresponding holomorphic
  triple  is
  $\alpha$-(semi)stable for $\alpha = 2g-2$.
\end{proposition}
\begin{proof}
   Let $T = (E_{1},E_{2},\phi)$
  be the triple corresponding to the Higgs bundle $(V\oplus W,\Phi)$.
  For definiteness we shall assume that $\gamma=0$ (of course, the same 
  argument applies if $\beta =0$). Thus $E_{1} = V 
  \otimes K$ and $E_{2} = W$ and, hence,
  \begin{displaymath}
    \deg(E_{1}) = \deg(V) + p(2g-2).
  \end{displaymath}
  Since $p = \rk(E_{1})$ and $q =
  \rk(E_{2})$ it follows that
  \begin{equation}
    \mu_{\alpha}(T)=
    \mu(E) + \frac{p}{p+q}(2g-2) + 
    \frac{q}{p+q}\alpha.
  \end{equation}
  If we set $\alpha = 2g-2$ we therefore have
  \begin{equation}
    \label{eq:mu-alpha}
    \mu_{\alpha}(T) = \mu(E) + 2g-2.
  \end{equation}
  Clearly the correspondence between holomorphic triples and
  $\U(p,q)$-Higgs bundles gives a correspondence between sub-triples
  $T' = (E_{1}',E_{2}',\phi')$ and $\Phi$-invariant subbundles of $E$
  which respect the decomposition $E = V \oplus W$ (i.e., subbundles
  $E' = V' \oplus W'$ with $V' \subseteq V$ and $W' \subseteq W$).
  Now, it follows from \eqref{eq:mu-alpha} that $\mu(E') < \mu(E)$ if
  and only if $\mu_{\alpha}(T') < \mu_{\alpha}(T)$ (and similarly for
  semistability), thus concluding the proof.
\end{proof}

We thus have the following important characterization of the subspace 
of local minima of $f$ on $\mathcal{M}(a,b)$.
\begin{theorem}
  \label{thm:minima=triple-moduli}
  Let $\mathcal{N}(a,b)$ be the subspace of local minima of
  $f$ on $\mathcal{M}(a,b)$ and let $\tau$ be the Toledo 
  invariant as defined in Definition \ref{defn:toledo}.
  
  If $a/p \leq b/q$, or equivalently if $\tau\le 0$,
   then $\mathcal{N}(a,b)$ can be 
  identified with the moduli space of $\alpha$-polystable triples of type
  $(p,q, a + p(2g-2), b)$, with $\alpha = 2g-2$. 
  
  If $a/p \geq b/q$,  or equivalently if $\tau\ge 0$,
  then $\mathcal{N}(a,b)$ can be 
  identified with the moduli space of $\alpha$-polystable triples of type
  $(q,p, b + q(2g-2),a)$, with $\alpha = 2g-2$.
  
  That is,
  \begin{displaymath}
    \mathcal{N}(a,b)
    \cong \begin{cases}
    \mathcal{N}_{2g-2}(p,q,a + p(2g-2),b)
    &\text{if $a/p \leq b/q$ (equivalently $\tau\le 0$)}\\ 
    \mathcal{N}_{2g-2}(q,p,b + q(2g-2),a)
    &\text{if $a/p \geq b/q$ (equivalently $\tau\ge 0$)}\end{cases}
  \end{displaymath}
\end{theorem}
\begin{proof}
This follows from Theorem~\ref{thm:minima}, 
Proposition~\ref{prop:vanishing}, and
Proposition~\ref{prop:triple-higgs-stability}.
\end{proof}

Thus, combining  Proposition~\ref{prop:N-M-list} and
Theorem~\ref{thm:minima=triple-moduli}, we get
 
\begin{theorem}
  \label{thm:connected-triples-connected-higgs}

  \begin{itemize}
    
  \item[$(1)$] Suppose $a/p \leq b/q$.
  If $\mathcal{N}_{2g-2}(p,q,a + p(2g-2),b)$
    is connected then $\mathcal{M}(a,b)$ is connected.
  If $\mathcal{N}^{s}_{2g-2}(p,q,a + p(2g-2),b)$ is connected
    then $\bar{\mathcal{M}}^s(a,b)$ is connected.
 
  \item[$(2)$] Suppose $a/p \geq b/q$.
If $\mathcal{N}_{2g-2}(q,p,b + q(2g-2),a)$ is connected then 
$\mathcal{M}(a,b)$ is connected.
If $\mathcal{N}^{s}_{2g-2}(q,p,b + q(2g-2),a)$ is connected
    then $\bar{\mathcal{M}}^s(a,b)$ is connected.
  \end{itemize}
\end{theorem}

%%%%%%%%%%%%%%%%%%%%%%%%%%%%%%%%%%%%%%%%%
\subsection{The Toledo invariant, $2g-2$, and $\alpha$-stability for triples}
%%%%%%%%%%%%%%%%%%%%%%%%%%%%%%%%%%%%%%%%%

In view of Theorems~\ref{thm:minima=triple-moduli} and
\ref{thm:connected-triples-connected-higgs}, it is important to
understand where $2g-2$ lies in relation to the range (given by
Proposition \ref{prop:alpha-range}) for the stability parameter
$\alpha$. Recall that for given $(p,q,a,b)$, the Toledo invariant
(Definition \ref{defn:toledo}) is constrained by $0\le|\tau|\leq
\tau_M$, where (see (\ref{eqn:tau-M})) $\tau_M=\min\{p,q\}(2g-2)$.
Recall also that $\alpha$ is constrained by the bounds given in
Proposition~\ref{prop:alpha-range}.  Whenever necessary we shal
indicate the dependence of $\alpha_m$ and $\alpha_M$ on $(p,q,a,b)$ 
by writing $\alpha_m=\alpha_m(p,q,a,b)$, and similarly for 
$\alpha_M$. 

\begin{lemma}\label{lemma:summary-alpha-mM}Fix $(p,q,a,b)$. Then
\begin{equation}\label{eqn:alpham}
\alpha_m(p,q,a,b)=(2g-2)-\frac{p+q}{2pq}|\tau|
\end{equation}

\noindent where $\tau$ is the Toledo invariant.  If $p\ne q$ then
\begin{equation}\label{eqn:alphaM}
\alpha_M(p,q,a,b)=\left(\ \frac{2\max\{p,q\}}{|p-q|}\right)
\alpha_m(p,q,a,b)\ .
\end{equation}
\noindent If  $p=q$ then $\alpha_M(p,q,a,b)=\infty$.
\end{lemma}

\begin{proof}By Theorem \ref{thm:minima=triple-moduli} the type of the
triple is determined by the sign of $\tau$. The result thus follows 
by applying (5.3) and (5.4) to triples of type 
$(p,q,a+p(2g-2),b)$ (if $\tau\le 0$) or type $(q,p,b+q(2g-2),a)$ (if 
$\tau\ge 0$).
\end{proof}

\begin{proposition}\label{prop:summary-range-2g-2} Fix $(p,q,a,b)$. 
Then  
\begin{equation}
0\le|\tau|\le\tau_M\Leftrightarrow
\begin{cases}
       0<\alpha_m(p,q,a,b)\le 2g-2 \le 
       \alpha_M(p,q,a,b) &\text{if $p\ne q$} \\
       0\le\alpha_m(p,q,a,b)\le 2g-2 &\text{if $p=q$} 
\end{cases}
\end{equation}
Furthermore, 
\begin{equation}
\tau=0\Leftrightarrow 2g-2=\alpha_m
\end{equation}
and 
\begin{equation}
|\tau|=\tau_M\Leftrightarrow \begin{cases}
       2g-2=\alpha_M &\text{if $p\ne q$} \\
       \alpha_m=0 &\text{if $p=q$} 
\end{cases}
\end{equation}

\end{proposition}
\begin{proof}
Using (\ref{eqn:alpham}) and (\ref{eqn:alphaM}) we see that 
$0\le|\tau|\le\tau_M$ is equivalent to
\begin{equation}\label{eqn:alphm-bounds}
2g-2\ge\alpha_m\ge \left(\frac{|p-q|}{2\max\{p,q\}}\right)(2g-2)\ , 
\end{equation}
and hence also (assuming $p\ne q$) to 
\begin{equation}\label{eqn:alphM-bounds}
\left(\frac{2\max\{p,q\}}{|p-q|}\right)(2g-2)\ge\alpha_M\ge (2g-2)\ .
\end{equation}
In both (\ref{eqn:alphm-bounds}) and (\ref{eqn:alphM-bounds}), we get 
equality in the first place if and only if 
$\tau=0$, and in the second place if and only if
$|\tau|=\tau_M$. Notice that $\frac{|p-q|}{2\max\{p,q\}}$ is strictly 
positive if $p\ne q$ and is zero if $p=q$. The results follow. 
\end{proof} 

These results are summarized in Figure 1, which can be used as 
follows. For any allowed value of $\tau$, draw a horizontal line at 
height 
$\tau$. The corresponding range for 
$\alpha$ and the relative location of 
$2g-2$ are then read off from the $\alpha$-axis.

\begin{remark}\label{remark:MWbound}
The above proposition gives another explanation for the Milnor--Wood 
inequality in Corollary \ref{cor:toledo}. Using the fact that the 
non-emptiness of 
$\mathcal{M}(a,b)$ is equivalent to the non-emptiness of $\mathcal{N}(a,b)$
and hence to that of either
$\mathcal{N}_{2g-2}(p,q,a+p(2g-2),b)$ or 
$\mathcal{N}_{2g-2}(q,p,b+q(2g-2),a)$, we see that the Milnor--Wood
inequality is equivalent to the condition that 
$2g-2$ lies within the range where  $\alpha$-polystable triples of the 
given kind exist.

\end{remark}

%%%%%%%%%%%%%%%%%%%%%%%%%%%%%%%%%%%%%%%%%
\subsection{Moduli spaces of triples}
\label{sect:N(n,d)}
%%%%%%%%%%%%%%%%%%%%%%%%%%%%%%%%%%%%%%%%%

Proposition \ref{prop:summary-range-2g-2} shows that in order to study
$\mathcal{N}(a,b)$ for different values of the Toledo invariant, we
need to understand the moduli spaces of triples for values of $\alpha$
that may lie anywhere (including at the extremes $\alpha_m$ and
$\alpha_M$) in the $\alpha$-range given in
Proposition~\ref{prop:alpha-range}. The information we need can be
found in \cite{bradlow-garcia-gothen:2002:triples}.
From the results in \cite{bradlow-garcia-gothen:2002:triples} we get the
following for triples of type $(n_1,n_2,d_1,d_2)$.

\begin{theorem}\label{thm:summary-Ns}{\rm [Theorem A in 
\cite{bradlow-garcia-gothen:2002:triples}]}

\begin{enumerate}
\item[$(1)$] A triple $T=(E_1,E_2,\phi)$ of type $(n_1,n_2,d_1,d_2)$
is $\alpha_m$-polystable if and only if 
$\phi=0$ and $E_1$ and $E_2$ are polystable. We thus have 
$$
\mathcal{N}_{\alpha_m}(n_1,n_2,d_1,d_2)\cong M(n_1,d_1) \times M(n_2,d_2).
$$
where $M(n,d)$ denotes the moduli space of polystable bundles of rank 
$n$ and degree $d$.  In particular, 
$\mathcal{N}_{\alpha_m}(n_1,n_2,d_1,d_2)$ is non-empty and irreducible. 

\item[$(2)$] If $\alpha>\alpha_m $ is any value such that 
$2g-2\leq\alpha$  (and $\alpha<\alpha_M$ if $n_1\ne n_2)$) then 
$\mathcal{N}^s_\alpha(n_1,n_2,d_1,d_2)$ is non-empty and
irreducible. Moreover: 
%\begin{itemize}
%\item 

\quad $\bullet$ If $n_1=n_2 =n$ then $\mathcal{N}^s_{\alpha}(n,n,d_1,d_2)$ is 
birationally  equivalent to a 
$\mathbb{P}^N$-fibration  over
$M^s(n,d_2)\times \mathrm{Sym}^{d_1-d_2}(X)$, where 
$M^s(n,d_2)$ denotes the subspace of stable bundles of type $(n,d_2)$,
$\mathrm{Sym}^{d_1-d_2}(X)$ is the symmetric product, and the fiber
 dimension is $N=n(d_1-d_2)-1$.

%\item 
\quad $\bullet$ If $n_1>n_2$ then 
$\mathcal{N}^s_\alpha(n_1,n_2,d_1,d_2)$ is birationally 
equivalent to a $\mathbb{P}^N$-fibration over 
$M^s(n_1-n_2,d_1-d_2) \times M^s(n_2,d_2)$, where the fiber dimension is
$N=n_2d_1-n_1d_2+n_2(n_1-n_2)(g-1)-1$.  The birational 
equivalence is an isomorphism if $\GCD(n_1-n_2,d_1-d_2)=1$ and 
$\GCD(n_2,d_2)=1$.

%\item 
\quad $\bullet$ If $n_1<n_2$ then 
$\mathcal{N}^s_\alpha(n_1,n_2,d_1,d_2)$ is birationally 
equivalent to a $\mathbb{P}^N$-fibration over 
$M^s(n_2-n_1,d_2-d_1) \times M^s(n_1,d_1)$, where 
the fiber dimension is $N=n_2d_1-n_1d_2+n_1(n_2-n_1)(g-1)-1$. The 
birational equivalence is an isomorphism if 
$\GCD(n_2-n_1,d_2-d_1)=1$ and $\GCD(n_1,d_1)=1$.
%\end{itemize}
\newline
In particular, if $n_1\ne n_2$ then 
$\mathcal{N}^s_\alpha(n_1,n_2,d_1,d_2)$ is a smooth manifold of 
dimension $(g-1)(n_1^2 + n_2^2 - n_1 n_2) - n_1 d_2 + n_2 d_1 + 1$. 

\item[$(3)$] If $n_1\ne n_2$ then $\mathcal{N}_{\alpha_M}(n_1,n_2,d_1,d_2)$ is 
non-empty and irreducible. Moreover 
\begin{equation}
\mathcal{N}_{\alpha_M}(n_1,n_2,d_1,d_2)\cong
\begin{cases}
M(n_2,d_2) \times M(n_1-n_2,d_1-d_2) &\text{if $n_1>n_2$}\\ 
M(n_1,d_1) \times M(n_2-n_1,d_2-d_1) &\text{if $n_1<n_2$.} 
\end{cases}
\end{equation}
\end{enumerate}
\end{theorem}

\begin{theorem}\label{thm:summary-n1=n2}{\rm [Corollary 8.2 and Theorem 8.10 in
\cite{bradlow-garcia-gothen:2002:triples}]}

If $n_1=n_2=n$ then:
\begin{itemize}

\item[$(1)$]If $\alpha_m=0$, i.e.\ if $d_1=d_2\ (=d)$, then 
$\mathcal{N}_\alpha(n,n,d,d)\cong M(n,d)$ for all $\alpha> 0$.
In particular 
$\mathcal{N}_\alpha(n,n,d,d)$ is non-empty and 
irreducible. 

\item[$(2)$] If $0<d_1-d_2<\alpha$, then $\mathcal{N}_{\alpha}(n,n,d_1,d_2)$ 
is non-empty and irreducible. 

\end{itemize}
\end{theorem}

\begin{remark}Notice that if $n_1=n_2$ and $\alpha_m=0$, then 
$\mathcal{N}_{\alpha}(n,n,d,d)\cong M(n,d)$ for all $\alpha>0$, while
$\mathcal{N}_{0}(n,n,d,d)\cong M(n,d)\times M(n,d)$. The picture is quite
different if we restrict to the stable points in the moduli spaces. 
In fact there are no stable points in $\mathcal{N}_{0}(n,n,d,d)$, 
i.e., $\mathcal{N}^s_{0}(n,n,d,d)$ is empty, while 
$\mathcal{N}^s_{\alpha}(n,n,d,d)\cong M^s(n,d)$ for 
$\alpha>0$.
\end{remark}

\begin{proposition}{\rm [Proposition 2.6 and Lemma 2.7 in 
\cite{bradlow-garcia-gothen:2002:triples}]}
\label{triples-critical-range} $\mathrm{}$

\begin{itemize}
\item[$(1)$]
If $\alpha\in [\alpha_m,\alpha_M]$ is generic and
$\GCD(n_1,n_1+n_2,d_1+d_2)=1$, then  
$$\mathcal{N}_\alpha(n_1,n_2,d_1,d_2)= 
\mathcal{N}^s_\alpha(n_1,n_2,d_1,d_2) .$$
\noindent In particular, the moduli space 
$\mathcal{N}_\alpha(n_1,n_2,d_1,d_2)$ is non-empty and irreducible if 
in addition $2g-2\leq\alpha$. 
 
\item[$(2)$]  Let $m\in \Z$ be such that $\GCD(n_1+n_2,d_1+d_2-mn_1)=1$. Then 
 $\alpha=m$ is not a critical value and there are no 
 $\alpha$-independent semistable triples.
\end{itemize}
\end{proposition}

%%%%%%%%%%%%%%%%%%%%%%%%%%%%%%%%%%%%%%%%%
\section{Main results}\label{sec:main-results}
%%%%%%%%%%%%%%%%%%%%%%%%%

We now use the results of Section \ref{sect:N(n,d)}, applied to the 
case $\alpha=2g-2$, to deduce our main results on the moduli spaces 
of $\U(p,q)$-Higgs bundles, and hence for the representation spaces 
$\mathcal{R}(\PU(p,q))$ and $\mathcal{R}_{\Gamma}(\U(p,q))$ 
(defined in section \ref{sec:background}). Recall that we identified 
components of 
$\mathcal{R}(\PU(p,q))$ labeled by $[a,b]\in \Z\oplus\Z/(p+q)\Z$,
and similarly identified components of 
$\mathcal{R}_{\Gamma}(\U(p,q))$ labeled by $(a,b)\in \Z\oplus\Z$. 
Our arguments proceed along the following lines: 
\begin{itemize}

\item By Proposition 
\ref{prop:principal-jac} 
$\mathcal{R}_{\Gamma}(a,b)$ is a $\U(1)^{2g}$-fibration over 
$\mathcal{R}[a,b]$. The number of connected components of 
$\mathcal{R}_{\Gamma}(a,b)$ is thus greater than or equal to that of 
$\mathcal{R}[a,b]$. In particular, $\mathcal{R}[a,b]$ is connected 
whenever $\mathcal{R}_{\Gamma}(a,b)$ is. 

\item By Proposition 
\ref{prop:R=M} there is a homeomorphism between 
$\mathcal{R}_\Gamma(a,b)$ and the moduli space $\mathcal{M}(a,b)$ of 
$\U(p,q)$-Higgs bundles. This restricts to give a homeomorphism 
between $\mathcal{R}^*_\Gamma(a,b)$ and 
$\mathcal{M}^s(a,b)$. 
\item By Proposition \ref{prop:topology-exercise} the number of 
connected components of $\mathcal{M}(a,b)$ is bounded above by the 
number of connected components in the subspace of local minima for 
the Bott-Morse function defined in Section \ref{subs:Morse}. By   
Proposition \ref{prop:N-M-list} the same conclusion holds for 
$\mathcal{M}^s(a,b)$.
\item By Theorems~\ref{thm:minima} and
\ref{thm:minima=triple-moduli} we can identify the subspace of local 
minima as a moduli space of 
$\alpha$-stable triples, with $\alpha=2g-2$. 
  
\end{itemize}
Summarizing, we have:
\begin{align*}
|\pi_0(\mathcal{R}[a,b])| \leq |\pi_0(\mathcal{R}_\Gamma(a,b))|&
=|\pi_0(\mathcal{M}(a,b))|\\& \leq |\pi_0(\mathcal{N}(a,b))|
=|\pi_0(\mathcal{N}_{2g-2}(n_1,n_2,d_1,d_2))|
\end{align*}
where $|\pi_0(\cdot)|$ denotes the number of components, and (in the 
notation of Section 
\ref{sec:stable-triples}) the moduli space of triples which appears
in the last line is either 
$\mathcal{N}_{2g-2}(p,q,a + p(2g-2),b)$ (if $a/p \leq b/q$) or
$\mathcal{N}_{2g-2}(q,p,b + q(2g-2),a)$.
Similarly, we get that 

\begin{align*}
|\pi_0(\bar{\mathcal{R}}^*[a,b])| \leq 
|\pi_0(\bar{\mathcal{R}}^*_\Gamma(a,b))|& 
=|\pi_0(\bar{\mathcal{M}}^s(a,b))|\\& 
\leq |\pi_0(\bar{\mathcal{N}}^s(a,b))|
=|\pi_0(\bar{\mathcal{N}}^s_{2g-2}(n_1,n_2,d_1,d_2))|
\end{align*}
In particular, if the moduli spaces of triples are connected, then so 
are the Higgs moduli spaces and the moduli spaces of representations.

%%%%%%%%%%%%%%%%%%%%%%%%%%%%%%%%%%%%%%%%%
\subsection{Moduli spaces of Higgs bundles}
\label{sec:results-U(p,q)-Higgs}
%%%%%%%%%%%%%%%%%%%%%%%%%%%%%%%%%%%%%%%%%

We begin with results for the $\U(p,q)$-Higgs moduli spaces.  Recall 
from Proposition 
\ref{prop:smoothness-higgs} that, whenever the moduli space 
$\mathcal{M}^s(a,b)$ of stable 
$\U(p,q)$-Higgs bundles with invariants $(a,b)$ is non-empty, it is a 
smooth complex manifold of dimension $1 + (p+q)^2(g-1)$.  We shall 
refer to this dimension as the \emph{expected dimension} in the 
following.

\begin{theorem}\label{thm:summary-M(a,b)}
 Let $(p,q)$ be any pair of positive integers and let 
  $(a,b)\in\Z\oplus\Z$ be such that $0\le|\tau(a,b)|\le\tau_M$. 

\begin{itemize}

\item[$(1)$] If either of the following sets of conditions apply,
then the moduli space $\mathcal{M}^s(a,b)$ is a non-empty smooth 
manifold of the expected dimension, with connected closure 
$\bar{\mathcal{M}}^s(a,b)$: 
\begin{enumerate}
\item[$(i)$] $0<|\tau(a,b)|<\tau_M$ ,
\item[$(ii)$] $|\tau(a,b)|=\tau_M$ and $p=q$.
\end{enumerate}

\item[$(2)$] If any one of the 
following sets of conditions apply, then the moduli space 
$\mathcal{M}(a,b)$ is non-empty and connected: 
\begin{enumerate}
\item[$(i)$] $\tau(a,b)=0$,
  
\item[$(ii)$] $|\tau(a,b)|=\tau_M$ and $p\neq q$,
  
\item[$(iii)$] $(p-1)(2g-2) < |\tau| \leq \tau_M = p(2g-2)$ and 
$p=q$.

\end{enumerate}
\end{itemize}
\end{theorem}

\begin{proof}   

$(2)$ By Proposition \ref{prop:summary-range-2g-2}
condition $(i)$ implies that $\alpha_m<2g-2<\alpha_M$ for the triples 
corresponding to points in $\mathcal{N}(a,b)$. Thus Theorem 
\ref{thm:summary-Ns}(2) (together with Theorem 
\ref{thm:minima=triple-moduli}) implies that 
$\mathcal{N}(a,b)$ is non-empty and connected.  Similarly, condition 
$(ii)$ implies that 
$\alpha_m=0$, and we can apply Theorem \ref{thm:summary-n1=n2}(1).
The rest follows from Theorem 
\ref{thm:connected-triples-connected-higgs}. 
$(3)$ By Proposition \ref{prop:summary-range-2g-2}, 
the conditions in $(i)$ and $(ii)$ are equivalent to 
$\alpha_m=2g-2$ and $\alpha_M=2g-2$ respectively. It follows by
parts (1) and (3) of Theorem \ref{thm:summary-Ns} (together with 
Theorem \ref{thm:minima=triple-moduli}) that $\mathcal{N}(a,b)$ is 
non-empty and connected. The rest follows from Theorem
\ref{thm:connected-triples-connected-higgs}. 
 
For $(iii)$, we use the fact that $|\tau|=|b-a|$ if $p=q$. The 
condition on $|\tau|$ is thus equivalent to $d_1-d_2<2g-2$ for the 
triples  corresponding to points in $\mathcal{N}(a,b)$. The result 
thus follows by  Theorem \ref{thm:summary-n1=n2}(2).

\end{proof}

\begin{remark}Combining (1) and $(i)$--$(ii)$ of (2) in 
Theorem \ref{thm:summary-M(a,b)}, we see that {\it the moduli space 
$\mathcal{M}(a,b)$ is non-empty for all $(p,q,a,b)$ such that
$0\le|\tau|\le\tau_M$}.
\end{remark}
\begin{remark} In Theorem \ref{prop:rigidity} we gave a detailed
description for $\mathcal{M}(a,b)$ in the case that $p\ne q$ and 
$|\tau(a,b)|=\tau_M$. The description was complete, provided that the 
space was non-empty. By the previous remark we can now remove this 
caveat. 
\end{remark}

In general, the stable locus $\mathcal{M}^s(a,b)$ is not the full 
moduli space and the full moduli space $\mathcal{M}(a,b)$ is not 
smooth. Singularities can occur at points representing strictly 
semistable objects, and these can also account for singularities in 
$\mathcal{N}(a,b)$, the space of local minima (as in Section 
\ref {sec:stable-triples}). These types of singularities are 
prevented by certain coprimality condition: 

\begin{proposition}\label{prop:gcd}
Suppose that $\GCD(p+q,a+b)=1$. Then: 
\begin{enumerate}
\item[$(1)$] $\mathcal{M}(a,b)$ is smooth.
\item[$(2)$] $\alpha=2g-2$ is not a critical value for triples of
type $(p,q,a + p(2g-2),b)$ or $(q,p,b + q(2g-2),a)$.
\item[$(3)$] The moduli spaces
$\mathcal{N}_{2g-2}(p,q,a + p(2g-2),b)$ and 
$\mathcal{N}_{2g-2}(q,p,b + q(2g-2),a)$ are non-empty, 
smooth and irreducible.
\end{enumerate}
\end{proposition}
\begin{proof}
\begin{enumerate}
\item[$(1)$] This is simply a re-statement of (2) in Proposition
\ref{prop:R=M}.
\item[$(2)$] Apply
Proposition \ref{triples-critical-range} (2) with  
$(n_1,n_2,d_1,d_2)$ equal to $(p,q,a + p(2g-2),b)$ or $(q,p,b + 
q(2g-2),a)$ and $m=2g-2$. 
\item[$(3)$] Since $\GCD(p+q,a+b)=1$ implies 
$\GCD(p,p+q,b+a+q(2g-2))=1$ (or $\GCD(q,p+q,b+a+p(2g-2))=1$), 
the result follows from (2) and  Proposition 
\ref{triples-critical-range} (1). 
\end{enumerate}
\end{proof}

\begin{theorem}\label{thm:coprime-M(a,b)}
Let $(p,q)$ be any pair of positive integers and let $(a,b)$ be such 
that $0\le|\tau(a,b)|\le\tau_M$. Suppose also that $\GCD(p+q,a+b)=1$. 
Then the moduli space $\mathcal{M}(a,b)$ is a (non-empty) smooth, 
connected manifold of the expected dimension. 
\end{theorem}

\begin{proof} Combine Proposition \ref{prop:gcd} and Theorem
\ref{thm:connected-triples-connected-higgs}.
\end{proof}

Theorems \ref{thm:summary-M(a,b)} plus 
\ref{thm:coprime-M(a,b)} are equivalent to Theorem A in the Introduction.

%%%%%%%%%%%%%%%%%%%%%%%%%%%%%%%%%%%%%%%%%%%%%%%%%%%%%%%%%%%%%%%%%%%%%%
\subsection{Moduli spaces of representations}
%%%%%%%%%%%%%%%%%%%%%%%%%%%%%%%%%%%%%%%%%%%%%%%%%%%%%%%%%%%%%%%%%%%%%%

Using  Proposition~\ref{prop:R=M} we can translate the results of 
Section~\ref{sec:results-U(p,q)-Higgs} into results about the 
representation spaces $\mathcal{R}_\Gamma(a,b)$ and 
$\mathcal{R}^*_\Gamma(a,b)$ (for $\U(p,q)$ representations of the 
surface group $\Gamma$). We denote  the closure of 
$\mathcal{R}_\Gamma^*(a,b)$ in 
$\mathcal{R}_\Gamma(a,b)$ by 
$\bar{\mathcal{R}}_\Gamma^*(a,b)$.

\begin{theorem}\label{thm:results-R-Gamma}

  Let $(p,q)$ be any pair of positive
  integers and let $(a,b)\in\Z\oplus\Z$ be such that $0\le|\tau(a,b)|\le\tau_M$. 

\begin{itemize}

\item[$(1)$] The moduli space $\mathcal{R}_{\Gamma}(a,b)$ is 
non-empty.

\item[$(2)$] If either of the following sets of conditions apply,
then the moduli space $\mathcal{R}_\Gamma^*(a,b)$ is a non-empty 
smooth manifold of the expected dimension, with connected closure 
$\bar{\mathcal{R}}_\Gamma^*(a,b)$ in $\mathcal{R}_{\Gamma}(a,b)$: 
\begin{enumerate}
\item[$(i)$] $0<|\tau(a,b)|<\tau_M$ ,
\item[$(ii)$] $|\tau(a,b)|=\tau_M$ and $p=q$.
\end{enumerate}

\item[$(3)$]  If any one of the 
following sets of conditions apply, then the moduli space 
  $\mathcal{R}_\Gamma(a,b)$ is connected:
\begin{enumerate}
\item[$(i)$] $\tau(a,b)=0$,
  
\item[$(ii)$] $|\tau(a,b)|=\tau_M$ and $p\neq q$,
  
\item[$(iii)$] $(p-1)(2g-2) < |\tau| \leq \tau_M = p(2g-2)$ and $p=q$,
  
\item[$(iv)$] $\GCD(p+q,a+b)=1$
\end{enumerate}

\item[$(4)$] If  $\GCD(p+q,a+b)=1$ then $\mathcal{R}_\Gamma(a,b)$ 
is a smooth manifold of the expected dimension. 
\end{itemize}
\end{theorem}

\begin{proof} By Proposition~\ref{prop:R=M}, 
this follows from Theorem 
\ref{thm:summary-M(a,b)} and \ref{thm:coprime-M(a,b)}.
\end{proof}

\begin{theorem}\label{thm:RgammaTmax}
Let $(p,q)$ be any pair of positive integers such that $p\ne q$, and 
let $(a,b)$ be such that 
$|\tau(a,b)|=\tau_M$.  Then every representation in 
$\mathcal{R}_\Gamma(a,b)$ is reducible  (i.e.\ 
$\mathcal{R}^*_\Gamma(a,b)$ is empty).  If $p < q$, then every such 
representation decomposes as a direct sum of a semisimple 
representation of 
$\Gamma$ in $\U(p,p)$ with maximal Toledo invariant and a semisimple 
representation in  $\U(q-p)$.  Thus, if $\tau = p(2g-2)$ then there 
is an isomorphism 
  \begin{displaymath}
    \mathcal{R}_\Gamma(p,q,a,b) \cong \mathcal{R}_\Gamma(p,p,a,a - p(2g-2))
    \times \mathcal{R}_{\Gamma}(q-p,b-a + p(2g-2)),
  \end{displaymath}
  where the notation $\mathcal{R}_\Gamma(p,q,a,b)$ indicates
  the moduli space of representations of $\Gamma$ in $\U(p,q)$ with
  invariants $(a,b)$, and $R_{\Gamma}(n,d)$ denotes the moduli space of
  degree $d$ representations of $\Gamma$ in $\U(n)$. 
  
  (A similar result holds if $p>q$, as well as if $\tau=-p(2g-2)$).
\end{theorem}

\begin{proof} Proposition~\ref{prop:R=M} and 
Theorem \ref{prop:rigidity}.
\end{proof}

{}As observed in Section~\ref{subs:invar} (cf.\ (\ref{eqn:RUpq})),  the 
spaces 
$\mathcal{R}(a) = \mathcal{R}_\Gamma(a,-a)$
can be identified with components of  $\mathcal{R}(\U(p,q))$, i.e.\  
with components of the moduli space for representations of 
$\pi_1 X$ in $\U(p,q)$. Applying Theorems
\ref{thm:results-R-Gamma} and \ref{thm:RgammaTmax},
together with the observation that $\tau(a,-a)=2a$ in the special 
case where $b=-a$, we thus obtain the following results for 
$\mathcal{R}(\U(p,q))$. Notice that the condition
$\GCD(p+q,a+b)=1$ is never satisfied if $a+b=0$.

% This explains the absence of a part (4) in 
%Corollary \ref{cor:results-R-SUpq}. 

\begin{theorem}\label{thm:results-R-SUpq1}
Let $(p,q)$ be any pair of positive integers and let 
$a\in\Z\oplus\Z$ be such that $|a|\le\min\{p,q\}(g-1)$. 

\begin{itemize}

\item[$(2)$] The moduli space $\mathcal{R}_{\Gamma}(a)$ is non-empty
  
\item[$(2)$] If either of the following sets of conditions 
apply, then the moduli space $\mathcal{R}^*(a)$ is a non-empty, 
smooth manifold of the expected dimension, with connected closure 
$\bar{\mathcal{R}}^*(a)$ in $\mathcal{R}(a)$:

\begin{enumerate}
\item[$(i)$] $0<|a|<\min\{p,q\}(g-1)$ , or
\item[$(ii)$] $|a|=p(g-1)$ and $p=q$,
\end{enumerate}

\item[$(3)$]  If any one of the 
following sets of conditions apply, then the moduli space 
$\mathcal{R}(a)$ is connected:
  
\begin{enumerate}
\item[$(i)$] $a=0$,
  
\item[$(ii)$] $|a|=\min\{p,q\}(g-1)$ and $p\neq q$,
  
\item[$(iii)$] $(p-1)(g-1) < |a| 
\leq p(g-1)$ and $p=q$,
  
\end{enumerate}
\end{itemize}
\end{theorem}

\begin{theorem}\label{thm:results-R-SUpq2}
Let $(p,q)$ be any pair of positive integers such that 
$p\ne q$. If  $|a|=\min\{p,q\}(g-1)$ then $\mathcal{R}^*(a)$ is empty and
every representation in $\mathcal{R}(a)$ is reducible. If 
$p < q$, then every such representation decomposes as a direct sum of 
a semisimple representation of 
$\Gamma$ in $\U(p,p)$ with maximal Toledo invariant and a semisimple 
representation in  $\U(q-p)$.  Thus, if 
$a=p(g-1)$ then there is an 
isomorphism 
  \begin{displaymath}
    \mathcal{R}(a) \cong \mathcal{R}_\Gamma(p,p,a,a - p(2g-2))
    \times \mathcal{R}_{\Gamma}(q-p,p(2g-2)),
  \end{displaymath}
  where the notation $\mathcal{R}_\Gamma(p,q,a,b)$ indicates
  the moduli space of representations of $\Gamma$ in $\U(p,q)$ with
  invariants $(a,b)$, and $R_{\Gamma}(n,d)$ denotes the moduli space of
  degree $d$ representations of $\Gamma$ in $\U(n)$. 
  
  (A similar result holds if $p>q$, as well as if $a=-p(g-1)$).
\end{theorem}

{}From Theorem~\ref{thm:results-R-Gamma} and 
Proposition~\ref{prop:principal-jac} we obtain the following theorem 
about the moduli spaces for $\PU(p,q)$ representations of $\pi_1 X$.  
Note that the closure $\bar{\mathcal{R}}^*[a,b]$ in 
$\mathcal{R}[a,b]$ is the image of $\bar{\mathcal{R}}_\Gamma^*(a,b)$ 
under the map of Proposition~\ref{prop:principal-jac}, hence these 
two spaces have the same number of connected components. 

\begin{theorem}
\label{thm:results-R1}

Let $(p,q)$ be any pair of positive integers and let 
$(a,b)\in\Z\oplus\Z$ be such that $0\le|\tau(a,b)|\le\tau_M$. 

\begin{itemize}

\item[$(1)$] The moduli space $\mathcal{R}[a,b]$ is non-empty.

\item[$(2)$] If either of the following sets of conditions apply,
then the moduli space $\mathcal{R}^*[a,b]$ is a non-empty smooth 
manifold of the expected dimension, with connected closure 
$\bar{\mathcal{R}}^*[a,b]$ in $\mathcal{R}[a,b]$: 

\begin{enumerate}
\item[$(i)$] $0<|\tau(a,b)|<\tau_M$ , or
\item[$(ii)$] $|\tau(a,b)|=\tau_M$ and $p=q$,
\end{enumerate}

 \item[$(3)$]  If any one of the 
following sets of conditions apply, then the moduli space 
  $\mathcal{R}[a,b]$ of all semi-simple representations is connected:

\begin{enumerate}
\item[$(i)$] $\tau(a,b)=0$,
  
\item[$(ii)$] $|\tau(a,b)|=\tau_M$ and $p\neq q$. , 
  
\item[$(iii)$] $(p-1)(2g-2) < |\tau| \leq \tau_M = p(2g-2)$ and 
$p=q$,

\item[$(iv)$]  $\GCD(p+q,a+b)=1$ 
\end{enumerate}
\end{itemize}
\end{theorem}

\begin{theorem}\label{thm:results-R2}
Let $(p,q)$ be any pair of positive integers such that $p\ne q$, and 
let $(a,b)$ be such that 
$|\tau(a,b)|=\tau_M$. Then $\mathcal{R}^*[a,b]$ is empty. If $p < q$, 
then every such representation reduces to a semisimple  
representation of $\pi_1 X$  in  $\mathrm{P}(\U(p,p) \times 
\U(q-p))$, such that the $\PU(p,p)$ representation induced via 
projection on the first factor has maximal Toledo invariant. (A 
similar result holds if $p>q$.) 
\end{theorem}

\begin{remark}
  As explained by Hitchin in \cite[Section 5]{hitchin:1987}, the
  moduli space of irreducible representations in the adjoint form of a
  Lie group is liable to acquire singularities, because of the
  existence of stable vector bundles which are fixed under the action
  of tensoring by a finite order line bundle.  For this reason we do
  not make any smoothness statements in Theorem~\ref{thm:results-R1}.
\end{remark}

%%%%%%%%%%%%%%%%%%%%%%%%%%%%%%
\subsection{Total number of components and coprimality conditions}
\label{subs:total-number}
%%%%%%%%%%%%%%%%%%%%%%%%%%%%%

We end with some elementary observations about the total number of 
components in the decomposition 
$\mathcal{R}(\PU(p,q))=\bigcup_{(a,b)}\mathcal {R}[a,b]$, and about
the number of such components for which the coprime condition 
$\GCD(p+q,a+b)=1$ apply. We begin with the number of components. 

By definition, $\tau(a,b)$ takes values in  $\frac{2}{n}\Z$, where 
$n=p+q$. 

\begin{proposition}\label{prop:[a,b]vsTau} Suppose that $\GCD(p,q)=k$.
Then the map 
\begin{align*}
\tau\colon \Z\oplus\Z/(p,q)\Z&\longrightarrow \frac{2}{n}\Z\ \\
[a,b]&\longmapsto \frac{2}{n}(aq-bp) 
\end{align*}
fits in an exact sequence 
\begin{displaymath}
0 \lto \Z/k\Z \overset{\sigma}{\lto}
\Z\oplus\Z/(p,q)\Z \overset{\tau}{\lto}
\frac{2k}{n}\Z \lto 0
\end{displaymath}
where the map $\sigma$ is $[t]\mapsto [t\frac{p}{k},t\frac{q}{k}]$. 
In particular, $\tau$ is a $k:1$ map onto the subset 
$\frac{2k}{n}\Z\subset \frac{2}{n}\Z$. 
\end{proposition}

\begin{proof} The map $\sigma$ is clearly injective, and 
$\tau\circ\sigma=0$. To see that $\ker(\tau)=\im(\sigma)$, observe 
that if $\tau[a,b]=0$ then either $a=b=0$ or 
$\frac{a}{b}=\frac{p}{q}$, i.e.\ $[a,b]=[t\frac{p}{k},t\frac{q}{k}]$
for some $t\in\Z$.  Finally, if 
$a_0q-b_0p=k$, then for any $l\in\Z$ we have $\tau[la_0,lb_0]=\frac{2k}{n}l$.
Thus $\tau$ is surjective onto $\frac{2k}{n}\Z$. 
\end{proof}

\begin{remark} Proposition \ref{prop:[a,b]vsTau} shows why\footnotemark  we 
must use 
$[a,b]$ rather than $\tau$ to label the components of $\mathcal 
{R}(\PU(p,q))$ or of $\mathcal {R}_{\Gamma}(\U(p,q))$. 
\footnotetext{Unless $p$ and $q$ are
coprime, in which case there is a bijective correspondence between 
$[a,b]$ and 
$\tau$.} 
\end{remark}

\begin{definition}\label{defn: C-N} Suppose that 
$\GCD(p,q)=k$. Define
\begin{equation}\label{eqn:C}
\mathcal{C}=\tau^{-1}([-\tau_M,\tau_M]\cap \frac{2k}{n}\Z)\ ,
\end{equation}
where $\tau$ is the map defined in Proposition \ref{prop:[a,b]vsTau}. 
\end{definition}

The following is then an immediate corollary of Proposition
\ref{prop:[a,b]vsTau}.
 
\begin{corollary}\label{cor:number-of-[a,b]}
 Suppose that $\GCD(p,q)=k$ and 
$\mathcal{C}$ is as above. Then 
$\mathcal{C}$ is precisely the 
set of all the points in $\Z\oplus\Z/(p,q)\Z$ which label components 
$\mathcal{R}[a,b]$ in $\mathcal{R}(\PU(p,q))$.  The
cardinality of $\mathcal{C}$ is 
\begin{align*}
|\mathcal {C}|&=2n\min\{p,q\}(g-1)+k\\ 
&=|([-\tau_M,\tau_M]\cap 
\frac{2}{n}\Z)|+\GCD(p,q)-1\ .
\end{align*} 
\end{corollary}

\begin{proof} The first statement is a direct consequence of 
Proposition 
\ref{prop:[a,b]vsTau} and the bound on $\tau$. Suppose for definiteness that
$\min\{p,q\}=p$. Then since 
$\tau_M=2\min\{p,q\}(g-1)=\frac{2k}{n}(n\frac{p}{k}(g-1))\in \frac{2k}{n}\Z$, 
the number of points in $[-\tau_M,\tau_M]\cap \frac{2k}{n}\Z$ is  
$2n\frac{p}{k}(g-1)+1$. The second statement now follows from the 
fact that $\tau$ is a $k:1$ map. The proof is similar if 
$\min\{p,q\}=q$.
\end{proof}

Finally, we examine the coprime condition 
$\GCD(p+q,a+b)=1$. We regard $p$ and $q$ as fixed, but allow 
$[a,b]$ to vary. The coprime 
condition $\GCD(p+q,a+b)=1$ can thus be satisfied on some components 
but not on others. 

\begin{definition}\label{prop:good-bad}Fix $p$ and $q$ and let 
$\mathcal {C}\subset \Z\oplus\Z/(p+q)\Z$ be as in 
Definition \ref{defn: C-N}. Define 
$\mathcal {C}_{1}$ to be the subset of classes 
$[a,b]\in\mathcal {C}$ for which the condition $\GCD(p+q,a+b)=1$
is satisfied. 
\end{definition}

\begin{proposition}Fix $p$ and $q$ and let $\mathcal{C}$ and 
$\mathcal {C}_{1}$ be as above.
Both $\mathcal {C}_{1}$ and its complement in 
$\mathcal {C}$ are non-empty. 
\end{proposition} 

\begin{proof}
If $a=p$ and $b=q-1$ then $\GCD(p+q,a+b)=1$. Also, 
$\tau(p,q-1)=\frac{2p}{p+q}$, which is in 
$[-\tau_M,\tau_M]\cap \frac{2k}{n}\Z$. Thus
$[p,q-1]$ is in $\mathcal {C}_{1}$. It is similarly straightforward
to see that $(a,b)= (0,0)$ defines an element in $\mathcal {C}- 
\mathcal {C}_{1}$, as does $(a,b)=(p,-p)$ if $p\leq q$
or $(a,b)= (q,-q)$ if $q\leq p$. 
\end{proof}
 
It seems somewhat complicated to go beyond this result and completely 
enumerate the elements in $\mathcal {C}_{1}$. The following result 
is, however, a step in that direction. 

\begin{definition} Let $\Omega\subset\R\oplus\R$ be the region depicted in
Figure 2, i.e.\ the region bounded by (i) the ray $ b=q$  and $a\le 
p$, (ii) the ray $ a=p$ and $b\le q$, (iii) the ray $ a=0$  and $b\le 
0$, (iv) the ray $ b=0$ and $a\le 0$, (v) the line 
$ aq-bp=\frac{n}{2}\tau_M$, and (vi) the line $ aq-bp=-\frac{n}{2}\tau_M$,
and including all the boundary lines except the first two rays. Let 
$\Omega_{\Z}$ be the set of integer points in $\Omega$, i.e.\ 
$\Omega_{\Z}=\Omega\bigcap\Z\oplus\Z$. We refer to $\Omega$  as the 
\emph{fundamental region for $(p,q)$} (see Figure 2). Then 
$\Omega_{\Z}$ is the integer lattice inside the fundamental region.
\end{definition}

\begin{proposition} Suppose that $p$ and $q$ are integers with
$\GCD (p,q)=k$.
\begin{enumerate}
\item[$(1)$] There is a bijection between $\mathcal{C}$ and
$\Omega_{\Z}$.
\item[$(2)$] If $(a,b)$ lies in $\Omega_{\Z}$ then $d=a+b$ satisfies the
bounds
\begin{equation}
-n(g-1)\le d< n\ .
\end{equation}
All values of $d$ in this range occur.
\item[$(3)$] Let $l_t$ denote the line $aq-bp=tk$. Then the points
on $l_t\bigcap\Omega_{\Z}$ define the locus of points $(a,b)$ for
which $\tau(a,b)=t\frac{2k}{n}$.
\item[$(4)$] The line $l_t$ intersects $\Omega_{\Z}$ for
$-\frac{n}{2k}\tau_M\le t
\le \frac{n}{2k}\tau_M$. For each integer $t$ in this range, there are
$k$ points on $l_t\bigcap\Omega_{\Z}$.
\item[$(5)$] For a fixed $t$, $\GCD(a+b,\frac{n}{k})$ is the same
for all integer points
$(a,b)$ on $l_t\bigcap\Omega_{\Z}$.
\item[$(6)$] Fix $t$ and let $(a,b)$ be any point in the set
$l_t\bigcap\Omega_{\Z}$.  If  $\GCD(a+b,\frac{n}{k})\ne 1$
then $\GCD(a'+b',n)\ne 1$ for all $(a',b')\in l_t\bigcap\Omega_{\Z}$.
%\item The condition for d to occur is that there is an allowed $\tau$ and
%an allowed $b$ such that
%$pb$ is congruent to $-n\tau/2$ mod $q$ YIKES!
\end{enumerate}
\end{proposition}

\begin{proof}
(1) Suppose first that $\frac{a}{p}\ge \frac{b}{q}$. Pick $l$ such
that $0\le a+lp\le p$. Then $b+lq\le q$, so that $(a+lp, b+lq)$ is in
the fundamental region. Similarly, if  $\frac{a}{p}\le \frac{b}{q}$
then we pick $l$ such that $0\le b+lq\le q$ and see that $a+lp\le p$.
In this way we get a well defined map from $\mathcal {C}$ to the
fundamental region. The map is clearly injective. To see that it is
surjective, notice that the boundary lines $
aq-bp=\pm\frac{n}{2}\tau_M$ correspond to the conditions
$\tau=\pm \tau_M$.

(2) This is clear from a sketch of the fundamental region (see Figure
2), in which the loci of points with constant value of
$d=a+b$ are straight lines of slope $-1$. The maximal value for $d$
corresponds to the line passing through the top right corner of the
region, i.e. through $(p,q)$. Thus $d_{max}=p+q=n$. The minimal value
for $d$ corresponds either to the line passing through
$(-\frac{n}{2q}\tau_M,0)$ or to the line through $(0,
-\frac{n}{2p}\tau_M)$, depending on which yields the smaller value
for $d$. Since $\tau_M=2\min\{p,q\}(g-1)$, we find in all cases that
$d_{min}= -n(g-1)$. It is straightforward to see that all intermediate
values for $d$ occur.

(3)-(4) This is simply a restatement of Proposition
\ref{prop:[a,b]vsTau}.

(5)-(6) Both follow from the fact that for any two points $(a,b)$ and
$(a',b')$ on $l_t$, we get $d'=d+s\frac{n}{k}$ for some $s\in\Z$.
\end{proof}

\begin{remark} Part (6) says that for fixed $t$, if
$\GCD(a+b,\frac{n}{k})\ne
1$ for {\it any} point $(a,b)\in l_t\bigcap\Omega_{\Z}$, then
$\GCD(a+b,n)\ne 1$ for {\it all} points $(a,b)\in l_t\bigcap\Omega_{\Z}$.
That is, we can detect the non-coprimality of $(a+b, n)$ for all
$(a,b)\in l_t\bigcap\Omega_{\Z}$ by checking the non-coprimality of
$(a+b, \frac{n}{k})$ at any one $(a,b)\in l_t\bigcap\Omega_{\Z}$. We
cannot however check for coprimality in the same way.  If
$\GCD(a+b,\frac{n}{k})= 1$, it is possible that $\GCD(a'+b',n)\ne 1$
for some $(a',b')\in l_t\bigcap\Omega_{\Z}$. For example, take $p=2,
q=4, a=-1, b=0, a'=0, b'=2 $, and $t=-2$. Then $\GCD(a'+b',n)=2$
while $\GCD(a+b,\frac{n}{k})= 1$.
\end{remark}

\newpage

\begin{picture}(400,550)(-20,-350)
\linethickness{1pt}
%%%%%%%%%%%%%%%%%%
%p<q label
%%%%%%%%%%%%%%%%%
\put (310, 160){$p\ne q$}
%\tau-axis
\put (0,0){\vector(0,1){200}}
\put (-10,195){$\tau$}
\put (-3,0){\line(1,0){6}}
\put (-36,-2){$-\tau_M$}
\put (-3,80){\line(1,0){6}}
\put (-20,78){$0$}
\put (-3,160){\line(1,0){6}}
\put (-32,159){$\tau_M$}
%\alpha-axis
\put (0,-20){\vector(1,0){370}}
\put (360,-30){$\alpha$}
\put (0,-23){\line(0,1){6}}
\put (-2,-35){$0$}
\put (10,-23){\line(0,1){6}}
\put (0,-40){$\frac{|p-q|}{2\max\{p,q\}}(2g-2)$}
\put (130,-23){\line(0,1){6}}
\put (113,-35){$2g-2$}
%\alpha_m
\put (10,0){\line(6, 4){128}}
\put (52,40){$\alpha_m$}
\put (10,160){\line(6, -4){128}}
%2g-2
\multiput (130,0)(0,5){32}{\line(0,1){2}}
%\alpha_M
\put (130,0){\line(6, 2){250}}
\put (270,40){$\alpha_M$}
\put (130,160){\line(6, -2){250}}
%boundaries
\put (10,0){\line(1,0){120}}
\put (10,-4){\line(0,1){8}}
%\put (9,-4){\line(0,1){8}}
\put (130,0){\circle*{7}}
\put (10,160){\line(1,0){119}}
\put (10,156){\line(0,1){8}}
\put (130,160){\circle*{7}}
%sample line at \tau
\put (70,40){\line(1,0){180}}
\put (70,36){\line(0,1){8}}
\put (130,40){\circle*{7}}
\put (251,36){\line(0,1){8}}
%sample line at 0
\put (130,80){\line(1,0){240}}
\put (373,76){\line(0,1){8}}
\put (130,80){\circle*{7}}
%%%%%%%%%%%%%%%%%%%%%
%p=q label
%%%%%%%%%%%%%%%%%%%%
\put (310, -110){$p=q$}
%\tau-axis
\put (0,-300){\vector(0,1){200}}
\put (-10,-105){$\tau$}
%\put (-3,-300){\line(1,0){6}}
\put (-36,-301){$-\tau_M$}
\put (-3,-220){\line(1,0){6}}
\put (-20,-222){$0$}
%\put (-3,-140){\line(1,0){6}}
\put (-33,-142){$\tau_M$}
%\alpha-axis
\put (0,-320){\vector(1,0){370}}
\put (360,-330){$\alpha$}
\put (0,-323){\line(0,1){6}}
\put (-2,-335){$0$}
\put (120,-323){\line(0,1){6}}
\put (103,-335){$2g-2$}
%\alpha_m
\put (0,-300){\line(6, 4){128}}
\put (52,-250){$\alpha_m$}
\put (0,-140){\line(6, -4){128}}
%d_1-d_2
\put (0,-300){\line(6, 2){245}}
\put (170,-250){$d_1-d_2$}
\put (0,-140){\line(6, -2){245}}
%2g-2
\multiput (120,-300)(0,5){32}{\line(0,1){2}}
%boundaries
\put (0,-300){\vector(1,0){350}}
\put (120,-300){\circle*{7}}
\put (0,-140){\vector(1,0){350}}
\put (120,-140){\circle*{7}}
%sample line at \tau
%\put (30,-280){\vector(1,0){320}}
%\put (30,-284){\line(0,1){8}}
%\put (120,-280){\circle*{7}}
%sample line at 0
\put (120,-220){\vector(1,0){230}}
\put (120,-220){\circle*{7}}

\end{picture}

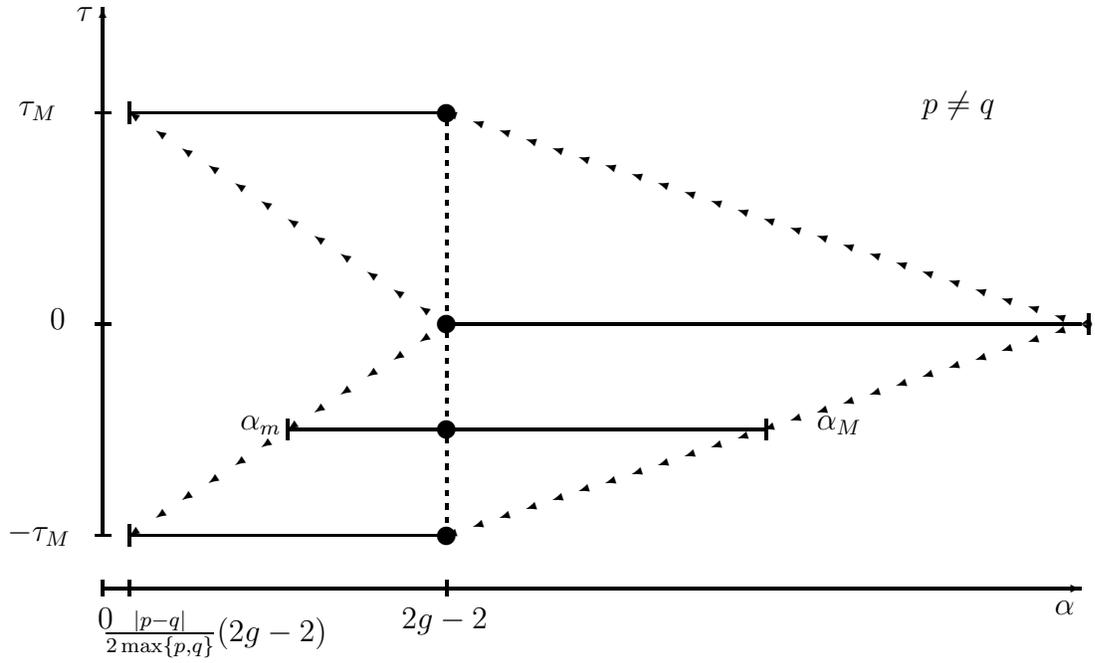
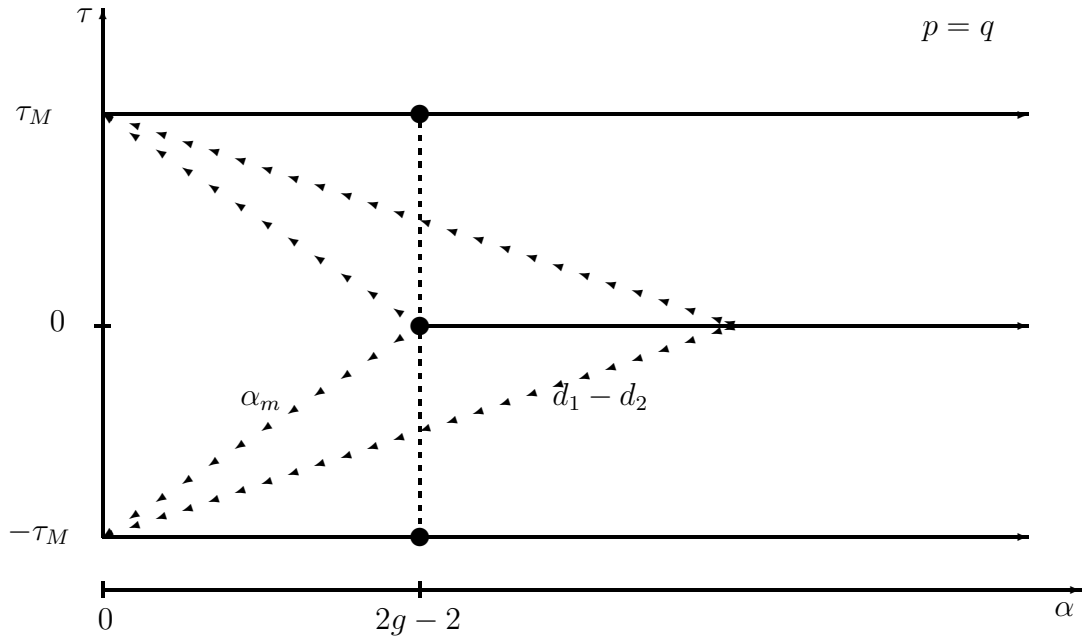
\noindent Figure 1: Range for the stability parameter $\alpha$ for
triples in $\mathcal{N}(a,b)$, displayed as functions of 
$\tau=\frac{2pq}{p+q}(\frac{a}{p}-\frac{b}{q})$, and showing the 
relative location of $2g-2$. 
\newpage

%%%%%%%%%%%
%%Figure 2
%%%%%%%%%%
\begin{picture}(400,400)(100,0)
%a-axis (with labels)
\multiput (270,220)(5,0){22}{\line(1,0){1}}
\put (380,220){\vector(1,0){2}}
\put (375,208){$a$}
\put (210,217){\line(0,1){6}}
\put (200,206){$-\frac{p+q}{2q}\tau_M$}
\put (268,208){$0$}
\put (312,210){$p$}
%b-axis (with labels)
\multiput (270,220)(0,5){33}{\line(0,1){1}}
\put (270,380){\vector(0,1){3}}
\put (260,375){$b$}
\put (267,70){\line(1,0){6}}
\put (226,70){$-\frac{p+q}{2p}\tau_M$}
\put (263,216){$0$}
\put (263,322){$q$}
%%%
\linethickness{1pt}
% b=q and a\le p
%\put (310,320){\line(-1,0){60}}
\multiput (308,320)(-3,0){20}{\line(1,0){2}}
% a=p and b\le q
\multiput (310,318)(0,-3){50}{\line(0,1){2}}
% b=0 and a\le 0
\put (270,220){\line(-1,0){60}}
% a=0 and b\le 0
\put (270,220){\line(0,-1){150}}
%lines aq-bp=\tau
\put (270,70){\line(2,5){40}}
\put (273,270){$\tau=0$}
\put (270,220){\line(2,5){40}}
\put (273,110){$\tau=\tau_M$}
\put (210,220){\line(2,5){40}}
\put (210,290){$\tau=-\tau_M$}
%lines a+b=d
\put (315,275){$d=q$}
\multiput (270,320)(5,-5){9}{\circle*{3}}
\put (315,175){$d=0$}
\multiput (310,180)(-5,5){17}{\circle*{3}}
\put (305,130){$d=-\frac{p+q}{2q}\tau_M$}
\multiput (270,160)(5,-5){6}{\circle*{3}}
%\multiput (265,165)(-5,5){12}{\circle*{1}}
%\put (100,200){$d=-\frac{1}{p}\tau_M$}
%\multiput (270,70)(-5,5){25}{\circle*{1}}
%\multiput (350,280)(-5,5){25}{\circle*{1}}
\end{picture}

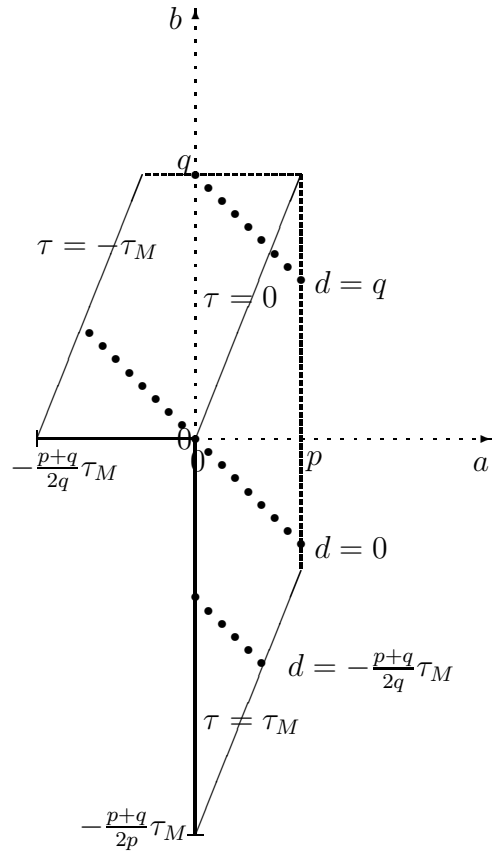
\noindent Figure 2: Fundamental region for $(a,b)$. Components of 
$\mathcal{R}(\PU(p,q))$ correspond to the integer points in this region. 
Illustrative lines of constant $\tau$ (at $\tau=-\tau_M, 0, \tau_M$) 
and lines of constant $d$ (at $d=-\frac{n}{2q}\tau_M, 0, q$) are 
shown.

\newpage

%%%%%%%%%%%%%%%%%%%%%%
%\section{References}
%%%%%%%%%%%%%%%%%%%%%%

%\bibliography{upq}
%\bibliographystyle{invmatplain}

\end{document}